\documentclass[american,oneside]{amsart}
\usepackage[T1]{fontenc}
\usepackage{geometry}
\usepackage{graphicx}
\usepackage{amssymb}

\makeatletter

\providecommand{\tabularnewline}{\\}

 \theoremstyle{plain}
\newtheorem{thm}{Theorem}[section]
\newenvironment{lyxlist}[1]
{\begin{list}{}
{\settowidth{\labelwidth}{#1}
 \setlength{\leftmargin}{\labelwidth}
 \addtolength{\leftmargin}{\labelsep}
 }}
{\end{list}}
  \theoremstyle{definition}
  \newtheorem{defn}[thm]{Definition}
  \theoremstyle{remark}
  \newtheorem{rem}[thm]{Remark}
  \theoremstyle{plain}
  \newtheorem{prop}[thm]{Proposition}
  \theoremstyle{plain}
  \newtheorem{lem}[thm]{Lemma}
  \theoremstyle{plain}
  \newtheorem{cor}[thm]{Corollary}
 \theoremstyle{definition}
  \newtheorem{example}[thm]{Example}

\usepackage{subfigure}
\usepackage{graphicx}
\numberwithin{equation}{section}
\numberwithin{figure}{section}

\usepackage{babel}
\makeatother

\begin{document}
\newcommand{\ww}[1]{\mathbb{#1}}

\newcommand{\germ}[1]{\ww{C}\left\{  #1\right\}  }

\newcommand{\pp}[1]{\frac{\partial}{\partial#1}}

\newcommand{\flw}[3]{\Phi_{#1}^{#2}#3}

\newcommand{\fol}[1]{\mathcal{F}_{#1}}

\newcommand{\eps}{\varepsilon}

\newcommand{\sect}[2]{V_{#1}^{#2}}

\newcommand{\ssect}[2]{\mathcal{V}_{#1}^{#2}}

\newcommand{\ddt}{\cdot}

\newcommand{\void}{\emptyset}

\newcommand{\ov}{\overline}

\title{Analytical moduli for unfoldings of saddle-node vector fields}

\author{Christiane ROUSSEAU$\,^{\star}$, Lo\"{i}c TEYSSIER$\,^{\dagger}$}

\date{March $\mbox{27}^{\mbox{th}}$ 2007, revised July $\mbox{25}^{\mbox{th}}$
2007}

\thanks{This work is supported by NSERC in Canada.}



\email{\textbf{Christiane Rousseau : }\texttt{\footnotesize rousseac@dms.umontreal.ca,
}\textbf{Loïc Teyssier :} \texttt{\footnotesize teyssier@math.u-strasbg.fr}}

\urladdr{\texttt{\footnotesize http://www.dms.umontreal.ca/\textasciitilde{}rousseac,
http://math.u-strasbg.fr/\textasciitilde{}teyssier}}

\address{\textbf{Christiane Rousseau :} Département de mathématiques et de
statistique, Université de Montréal CP 6128, succ. Centre-ville Montréal,
Québec H3C 3J7, Canada\\
 \textbf{Loïc Teyssier :} Laboratoire IRMA, 7 rue René Descartes F-67084
Strasbourg Cedex, France}

\maketitle
\begin{tabular}{crc}
{\footnotesize $\,^{\star}$Département de Mathématiques et de Statistique} & {\footnotesize ~\hfill{}~} & {\footnotesize $\,^{\dagger}$Institut de Recherche Mathématique Avancée}\tabularnewline
{\footnotesize Université de Montréal, Canada} & {\footnotesize ~\hfill{}~} & {\footnotesize Université Louis Pasteur, Strasbourg, France}\tabularnewline
\end{tabular}

\begin{abstract}
In this paper we consider germs of $k$-parameter generic families
of analytic 2-dimensional vector fields unfolding a saddle-node of
codimension $k$ and we give a complete modulus of analytic classification
under orbital equivalence and a complete modulus of analytic classification
under conjugacy. The modulus is an unfolding of the corresponding
modulus for the germ of a vector field with a saddle-node. The point
of view is to compare the family with a {}``model family'' \emph{via}
an equivalence (conjugacy) over canonical sectors. This is done by
studying the asymptotic homology of the leaves and its consequences
for solutions of the cohomological equation. 
\end{abstract}
\begin{center}
\emph{This paper is dedicated to the memory of Adrien Douady.} 
\par\end{center}

\bigskip{}

\bigskip{}

\section{\label{sec:Introduction}Introduction}

We consider germs of generic unfoldings of holomorphic vector fields
$Z_{0}$ in $\ww{C}^{2}$ near an isolated singularity which is a
saddle-node of codimension $k\in\ww{N}_{>0}$ (\emph{i.e.} of multiplicity
$k+1$). For such germs there exist polynomial normal forms under
orbital equivalence (\emph{resp.} conjugacy) but generically there
exist no analytic change of coordinates to these normal forms: if
we restrict to real variables in the case of real vector fields the
change of coordinates is $C^{\infty}$ in the case of a single vector
field and only $C^{N}$ for arbitrarily high $N$ in the case of an
unfolding.

A modulus space has been given for a single vector field by Martinet-Ramis
\cite{MR1} for the problem of orbital equivalence and by Teyssier
\cite{T2} and Meshcheryakova-Voronin \cite{MV} for the problem of
conjugacy (\cite{MV} treats the codimension 1 case). In both cases
the modulus is functional and the modulus space is huge. In this paper
we address the same problem for germs of families unfolding a germ
of vector field with a saddle-node at the origin. We could complete
the first part of the program. We prove a theorem allowing to prepare
a family and we identify two complete moduli of analytic classification
for prepared families: one under orbital equivalence and one under
conjugacy. These moduli are unfoldings of the corresponding moduli
for the associated germs of vector fields with a saddle node obtained
by Martinet-Ramis in the orbital case and Teyssier and Meshcheryakova-Voronin
for the conjugacy case. In each case the identification of the modulus
space is still an open problem. Our approach enlightens why the modulus
spaces for the case of a single vector field are so large. Indeed
a saddle-node of codimension $k$ is the confluence of $k+1$ simple
singular points. Each singular point is an organizing locus for the
space of leaves in its neighborhood. The space of leaves restricted
to special domains (canonical sectors) have a rigid complex structure:
they are parameterized by $\mathbb{C}$ with one special leaf, the
{}``center leaf'' parameterized by $0$. Hence the only changes
of parameterization of the space of leaves are the linear maps. In
the global family these local spaces of leaves generically glue in
a non trivial way. When this persists until the limit case where all
$k+1$ singular points merge together this yields divergence of the
normalizing change of coordinates for a single vector field with a
saddle-node.

The polynomial normal form for the family is what we can call the
{}``model family''. We can bring the family into this form using
a formal transformation near $\left(0,0,0\right)\in\ww{C}^{k+2}$.
In the model family all spaces of leaves glue trivially, so the model
family is too poor to encode all the rich dynamics of an arbitrary
analytic family of vector fields. Hence there exists in general no
analytic family of changes of coordinates (and time scalings in the
case of orbital equivalence) to the model family. However there exist
analytic families of changes of coordinates (and time scalings in
the case of orbital equivalence) to the model family over canonical
sectors. The modulus measures the obstruction to gluing the different
changes of coordinates into a global change of coordinates.

\medskip{}

There are at least two different approaches to the modulus of a single
vector field with a saddle-node at the origin. The first approach
by Martinet-Ramis \cite{MR1} characterizes the vector field under
orbital equivalence by identifying the divergence of the normalizing
formal power series with a co-chain in the ring of summable power
series, which in turn can be understood geometrically as a collection
of transition diffeomorphisms between consecutive sectorial spaces
of leaves. It turns out that these invariants coincide with the Écalle-Voronin
invariants of the induced holonomy of the strong separatrix. Meshcheryakova-Voronin
added the first-return time needed to compute the holonomy to identify
classes under conjugacy for vector fields. The second approach, by
Teyssier \cite{T2} uses the geometry of the leaves in the neighborhood
of the saddle-node, which is described in terms of asymptotic homology.
Both approaches could have been generalized (unfolded) to the family
case. A treatment with the first approach would have been similar
to \cite{R2} and \cite{RC}. We have chosen to use the second approach
so as to enlighten the asymptotic homology of the leaves and the special
geometry of the space of leaves. Solving the conjugacy problem is
then equivalent to solving some cohomological equations.

\medskip{}

For convenience we will locate the singularity of $Z_{0}$ at $\left(0,0\right)$.
An \textbf{unfolding} $(Z_{\eps})_{\varepsilon}$ of $Z_{0}$ is a
germ of analytic family of analytic vector fields. It has a representative
for $\left|\left|\varepsilon\right|\right|<\rho$ and $(x,y)\in r\ww{D}\times r'\ww{D}$,
where $\ww{D}:=\left\{ \left|z\right|<1\right\} \subset\ww{C}$. We
want to study the space of all such families or, more precisely, its
quotient under the action of local changes of coordinates (and time
scalings in the case of orbital classification).

\medskip{}

The strategy is the following. We first {}``prepare'' the family
to a preliminary prenormal form and we identify for each family the
{}``model family'' to which it will be compared. In particular we
show that in this prenormal form the parameters are analytic invariants
and hence that any equivalence or conjugacy preserves the parameters.
This allows to work for each fixed value of the parameter (but on
a neighborhood of the singular point independent of the chosen parameter).
We then determine canonical sectors over which the space of leaves
has a canonical structure. Over each canonical sector we get an equivalence
between the original family and the model family. An equivalence between
any two families over a canonical sector is obtained by composing
the equivalence of the first family to the model with the equivalence
of the model to the second family. The modulus is the obstruction
to gluing the equivalences to the model family over the canonical
sectors into a global equivalence. If two families have the same modulus
it is then possible to glue together the equivalences over canonical
sectors into a global equivalence between the two families.

The program above requires first to study in detail the model family.
This is started in Section~\ref{sec:preparation} and finished in
Section~\ref{sec:model-family}. These two parts are quite long,
but are likely to be used in further work on the realization part.
Because this preliminary part is long we have added a Section~\ref{sec:results}
with the statements of the results. In Section~\ref{sec:center}
we show how the $k$ sectorial center manifolds of a saddle-node of
codimension $k$ unfold as $k$ special leaves over $k$ canonical
sectors. In the model family these special leaves glue together as
a global leaf; measuring the obstruction to a global gluing is the
first part of the orbital modulus. In Section~\ref{sec: asymptotic}
we introduce the notion of asymptotic homology and we build the canonical
sectors. In Section~\ref{sec:cohomological} we discuss solutions
of the cohomological equation, as these will be the tool for the classification
problem. In Section~\ref{orb_classification} we give a new proof
of the Hukuhara-Kimura-Matuda sectorial normalization theorem together
with a generalization to unfoldings restricted to canonical sectors.
Sections~\ref{sec:orbital}, \ref{sec:conjugacy} and \ref{Proofs}
contain the full definitions of the modulus of an analytic family
under orbital equivalence and under conjugacy and the proof that the
modulus is indeed a complete modulus of analytic classification. Section~\ref{sec:Questions}
contains questions for future research and applications.

\medskip{}

We were precisely in the final stage of writing this paper when we
learned the death of Adrien Douady. Clearly his heritage in the subject
is immense. Although many people has conjectured the Stokes phenomena
coming from $k$-summability to be the limits of transitions when
all singular points of an unfolding were in the Poincaré domain, no
one knew how to deal with the Siegel direction. It is the visionary
geometric ideas of Douady and the thesis of his student Lavaurs which
opened the subject and the hope to derive complete invariants of analytic
classification for germs of families of vector fields. We dedicate
this paper to his memory.

\newpage{} \tableofcontents{}

\newpage{}

\section*{Index of notations}

\label{notations}

\begin{lyxlist}{00.00.0000}
\item [{{{{$k$:}}}}] a positive integer. 
\item [{{{{{\textbf{$\germ{x_{1},\ldots,x_{n}}$}:}}}}}] ~\\
 the algebra of germs of holomorphic functions on $\mathbb{C}^{n}$
at $0\in\ww{C}^{n}$. 
\item [{{{{{$X\cdot F$:}}}}}] the Lie derivative of the function
$F$ along the vector field $X$. 
\item [{{{{{$\varepsilon=\left(\varepsilon_{0},\ldots,\varepsilon_{k-1}\right)\in\ww{C}^{k}$:}}}}}] ~\\
 the \textbf{canonical multi-parameter} of a prepared unfolding, see
Definition~\ref{def:canonic_param}. 
\item [{{{{{$\left(P_{\varepsilon}\right)_{\varepsilon}$:}}}}}] the
analytical family of polynomials of degree $k+1$ unfolding $x^{k+1}$,
namely $P_{\eps}(x)=x^{k+1}+\eps_{k-1}x^{k-1}+\dots\eps_{1}x+\eps_{0}$. 
\item [{{{{{$\Sigma_{0}$:}}}}}] the semi-algebraic open set in
$\varepsilon$-space defined by the condition that $P_{\varepsilon}$
has $k+1$ distinct roots, see Section~\ref{sub:generic_eps}. 
\item [{{{{{$X_{\varepsilon}^{M}=P_{\varepsilon}\pp{x}+y\left(1+a\left(\varepsilon\right)x^{k}\right)\pp{y}$:}}}}}] ~\\
 the orbital model family; $a\in\germ{\varepsilon}$ is fixed once
and for all. The singular set of the prepared vector field coincides
with $P_{\varepsilon}^{-1}\left(0\right)\times\left\{ 0\right\} $.
See Section~\ref{sec:model-family}. 
\item [{{{{{$Z_{\varepsilon}^{M}=Q_{\varepsilon}X_{\varepsilon}^{M}$:}}}}}] ~\\
 the \textbf{model family} ; $Q_{\varepsilon}=C_{0,\eps}+C_{1,\eps}x+\dots+C_{k,\eps}x^{k}$
with $C_{0,\eps}\neq0$ and $\varepsilon\mapsto C_{j,\varepsilon}\in\germ{\varepsilon}$.
See Section~\ref{sec:model-family} and Theorem~\ref{pro:eqv_eps}. 
\item [{{{{{$\tau_{\varepsilon}=\frac{dx}{P_{\varepsilon}}$:}}}}}] the
canonical closed \textbf{time-form} associated to $X_{\varepsilon}^{M}$,
that is $\tau_{\varepsilon}\left(X_{\varepsilon}^{M}\right)=1$. 
\item [{{{{{$\left(X_{\varepsilon}\right)_{\varepsilon}$:}}}}}] a
prepared unfolding with $\tau_{\varepsilon}\left(X_{\varepsilon}\right)=1$.
This means the family only unfolds the foliation defined by $X_{0}$.
We write (see Proposition~\ref{pro:straighten} and Definition~\ref{def:prepared})
\begin{eqnarray*}
X_{\varepsilon}\left(x,y\right) & = & X_{\varepsilon}^{M}\left(x,y\right)+\left[P_{\varepsilon}\left(x\right)R_{0,\varepsilon}\left(x\right)+y^{2}R_{2,\varepsilon}\left(x,y\right)\right]\pp{y}\,.\end{eqnarray*}

\item [{{{{{$\left(Z_{\varepsilon}\right)_{\varepsilon}=\left(U_{\varepsilon}X_{\varepsilon}\right)_{\varepsilon}$:}}}}}] ~\\
 a \textbf{prepared unfolding} with $\left(U_{\varepsilon}\right)_{\varepsilon}\in\germ{x,y,\varepsilon}$
and $U_{\varepsilon}=Q_{\varepsilon}+O\left(P_{\varepsilon}\left(x\right)\right)+O\left(y\right)$
where $GCD\left(Q_{\varepsilon},P_{\varepsilon}\right)=1$. The function
$U_{\varepsilon}$ is called the \textbf{time part} of $Z_{\varepsilon}$,
whereas $X_{\varepsilon}$ is the \textbf{orbital part}. The modulus
of the orbital part is analyzed on $R_{\eps}(x,y)=P_{\varepsilon}(x)R_{0,\varepsilon}(x)+y^{2}R_{2,\varepsilon}(x,y)$. 
\item [{{{{{$r$,~$r'$,~$\rho$:}}}}}] the radii of the open
domain $r\ww{D}\times r'\ww{D}\times\left\{ \left|\left|\varepsilon\right|\right|\leq\rho\right\} $
considered in $\left(x,y,\varepsilon\right)$-space. Here $\left|\left|\varepsilon\right|\right|:=\max\left(\left|\varepsilon_{0}\right|^{1/(k+1)},\ldots,\left|\varepsilon_{k-1}\right|^{1/2}\right)$
and $\ww{D}=\left\{ \omega\in\ww{C}\,:\,\left|\omega\right|<1\right\} $.
What we mean by $\{||\eps||\leq\rho\}$ is $\{||\eps||<\rho'\}$ for
some $\rho'>\rho$. 
\item [{$j$:}] an element of $\ww{Z}/k$.
\item [{{{{{$\sect{j,\varepsilon}{\#}$:}}}}}] a \textbf{squid-sector}
in the $x$-variable. Here $\#$ may be $+$,$-$, $s$ or $n$. See
Definition~\ref{def:squid-sectors} and Lemma~\ref{lem:intersect}. 
\item [{{{{{$p_{j,n}$,~$p_{j,s}$~(or}}}}}] $p_{j,n}^{\pm}$,~$p_{j,s}^{\pm}$):
\\
 the singular points of $Z_{\varepsilon}$ over the closure of a sector
$\sect{j,\varepsilon}{\#}$. Here {}``$n$'' and {}``$s$'' stand
for {}``node type'' and {}``saddle type'' in the generic case
$\varepsilon\in\Sigma_{0}$. See Definition~\ref{def:pt_sing}. 
\item [{{{{{$\sigma$:}}}}}] the one-to-one correspondence associating
to a sector $V_{j,\eps}^{+}$ a sector $V_{\sigma(j),\eps}^{-}$,
where $V_{j,\eps}^{+}$ and $V_{\sigma(j),\eps}^{-}$ share the same
singular points of saddle and node type. (See Lemma~\ref{lemma_DS}
and \eqref{def_sigma}.) 
\item [{{{{{$V_{j,\sigma(j),\varepsilon}^{g}$~or~$V_{j,\varepsilon}^{g}$:}}}}}] ~\\
 the \textbf{gate sector} which is the intersection of two non consecutive
squid sectors $V_{j,\eps}^{+}$ and $V_{\sigma(j),\eps}^{-}$ sharing
the same singular points $p_{j,n}^{+}=p_{\sigma\left(j\right),n}^{-}$
and $p_{j,s}^{+}=p_{\sigma\left(j\right),s}^{-}$. 
\item [{{{{{$\ssect{j,\varepsilon}{\#}$:}}}}}] the \textbf{canonical
sector} of the foliation corresponding to $\sect{j,\varepsilon}{\#}$,
obtained by considering all points $\left(\overline{x},\overline{y}\right)\in\sect{j,\varepsilon}{\#}\times r'\ww{D}$
which can be linked to the singular point $p_{j,n}$ by a tangent
asymptotic path. See Theorem~\ref{thm:secto-trivial} and Definition~\ref{def:canonical-sector}. 
\item [{{{{{$H_{j,\varepsilon}^{\pm}$:}}}}}] the corresponding
\textbf{canonical first integral} over $\ssect{j,\varepsilon}{\pm}$
whose level sets coincide with the leaves of the foliation induced
by $Z_{\varepsilon}$ over $\mathcal{V}_{j,\varepsilon}^{\pm}$, see
Definition~\ref{def:secto-first-integ}. 
\item [{{{{{$\gamma_{j,\varepsilon}^{s}\left(p\right)$:}}}}}] an
\textbf{asymptotic path} passing through $p\in\ssect{j,\varepsilon}{s}$
linking $p_{j,n}$ and $p_{j+1,n}$. The upcoming analytic invariants
of the family will be obtained as integrals over these asymptotic
paths. See Definition~\ref{def:asy_path}. 
\item [{{{{{$\{W_{i}\}_{1\leq i\leq d}$:}}}}}] an open finite
covering of $\Sigma_{0}$ with good sectors. See Definition~\ref{good_covering}.
\item [{$\mathcal{O}_{b}\left(W\right)$:}] the set of functions $f\,:\,\ov{W}\rightarrow\ww{C}\backslash\left\{ 0\right\} $
which are continuous on the closure $\ov{W}$ of a good sector $W\subset\Sigma_{0}$
and analytic on $W$.
\item [{$\mathcal{N}_{\varepsilon}^{i}=\left(a,\psi_{0,\varepsilon}^{\infty,i},\ldots,\psi_{k-1,\varepsilon}^{\infty,i},\phi_{0,\varepsilon}^{0,i},\ldots,\phi_{k-1,\varepsilon}^{0,i}\right)$:}] ~\\
 it is defined for $\eps\in W_{i}$. The $d$-uple $\{\mathcal{N}_{\varepsilon}^{i}\}_{1\leq i\leq d}$
forms the \textbf{orbital part of the modulus} associated to $\left(X_{\varepsilon}\right)_{\varepsilon}$,
which provides a complete set of invariants of the unfolding under
orbital equivalence. The $\psi_{j,\varepsilon}^{\infty,i}$ are affine
maps and $\phi_{j,\varepsilon}^{0,i}\in\germ{h}$ with $\phi_{j,\varepsilon}^{0,i}\left(0\right)=0$.
They correspond to changes of coordinate in the space of leaves over
the intersections $\ssect{j,\varepsilon}{n}$ and $\ssect{j,\varepsilon}{s}$
respectively. See Section~\ref{sec:orbital}. 
\item [{{{{{$\mathcal{T}_{\varepsilon}^{i}=\left(C_{0,\varepsilon},\ldots,C_{k,\varepsilon},\zeta_{0,\varepsilon}^{i},\ldots,\zeta_{k-1}^{i},\varepsilon\right)$:}}}}}] ~\\
 it is defined for $\eps\in W_{i}$. The collection $\{\mathcal{T}_{\varepsilon}^{i}\}_{1\leq i\leq d}$
forms the \textbf{time part of the modulus} associated to $\left(Z_{\varepsilon}\right)_{\varepsilon}$.
Together with $a\left(\varepsilon\right)$ and $\{\mathcal{N}_{\varepsilon}^{i}\}_{1\leq i\leq d}$,
it provides a complete set of invariants of the unfolding under conjugacy.
The $C_{j,\varepsilon}$ are simply the coefficients of the polynomial
$Q_{\varepsilon}$ whereas the $\zeta_{j,\varepsilon}^{i}\in\germ{h}$
represent time scalings over $\ssect{j,\varepsilon}{s}$. See Section~\ref{sec:conjugacy}. 
\end{lyxlist}
\newpage{}

\section{\label{sec:results} Statement of results}

This section is informal. For more precise statements and definitions
we refer to the corresponding sections indicated between parentheses.

Our main goal is to provide invariants for classification of germs
of an analytic family $\left(Z_{\eps}\right)_{\varepsilon}$ under
both orbital equivalence and conjugacy. Let us define these terms.

\begin{defn}
(see Section~\ref{Proofs}) 
\begin{enumerate}
\item Two analytic vector fields (\emph{resp.} germs of analytic vector
fields) $X$ and $Y$ are \textbf{conjugate} if there exists an analytic
diffeomorphism (\emph{resp.} a germ of analytic diffeomorphism) $\Psi$
such that $\Psi^{*}X=Y$, that is $X\circ\Psi=D\Psi\left(Y\right)$. 
\item $X$ and $Y$ are \textbf{orbitally equivalent} under $\Psi$ if there
exists an analytic non-vanishing function (\emph{resp.} germ) $U$
such that $X$ and $UY$ are conjugate under $\Psi$. Equivalently
this means that the image by $\Psi$ of any integral curve of $X$
is an integral curve of $Y$. We also speak of \textbf{equivalence}
of the underlying foliations. 
\item Two analytic families (\emph{resp.} germs of analytic families of
vector fields) $\left(Z_{\varepsilon}\right)_{\varepsilon}$ and $\left(\overline{Z}_{\overline{\varepsilon}}\right)_{\overline{\varepsilon}}$
are conjugate (\emph{resp.} orbitally equivalent) by a change of coordinates
and parameters if there exists an analytic diffeomorphism (\emph{resp.}
germ of analytic diffeomorphism) $\left(x,y,\varepsilon\right)\mapsto\left(\Psi_{\varepsilon}\left(x,y\right),\varphi\left(\varepsilon\right)\right)$
such that

\begin{enumerate}
\item $\overline{\varepsilon}=\varphi\left(\varepsilon\right)$ 
\item for fixed $\varepsilon$ the vector fields $Z_{\varepsilon}$ and
$\overline{Z}_{\overline{\varepsilon}}$ are conjugate (\emph{resp.}
orbitally equivalent) under $\Psi_{\varepsilon}$. 
\end{enumerate}
\end{enumerate}
\end{defn}

\subsection{Preparation}

In order to study the analytic classification of families unfolding
a saddle-node it is necessary to \lq\lq prepare them'', so that
the singular points are located on the $x$-axis and their eigenvalues
easily computed from the prepared form.

The following preparation theorem is proved : \\

\textbf{Preparation Theorem.} (see Section~\ref{sec:preparation})
\emph{A representative of a germ of analytic $k$-parameter family
of vector fields unfolding a saddle-node of codimension $k$ is conjugate
by an analytic change of coordinates and parameters over a neighborhood
of the origin in $\mathbb{C}^{2+k}$ to a family of the prepared form
\begin{equation}
Z_{\varepsilon}=U_{\varepsilon}X_{\eps}\label{eq_Z1}\end{equation}
 where \begin{equation}
X_{\eps}(x,y)=P_{\eps}(x)\pp{x}+\left(P_{\eps}(x)R_{0,\varepsilon}\left(x\right)+y\left(1+a(\eps)x^{k}\right)+y^{2}R_{2,\eps}(x,y))\right)\pp{y},\label{eq_X1}\end{equation}
 \[
U_{\varepsilon}(x,y)=Q_{\eps}(x)+P_{\eps}(x)q_{\eps}(x)+O(y),\]
 $\eps=(\eps_{0},\dots,\eps_{k-1})$ is a multi-parameter and \[
\begin{cases}
P_{\eps}(x)=x^{k+1}+\eps_{k-1}x^{k-1}+\dots+\eps_{1}x+\eps_{0}\\
Q_{\eps}(x)=C_{0,\varepsilon}+C_{1,\varepsilon}x+\dots C_{k,\varepsilon}x^{k}.\end{cases}\]
 Here $\varepsilon\mapsto a\left(\varepsilon\right)$, $\varepsilon\mapsto C_{j,\varepsilon}$,
$\left(x,\varepsilon\right)\mapsto q_{\varepsilon}\left(x\right)$,
$\left(x,\varepsilon\right)\mapsto R_{0,\varepsilon}\left(x\right)$
and $\left(x,y,\varepsilon\right)\mapsto R_{2,\varepsilon}\left(x,y\right)$
are germs of holomorphic function and $GCD\left(P_{\varepsilon},Q_{\varepsilon}\right)=1$.
If $\eps=(\eps_{0},\dots,\eps_{k-1})$ and} $\overline{\eps}=(\overline{\eps}_{0},\dots,\overline{\eps}_{k_{1}})$
\emph{we define the equivalence relation over couples} $\left(\varepsilon,a\right)$
with $a\in\germ{\varepsilon}$ : \begin{eqnarray*}
\left(\eps,a\right)\sim\left(\overline{\eps},\overline{a}\right)\,\,\Longleftrightarrow\,\left(\forall j\right)\,\overline{\eps}_{j}=\exp(-2\pi im(j-1)/k)\eps_{j} &  & \,\,\mbox{and }\, a\left(\varepsilon\right)=\overline{a}\left(\overline{\varepsilon}\right).\end{eqnarray*}
 \emph{Then} $\left(\varepsilon,a\right)/\sim$ \emph{is an analytic
invariant.}

\bigskip{}
 This allows to define the model family \emph{\begin{eqnarray}
Z_{\varepsilon}^{M} & := & Q_{\eps}X_{\eps}^{M}\label{eq:model_Z1}\end{eqnarray}
}  with\emph{\begin{eqnarray}
X_{\eps}^{M}\left(x,y\right) & := & P_{\eps}(x)\frac{\partial}{\partial x}+y(1+a(\eps)x^{k})\frac{\partial}{\partial y}\,.\label{eq:model_X1}\end{eqnarray}
} 

\begin{rem}
The germs $\varepsilon\mapsto a\left(\varepsilon\right)$, $\varepsilon\mapsto C_{0,\varepsilon}$,
..., $\varepsilon\mapsto C_{k,\varepsilon}$ are the formal invariants
(they are invariant under formal changes of coordinates in $\left(x,y,\varepsilon\right)$
fibered in the parameter). We can explain their presence in the following
way. When $k+1$ singular points merge in a saddle-node of codimension
$k$ we could expect that all combinations of eigenvalues $(\lambda_{i},\mu_{i})$
would be permitted. As there are only $k$ parameters $\eps_{i}$,
the other degrees of freedom are provided by the formal invariants.
In the case of orbital equivalence it is not the eigenvalues that
are relevant but only their quotients: there are $k+1$ of these,
hence the presence of the formal parameter $a(\eps)$. Of course not
all combinations are possible in a given family, but the class of
families allows for all possibilities. In the conjugacy case there
are $2(k+1)$ eigenvalues, so we need to add the $k+1$ additional
degrees of freedom with the constants $C_{j,\eps}$.

In the case $k=1$, $a(\eps)$ allows for a shift between the quotient
of the eigenvalues at the singular points, one being not necessarily
the inverse of the other. Two additional constants $C_{0,\eps}$ and
$C_{1,\eps}$ allow to determine $\lambda_{0}$ and $\lambda_{1}$,
from which $\mu_{0}$ and $\mu_{1}$ can be found. 
\end{rem}

\subsection{Sectorial decomposition and study of the model}

See Section~\ref{sec:model-family}. The general purpose is to describe
the family of vector fields on a fixed neighborhood $r\ww{D}\times r'\ww{D}$
of the origin in $(x,y)$-space for all values of the parameters in
a fixed neighborhood $\left\{ \left|\left|\varepsilon\right|\right|\leq\rho\right\} $
of the origin in parameter space. In the whole paper we will suppose
that $\rho$ is sufficiently small so that the $k+1$ singular points
coming from the unfolding of the saddle-node remain in $r\ww{D}\times r'\ww{D}$.
We will also suppose that $\rho$ is sufficiently small so that the
whole study is valid on a domain $\left\{ \left|\left|\varepsilon\right|\right|<\rho'\right\} $
with $\rho<\rho'$. The reason for this is that we want to introduce
a conic structure on $\eps$-space from a partition of the sphere
$\left\{ \left|\left|\varepsilon\right|\right|=\rho\right\} $. In
practice we will simply write $\left|\left|\eps\right|\right|\leq\rho$.

The idea is to work with generic $\eps$ for which $P_{\eps}$ has
distinct zeroes and to use the boundedness of the construction to
fill the holes for the other values of $\eps$. If $\Sigma_{0}$ is
the set of generic $\eps$ in a ball of radius $\rho$ where the discriminant
of $P_{\eps}$ does not vanish, then we give a finite covering $\{W_{i}\}_{1\leq i\leq d}$
of $\Sigma_{0}$ with {}``sectors'' $W_{i}$, such that a uniform
treatment can be done over each $W_{i}$ (yielding analytic objects
with respect to $\varepsilon$) and the treatments over different
sectors have the same limit for $\eps=0$.

For a fixed $\eps$ in a given sector $W_{i}$ we divide the phase
space minus the strong separatrices as the union of $2k$ simply connected
domains of the form $V\times r'\mathbb{D}$, where $V$ is a spiraling
sector in $x$-space, which we call \lq\lq squid sector''. Roughly
speaking a good sector $W_{i}$ is defined by the condition that the
length of the spiral is uniformly bounded. The union of the squid
sectors and the singular points is a ball $r\mathbb{D}$ in $x$-space.
The construction of the sectors $V$ is greatly inspired by the work
of Douady and Sentenac \cite{DS}. Each sector is associated to a
sector of the boundary of $r\mathbb{D}$. We expect this construction
to be useful for other problems of moduli of analytic classification,
for instance the problem of the classification of a codimension $k$
parabolic fixed point of a diffeomorphism.

In this partition process, each squid sector $V$ is adherent to two
singular points, one of \lq\lq node type'' (all leaves over the
sector are asymptotic to the point) and one of \lq\lq saddle type''
(a unique leaf is asymptotic to the point over the sector). Note that
for $\eps$ in different sectors $W_{i}$ and for the same sector
of the boundary of $r\mathbb{D}$ we obtain in general different families
of adherent singular points of node and saddle types. As noted by
Douady and Sentenac, the construction could be generalized to the
case of multiple points ($\eps\notin\Sigma_{0}$). In that case the
two adherent singular points of a squid sector of saddle and node
type could be saddle-nodes (and even the same saddle-node), but then
only a saddle sector or a node sector of the saddle-node(s) is included
in the squid sector $V$.

Because of the preparation theorem we have the same squid sectors
for a prepared family and for the associated model family. We prove
a sectorial normalization theorem which is the generalization (an
unfolding) of the theorem of Hukuhara-Kimura-Matuda and show that
over a sector $V$ the foliation is biholomorphic to the model restricted
to the same $V$ (the size of the disk $r'\mathbb{D}$ in $y$-coordinate
has to be adjusted a little). We then show that over these sectors
the space of leaves of the model vector field and of the original
vector field are $\mathbb{C}$. This allows to define almost rigid
coordinates on them, the leaf-coordinates.

We also show the existence of a marked leaf over each squid sector,
corresponding to the weak separatrix of the saddle point attached
to the sector (see Section~\ref{sec:center}). These leaves are called
center manifolds, a name justified by the fact that for $\eps=0$
they indeed coincide with a sectorial center manifold. The leaf-coordinates
are adjusted so as to vanish on the center manifolds.

\subsection{The moduli of analytic classification}

For a given good sector $W_{i}$ in parameter space and the associated
squid sectors $V$ in $x$-space we compare the leaf-coordinates by
means of diffeomorphisms on the intersection of two domains of the
form $V\times r'\mathbb{D}$ , which will allow to define the modulus
for a given $\eps\in W_{i}$. The connected components of the intersection
of two such domains can be of three forms: a sector $V^{s}$ adherent
to a point of saddle type, a sector $V^{n}$ adherent to a point of
node type and a sector $V^{g}$ (for gate) adherent to both. Over
the sectors $V^{g}$ the change of leaf-coordinates is linear. Over
the sector $V^{s}$ the space of leaves is biholomorphic to a disk
and the changes of leaf-coordinates are diffeomorphisms of the form
$h\mapsto h\exp(\phi^{0,i}(h))$, with $\phi^{0,i}\in\germ{h}$ vanishing
at $0$. As there are $k$ sectors this yields $k$ analytic germs
$\phi_{0,\eps}^{0,i},\dots,\phi_{k-1,\eps}^{0,i}$. Over the sectors
$V^{n}$ the changes of leaf-coordinates are given by affine maps
$\psi^{\infty}$, corresponding in particular to changes of center
manifolds. Again there are $k$ such affine maps $\psi_{0,\eps}^{\infty,i},\dots,\psi_{k-1,\eps}^{\infty,i}$.
These maps are defined up to the choice of leaf-coordinates on each
sector, \emph{i.e.} up to linear changes of coordinates. We choose
convenient leaf-coordinates for which the derivative at $0$ of $\psi_{j,\varepsilon}^{\infty,i}$
is $e^{2i\pi a\left(\varepsilon\right)/k}$. This condition forces
the possible changes of leaf-coordinates to be of the special form
$h_{j,\eps}\mapsto c_{\varepsilon}h_{j,\eps}$ with $c_{\varepsilon}\in\ww{C}_{\neq0}$
independent on $j$. On a good sector $W_{i}$ it is possible to choose
the leaf-coordinates so as to depend analytically on $\eps\in W_{i}$
with continuous limit at $\eps=0$. Then the $\phi_{j,\eps}^{0,i}$
and the $\psi_{j,\eps}^{\infty,i}$ can be chosen depending analytically
on $\eps\in W_{i}$ with continuous limit at $\eps=0$. We will always
limit ourselves to this choice. So it is natural to limit the changes
of leaf-coordinates $c_{\eps}$ to belong to the set of functions
\begin{equation}
\mathcal{O}_{b}\left(W_{i}\right)=\{f\in C^{0}(\ov{W}_{i},\ww{C}_{\neq0})\,:\, f\,\text{analytic on }W_{i}\},\label{O_W_i}\end{equation}
 where $\ov{W}_{i}$ is the (topological) closure of $W_{i}$. Besides
a recent result of J. Rib\'{o}n \cite{R} shows this requirement
is necessary : two families of foliations are conjugate if, and only
if, for fixed $\varepsilon$ the two corresponding foliations are
conjugate on a polydisk of size independent on $\varepsilon$.

This allows to state the theorem giving the modulus of analytic classification.
Let us define for any good sector $W_{i}$\begin{eqnarray*}
\mathcal{N}_{\varepsilon}^{i} & := & \left(a,\psi_{0,\eps}^{\infty,i},\dots,\psi_{k-1,\eps}^{\infty,i},\phi_{0,\eps}^{0,i},\dots,\phi_{k-1,\eps}^{0,i}\right)\end{eqnarray*}
 and the equivalence relation \begin{eqnarray*}
\mathcal{N}_{\eps}^{i}\sim\overline{\mathcal{N}}_{\overline{\eps}}^{i} & \Longleftrightarrow & \left(\varepsilon,a\right)\sim\left(\overline{\varepsilon},\overline{a}\right)\mbox{ and for the same }m\in\ww{Z}/k\,:\,\\
 &  & \,\,\,\,\,\,\,\,\,\,\left(\exists c_{\eps}^{i}\in\mathcal{O}_{b}\left(W_{i}\right)\right)\left(\forall j,h\right)\,\,\,\,\begin{cases}
\psi_{j+m,\eps}^{\infty,i}(c_{\eps}^{i}h) & =c_{\eps}^{i}\overline{\psi}_{j,\overline{\eps}}^{\infty,i}(h)\\
\phi_{j+m,\eps}^{0,i}(c_{\eps}^{i}h) & =\overline{\phi}_{j,\overline{\eps}}^{0,i}(h)\,\,.\end{cases}\end{eqnarray*}
 In the work of Martinet-Ramis the modulus for orbital equivalence
of $X_{0}$ corresponds to some $\mathcal{N}_{0}$ satisfying the
same properties. \\

\noindent \textbf{Theorem I.} (see Section~\ref{sec:orbital}) \emph{The
$d$ families of equivalence classes of ($2k+1$)-tuples $\left\{ \mathcal{N}_{\varepsilon}^{i}/\sim\right\} _{\eps\in W_{i}}$,
$1\leq i\leq d$ form a complete modulus of analytic classification
under orbital equivalence for a prepared family $\left(X_{\eps}\right)_{\varepsilon}$
given in \eqref{eq_X1}. If $W_{i}\subset\Sigma_{0}$ is a good sector
then $\mathcal{N}_{\eps}^{i}$ can be chosen bounded and holomorphic
with respect to $\eps\in W_{i}$ and such that its limit for $\eps\to0$
is a fixed $\mathcal{N}_{0}$ independent of the sector $W_{i}$.}
\\

The modulus of analytic classification under conjugacy is composed
of the modulus under orbital equivalence plus a time part. The time
part is formed of the coefficients $C_{j,\eps}$ plus analytic functions
$\zeta_{j,\varepsilon}^{i}\in\germ{h}$. As before the latter functions
measure the obstructions to glue together the transformations of the
system to the model over the squid sectors. The only obstructions
appear on the sectors $V^{s}$. We then build the modulus for $\eps\in W_{i}$
\begin{eqnarray*}
\mathcal{T}_{\varepsilon}^{i} & := & \left(C_{0,\varepsilon},\ldots,C_{k,\varepsilon},\zeta_{0,\eps}^{i},\dots,\zeta_{k-1,\eps}^{i}\right)\end{eqnarray*}
 and extend the equivalence relation $\sim$ to couples\begin{eqnarray*}
\left(\mathcal{N}_{\varepsilon}^{i},\mathcal{T_{\varepsilon}}^{i}\right)\sim\left(\overline{\mathcal{N}}_{\overline{\varepsilon}}^{i},\overline{\mathcal{T}}_{\overline{\varepsilon}}^{i}\right) & \Longleftrightarrow & \mathcal{N}_{\varepsilon}^{i}\sim\overline{\mathcal{N}}_{\overline{\varepsilon}}^{i}\mbox{ and for the same }c_{\varepsilon}^{i}\mbox{ and }m\,:\,\\
 &  & \,\,\,\,\,\,\,\,\,\,\left(\forall j,h\right)\,\,\begin{cases}
C_{j,\varepsilon}e^{2i\pi mj/k} & =\overline{C}_{j,\overline{\varepsilon}}\\
\zeta_{j+m,\eps}^{i}(c_{\varepsilon}^{i}h) & =\overline{\zeta}_{j,\overline{\eps}}^{i}(h)\,.\end{cases}\end{eqnarray*}
 Note that for $\varepsilon=0$ the modulus of conjugacy of Teyssier
and Mershcheryakova-Voronin for the vector field $Z_{0}$ satisfies
the same kind of properties.

This yields the theorem:\\

\noindent \textbf{Theorem II.} (see Section~\ref{sec:conjugacy})
\emph{The $d$ families of equivalence classes of $\left(4k+2\right)$-tuples
$\left\{ \left(\mathcal{N}_{\varepsilon}^{i},\mathcal{T}_{\varepsilon}^{i}\right)/\sim\right\} _{\eps\in W_{i}}$,
$1\leq i\leq d$ is a complete modulus of analytic classification
under conjugacy for a prepared family $Z_{\eps}$ given in \eqref{eq_Z1}.
If $W_{i}\subset\Sigma_{0}$ is a good sector then $\left(\mathcal{N}_{\eps}^{i},\mathcal{T}_{\varepsilon}^{i}\right)$
can be chosen bounded and holomorphic with respect to $\eps\in W_{i}$
and such that the limit for $\eps\to0$ is a fixed $\left(\mathcal{N}_{0},\mathcal{T}_{0}\right)$
independent of the sector $W_{i}$.}

\subsection{The cohomological equation}

See Section~\ref{sec:cohomological}. The tool to prove Theorems
I and II is the solution of a cohomological equation, namely to find
a family of functions $\left(F_{\eps}\right)_{\varepsilon}$ such
that \[
X_{\eps}\cdot F_{\eps}=G_{\eps},\]
 where $\left(G_{\eps}\right)_{\varepsilon}$ is an analytic family
of functions and $X_{\varepsilon}\cdot F_{\varepsilon}$ is the Lie
derivative of $F_{\varepsilon}$ along the vector field $X_{\varepsilon}$.
We solve such an equation over the sectors $V\times r'\mathbb{D}$:
the solution $F_{\eps}$ is given by the integration of $G_{\eps}\tau_{\varepsilon}$
on asymptotic paths with starting point at the singular point of node
type. The difference of two sectorial solutions over non-void intersections
is a first integral of $X_{\varepsilon}$ and thus is constant on
leaves. This allows to write the obstructions to a global solution
as functions of the leaf-coordinates over the intersections $V^{s}$,
$V^{n}$ and $V^{g}$.

To bring the system $X_{\eps}$ of \eqref{eq_X1} to the model \eqref{eq:model_X1}
over $V\times r'\mathbb{D}$ we must bring the weak invariant manifold
of the point of saddle type to the horizontal axis $y=0$. We then
use a change of coordinates preserving $y=0$ in the form of the flow
of $y\frac{\partial}{\partial y}$ for some time $N_{\eps}$. This
time $N_{\eps}$ is found as a solution of a first cohomological equation
of the form \[
X_{\eps}\cdot N_{\eps}=\tilde{R}_{\eps}\]
 for an appropriate function $\tilde{R}_{\eps}$.

As for bringing the system $Z_{\eps}=U_{\varepsilon}X_{\varepsilon}$
of \eqref{eq_Z1} to the model \eqref{eq:model_X1} over $V\times r'\mathbb{D}$
we compose the previous change of coordinates with a change of coordinates
taking care of the time part. This change of coordinates is given
by the flow of the vector field $Q_{\eps}X_{\eps}$ for some time
$T_{\eps}$. The time $T_{\eps}$ is the solution of a second cohomological
equation of the form \[
X_{\eps}\cdot T_{\eps}=\frac{1}{U_{\eps}}-\frac{1}{Q_{\eps}}.\]

\section{\label{sec:preparation}Preparation of the family}

This section deals with a germ of analytic family unfolding a germ
of saddle-node of codimension $k\in\ww{N}_{\neq0}$. We consider a
germ of generic (to be defined below) analytic $k$-parameter family
unfolding a germ of vector field with a saddle-node of codimension
$k$. It is known that, up to a local analytic change of coordinates,
a representative of the germ of vector field can be taken in Dulac's
prenormal form \begin{equation}
\begin{array}{lll}
Z_{0} & := & U_{0}X_{0}\\
X_{0}\left(x,y\right) & := & x^{k+1}\pp{x}+\left(y\left(1+ax^{k}\right)+x^{k+1}R(x,y)\right)\pp{y}\end{array}\label{sad.node}\end{equation}
 with $U_{0}(x,y)=C_{0}+C_{1}x+\dots+C_{k}x^{k}+O(x^{k+1})+O(y)$
and $a,\, C_{j}\in\ww{C}$ with $C_{0}\neq0$.

\bigskip{}
 We consider an unfolding \begin{equation}
Z_{0}\left(x,y\right)+H_{1}\left(x,y,\eta\right)\pp{x}+H_{2}\left(x,y,\eta\right)\pp{y}\label{unfold1}\end{equation}
 where $H_{j}(x,y,\eta)=O(\eta)$ is a germ at $\left(0,0,0\right)$
of a holomorphic function and $\eta=(\eta_{0},\dots,\eta_{k-1})$
is a multi-parameter in a neighborhood of $0$ in $\mathbb{C}^{k}$.
We make the change of coordinates \begin{eqnarray*}
\left(\tilde{x},\tilde{y}\right) & := & \left(x,\left(y(1+ax^{k})+x^{k+1}R(x,y)\right)U_{0}(x,y)+H_{2}(x,y,\eta)\right)\,.\end{eqnarray*}
 Then all singular points occur on $\tilde{y}=0$. Given that the
origin of \eqref{sad.node} is of multiplicity $k$ there are $k+1$
small singular points of \eqref{unfold1} on $\tilde{y}=0$. \emph{Modulo}
a translation in the variable $\tilde{x}$ we can suppose that they
are roots of \begin{eqnarray}
P_{\eps}(\tilde{x}) & := & \tilde{x}^{k+1}+\eps_{k-1}\tilde{x}^{k-1}+\dots+\eps_{1}\tilde{x}+\eps_{0}.\end{eqnarray}
 The family is \textbf{generic} if the change of parameters $(\eta_{0},\dots,\eta_{k-1})\mapsto(\eps_{0},\dots,\eps_{k-1})$
is an analytic isomorphism in a neighborhood of the origin. Then the
family of vector fields has the form \begin{equation}
\left(P_{\eps}\left(\tilde{x}\right)h_{1}(\tilde{x},\eps)+\tilde{y}h_{3}(\tilde{x},\tilde{y},\eps)\right)\pp{\tilde{x}}+\left(P_{\eps}\left(\tilde{x}\right)h_{2}(\tilde{x},\eps)+\tilde{y}k_{1}(\tilde{x},\tilde{y},\eps)\right)\pp{\tilde{y}}\label{unfold2}\end{equation}
 with $h_{1}(0,0)k_{1}(0,0,0)\neq0$, for some $h_{1},h_{2}\in\germ{\tilde{x},\varepsilon}$
and $h_{3},k_{1}\in\germ{\tilde{x},\tilde{y},\varepsilon}$.

\bigskip{}
 It is possible to adapt the technique of Glutsyuk and straighten
all strong manifolds of singular points uniformly on a neighborhood
of the origin.

\begin{prop}
\label{pro:straighten}There exists an analytic change of coordinates
$\left(\tilde{x},\tilde{y}\right)\mapsto\left(x,y\right)$ on a neighborhood
of the origin in $\mathbb{C}^{2}$ depending holomorphically on $\eps$
for $\eps$ in a neighborhood of the origin and conjugating the vector
field \eqref{unfold2} to \begin{equation}
Z_{\varepsilon}=U_{\varepsilon}X_{\eps}\label{unfold3}\end{equation}
 where \begin{equation}
X_{\eps}\left(x,y\right)=P_{\eps}(x)\pp{x}+\left(P_{\eps}(x)R_{0,\eps}(x)+y\left(1+a(\eps)x^{k}\right)+y^{2}R_{2,\varepsilon}(x,y))\right)\pp{y}\label{eq_X0}\end{equation}
 and \begin{equation}
U_{\varepsilon}(x,y)=Q_{\eps}(x)+O(P_{\eps}(x))+O(y),\label{eq_U}\end{equation}
 with \begin{equation}
Q_{\eps}\left(x\right):=C_{0,\varepsilon}+C_{1,\varepsilon}x+\dots C_{k,\varepsilon}x^{k}.\label{eq_Q}\end{equation}
 Here $\varepsilon\mapsto a\left(\varepsilon\right)$, $\varepsilon\mapsto C_{j,\varepsilon}$,
$\left(x,\varepsilon\right)\mapsto R_{0,\varepsilon}\left(x\right)$,
$\left(x,y,\varepsilon\right)\mapsto U_{\varepsilon}\left(x,y\right)$
and $\left(x,y,\varepsilon\right)\mapsto R_{2,\varepsilon}\left(x,y\right)$
are germs of holomorphic functions at the origin. Moreover $GCD\left(Q_{\varepsilon},P_{\varepsilon}\right)=1$,
which in particular means \emph{$C_{0,0}\neq0$} when \emph{$\varepsilon=0$.} 
\end{prop}
\begin{proof}
The proof contains several steps. For the first step we write \eqref{unfold2}
as $h_{1}(\tilde{x},\eps)\hat{X}_{\eps}$. Then \begin{equation}
\hat{X}_{\eps}\left(\tilde{x},\tilde{y}\right)=\left(P_{\eps}(\tilde{x})+\tilde{y}h_{4}(\tilde{x},\tilde{y},\eps)\right)\frac{\partial}{\partial\tilde{x}}+\left(P_{\eps}(\tilde{x})h_{3}(\tilde{x},\eps)+\tilde{y}(1+O(\tilde{x},\tilde{y},\eps))\right)\frac{\partial}{\partial\tilde{y}}\label{eq_hat_X}\end{equation}
 where $h_{3},h_{4}\in\mathbb{C}\{\tilde{x},\tilde{y},\eps\}$ and
we look for a change of coordinates for $\hat{X}_{\eps}$ straightening
simultaneously all separatrices. We consider a small neighborhood
$W$ of the origin in $\eps$-space and the open subset $\Sigma_{0}$
of generic $\eps$ such that the discriminant of $P_{\eps}$ does
not vanish. Then $P_{\eps}$ has $k+1$ distinct roots $x_{0}(\eps),\dots,x_{k}(\eps)$.
There exists a neighborhood $V_{y}$ of the origin in $\tilde{y}$-space,
independent of $\eps$, such that the strong manifold of $(x_{j},0)$
is given by $\tilde{x}=F_{j}(\tilde{y})$ over $V_{y}$. The proof
of that later fact is done as in \cite{G}. The idea is the following:
we consider the cones $K_{j}=\left\{ (\tilde{x},\tilde{y})\,:\,|\tilde{y}|\geq|\tilde{x}-x_{j}(\eps)|\right\} $.
On such a cone we have $|\frac{dy}{dx}|>1$ for $\eps$ sufficiently
small (as $|\dot{\tilde{x}}|<1/2|\tilde{x}-x_{j}(\eps)|$ and $|\dot{\tilde{y}}|>1/2|\tilde{x}-x_{j}(\eps)|$
for $(\tilde{x},\tilde{y})$ in a small neighborhood of the origin
and $\eps$ sufficiently small). If $\chi(t):=|\tilde{y}(t)|$ we
also have that $\dot{\chi}>\frac{1}{2}\chi>0$. As the local invariant
manifold is given by $\tilde{x}=F_{j}(\tilde{y})=x_{j}(\eps)+c(\eps)\tilde{y}+o(\tilde{y})$,
with $F_{j}$ analytic and $c(\eps)$ small, then the graph of $F_{j}$
is contained in the cone $K_{j}$ for small $\tilde{y}$. Now the
extension of the invariant manifold is the union of all real trajectories
with positive time starting on points $(F_{j}(\delta e^{i\theta}),\delta e^{i\theta})$
for $\theta\in[0,2\pi]$ and $\delta>0$ small. These trajectories
remain in $K_{j}$, so there is no obstruction to extend the graph
of $F_{j}$ to a fixed neighborhood $V_{y}$. \medskip{}
 The straightening change of coordinates is then given by $\left(\tilde{x},\tilde{y}\right)\mapsto\left(G(\tilde{x},\tilde{y},\eps),\tilde{y}\right)=\left(\hat{x},\tilde{y}\right)$
where\begin{equation}
G\left(\tilde{x},\tilde{y},\varepsilon\right):=\sum_{j=0}^{k}x_{j}(\eps)\prod_{l\neq j}\frac{\tilde{x}-F_{l}(\tilde{y})}{F_{j}(\tilde{y})-F_{l}(\tilde{y})}\,.\label{change_X}\end{equation}
 It is holomorphic for $(\tilde{x},\eps)$ small and $\tilde{y}\in V_{y}$.
The holomorphy in $\eps$ follows from the invariance under permutations
of the $x_{j}$. Moreover it has a holomorphic extension to $\Sigma_{0}\cup\Sigma_{1}$
where $\Sigma_{1}$ is the set of $\eps$ for which $P_{\eps}(x)$
has exactly one double root. Indeed $G\left(\tilde{x},\tilde{y},\varepsilon\right)$
is given by the following formula: \begin{equation}
G\left(\tilde{x},\tilde{y},\varepsilon\right)=\frac{\left|\begin{array}{lllll}
0 & 1 & \tilde{x} & \dots & \tilde{x}^{k}\\
x_{0}(\eps) & 1 & F_{0}(\tilde{y}) & \dots & F_{0}^{k}(\tilde{y})\\
\vdots & \vdots & \vdots & \vdots & \vdots\\
x_{k}(\eps) & 1 & F_{k}(\tilde{y}) & \dots & F_{k}^{k}(\tilde{y})\end{array}\right|}{\left|\begin{array}{llll}
1 & F_{0}(\tilde{y}) & \dots & F_{0}^{k}(\tilde{y})\\
\vdots & \vdots & \vdots & \vdots\\
1 & F_{k}(\tilde{y}) & \dots & F_{k}^{k}(\tilde{y})\end{array}\right|}.\label{interpolation}\end{equation}
 From now on we write $x_{j}$ for $x_{j}(\eps)$. If $\ov{\eps}\in\Sigma_{1}$
and $x_{0}$ and $x_{1}$ coalesce as $\eps\to\ov{\eps}$, then the
limit of $G\left(\tilde{x},\tilde{y},\varepsilon\right)$ exists provided
that $\lim_{\eps\to\ov{\eps}}\frac{F_{0}(\tilde{y})-F_{1}(\tilde{y})}{x_{0}-x_{1}}$
exists. This can be seen by subtracting the row corresponding to $x_{0}$
to the row corresponding to $x_{1}$, both in the numerator and the
denominator of \eqref{interpolation}, and then by dividing both the
numerator and the denominator by $x_{0}-x_{1}$, thus removing the
indeterminacy. Here we have $\lim_{\eps\to\ov{\eps}}\frac{F_{0}(\tilde{y})-F_{1}(\tilde{y})}{x_{0}-x_{1}}=0$.
The latter has been proved by Glutsyuk \cite{G}. It comes from showing
that $\lim_{\eps\to\ov{\eps}}F_{j}'(\tilde{y})=0$, $j=0,1$, uniformly
over $V_{y}$. Let us take the case $j=1$. The details are as follows.
Without loss of generality we can suppose that one separatrix, for
instance that of $x_{0}$, has been straightened to $F_{0}(\tilde{y})\equiv x_{0}$,
which implies that $h_{4}(\tilde{x},\tilde{y},\eps)=(\tilde{x}-x_{0})h_{5}(\tilde{x},\tilde{y},\eps)$
in \eqref{eq_hat_X}. To prove that $\lim_{\eps\to\ov{\eps}}F_{1}'(\tilde{y})=0$,
Glutsyuk proves that $\lim_{\eps\to\ov{\eps}}\frac{F_{1}'(\tilde{y})}{F_{1}(\tilde{y})-x_{0}}$
is bounded.

\begin{eqnarray}
\begin{array}{lll}
F_{1}'(\tilde{y}) & = & \left.\frac{d\tilde{x}}{d\tilde{y}}\right|_{\tilde{x}=F_{1}(\tilde{y})}\\
 & = & \frac{(F_{1}(\tilde{y})-x_{0})(F_{1}(\tilde{y})-x_{1})\prod_{k\neq0,1}(F_{1}(\tilde{y})-x_{k})+\tilde{y}(F_{1}(\tilde{y})-x_{0})h_{5}(\tilde{x},\tilde{y},\eps)}{\tilde{y}(1+O(|F_{1}(\tilde{y}),\tilde{y},\eps|))+P_{\eps}(F_{1}(\tilde{y}))h_{3}(F_{1}(\tilde{y}),\eps))},\end{array}\end{eqnarray}
 with $h_{5}\in\mathbb{C}\{\tilde{x},\tilde{y}\}$ depending continuously
on $\eps$. Then \begin{eqnarray}
\frac{F_{1}'(\tilde{y})}{F_{1}(\tilde{y})-x_{0}}=\frac{\frac{F_{1}(\tilde{y})-x_{1}}{\tilde{y}}\prod_{k\neq0,1}(F_{1}(\tilde{y})-x_{k})+h_{5}(\tilde{x},\tilde{y},\eps)}{1+O(|F_{j}(\tilde{y}),\tilde{y},\eps|)+\frac{F_{1}(\tilde{y})-x_{1}}{\tilde{y}}\prod_{k\neq1}(F_{1}(\tilde{y})-x_{k})h_{3}(F_{1}(\tilde{y}),\eps)}\end{eqnarray}
 The conclusion follows as $\left|\frac{F_{1}(\tilde{y})-x_{1}}{\tilde{y}}\right|<1$
since we are in the cone $K_{1}$ (details in \cite{G}).

The extension of $G$ to $\Sigma_{0}\cup\Sigma_{1}$ depends analytically
on $\eps$ as it is again invariant under permutations of the $x_{j}$.
Since the complement of $\Sigma_{0}\cup\Sigma_{1}$ is of codimension
2 then, by Hartog's theorem, we can extend $G$ to all values of $\eps$,
satisfying $|\eps|\leq\rho$ for some positive $\rho$.

The change of coordinate \eqref{change_X} allows to factor $P_{\eps}(\hat{x})$
in the first component of the vector field. Hence it has the form
\[
P_{\eps}(\hat{x})U(\hat{x},\tilde{y},\eps)\pp{\hat{x}}+\left(P_{\eps}(\hat{x})h_{6}(\hat{x},\eps)+\tilde{y}h_{7}(\hat{x},\eps)+O(\tilde{y}^{2})\right)\pp{\tilde{y}}\]
 where $U(0,0,0)=h_{7}(0,0)\neq0$ and $h_{6},h_{7}\in\mathbb{C}\{\tilde{x},\eps\}$.
We factorize $U(\hat{x},y,\eps)$ in the vector field which then has
the form \[
U(\hat{x},\tilde{y},\eps)\left[P_{\eps}(\hat{x})\pp{\hat{x}}+\left(P_{\eps}(\hat{x})h_{8}(\hat{x},\eps)+\tilde{y}h_{9}(\hat{x},\eps)+O(\tilde{y}^{2})\right)\pp{\tilde{y}}\right]\]
 with $h_{8},h_{9}\in\mathbb{C}\{\tilde{x},\eps\}$ and $h_{9}(0,0)=1$.

As in \cite{K} we use a change of coordinate and parameter $(\hat{x},\eps)\mapsto(x,\tilde{\eps})$
to transform $\frac{P_{\eps}(\hat{x})}{h_{9}(\hat{x},\eps)}\frac{\partial}{\partial\hat{x}}$
into $\frac{P_{\tilde{\eps}}(x)}{1+a(\tilde{\varepsilon})x^{k}}\frac{\partial}{\partial x}$
to bring the vector field to the final form \eqref{eq_X0}. This ends
the preparation of the orbital part.

For the time part, the vector field has the form $X_{\tilde{\eps}}\hat{U}_{\tilde{\eps}}$.
We simply divide the $x$-part of $\hat{U}_{\tilde{\eps}}$ by $P_{\tilde{\eps}}(x)$:
\[
\hat{U}_{\tilde{\eps}}=Q_{\tilde{\eps}}(x)+P_{\tilde{\eps}}(x)q_{\tilde{\eps}}(x)+O(y).\]

\end{proof}
\bigskip{}

\begin{defn}
\label{def:prepared}From now on we always work with germs of analytic
families \begin{equation}
Z_{\varepsilon}=U_{\varepsilon}X_{\eps}\label{unfold}\end{equation}
 where \begin{equation}
X_{\eps}(x,y)=P_{\eps}(x)\pp{x}+\left(P_{\eps}(x)R_{0,\varepsilon}\left(x\right)+y\left(1+a(\eps)x^{k}\right)+y^{2}R_{2,\eps}(x,y))\right)\pp{y}\label{eq_X}\end{equation}
 and \begin{equation}
U_{\varepsilon}(x,y)=Q_{\eps}(x)+O\left(P_{\varepsilon}\left(x\right)\right)+O(y),\end{equation}
 which we call \textbf{\textit{\emph{prepared families}}}. Here $Q_{\varepsilon}$
and $R_{j,\varepsilon}$ are characterized in Proposition~\ref{pro:straighten}. 
\end{defn}
\begin{thm}
\label{pro:eqv_eps}~ 
\begin{enumerate}
\item A transformation $x\mapsto\tilde{x}=\exp(2\pi i\, m/k)x$ with $m=1,\dots,k-1$,
transforms a prepared family into a prepared family with corresponding
polynomials $P_{\tilde{\eps}}(\tilde{x})=\tilde{x}^{k+1}+\tilde{\eps}_{k-1}\tilde{x}^{k-1}+\dots+\tilde{\eps}_{1}\tilde{x}+\tilde{\eps}_{0}$,
where $\tilde{\eps}_{j}=\exp(-2\pi i\, m(j-1)/k)\eps_{j}$ and $\tilde{Q}_{\tilde{\eps}}(\tilde{x})=\tilde{C}_{0,\tilde{\varepsilon}}+\dots+\tilde{C}_{k,\tilde{\varepsilon}}\tilde{x}^{k}$
where $\tilde{C}_{j,\tilde{\varepsilon}}=\exp(-2\pi i\, mj/k)C_{j,\varepsilon}$. 
\item The $2(k+1)$ eigenvalues of the $k+1$ singular points of \eqref{unfold}
given by $(x_{j},0)$, $j=0,\dots,k$, where $x_{j}$ are the roots
of $P_{\eps}$, coincide with that of the \textbf{model family} \begin{equation}
Z_{\eps}^{M}(x,y):=Q_{\eps}(x)X_{\eps}^{M}(x,y)=Q_{\eps}(x)\left[P_{\eps}(x)\frac{\partial}{\partial x}+y(1+a(\eps)x^{k})\frac{\partial}{\partial y}\right]\label{model}\end{equation}
 where $Q_{\eps}$ is given in \eqref{eq_Q}. 
\item Suppose that two prepared families $\left(X_{\eps}\right)_{\varepsilon}$
and $\left(\tilde{X}_{\tilde{\eps}}\right)_{\tilde{\varepsilon}}$
are conjugate. We define the equivalence relations \begin{eqnarray}
\eps\cong\tilde{\eps}\Longleftrightarrow\left(\exists m\in\ww{Z}/k\right)\,\,\tilde{\eps}_{j} & = & \exp(-2\pi i\, m(j-1)/k)\eps_{j}\qquad j=0,\dots,k-1.\label{cong1}\\
\left(\varepsilon,a\right)\cong\left(\tilde{\varepsilon},\tilde{a}\right)\Longleftrightarrow\varepsilon\cong\tilde{\varepsilon} &  & \mbox{and}\,\,\tilde{a}\left(\tilde{\varepsilon}\right)=a\left(\varepsilon\right)\end{eqnarray}
 where $a$ is given in \eqref{eq_X0}. The equivalence classes $[\left(\varepsilon,a\right)]/\sim$
are analytic invariants. 
\item Let $\eps:=\left(\eps_{0},\dots,\eps_{k-1}\right)$ and $C_{\varepsilon}:=\left(C_{0,\varepsilon},\dots,C_{k,\varepsilon}\right)$,
where the $C_{j,\varepsilon}$'s are the coefficients of $Q_{\varepsilon}$.
We define the equivalence relation \begin{equation}
\left(\eps,C_{\varepsilon}\right)\cong\left(\tilde{\eps},\tilde{C}_{\tilde{\varepsilon}}\right)\Longleftrightarrow\left(\exists m\in\ww{Z}/k\right)\,\,\begin{cases}
\tilde{\eps}_{j}=\exp(-2\pi i\, m(j-1)/k)\eps_{j} & \,\,\,\,\, j=0,\dots,k-1,\\
\tilde{C}_{j,\tilde{\varepsilon}}=\exp(-2\pi i\, mj/k)C_{j,\varepsilon} & \,\,\,\,\, j=0,\dots,k.\end{cases}\label{cong2}\end{equation}
 The equivalence classes $\left[\left(\eps,C_{\varepsilon}\right)\right]/\cong$
are analytic invariants of \eqref{unfold}. 
\end{enumerate}
\end{thm}
\bigskip{}
 \emph{Proof of Theorem~\ref{pro:eqv_eps}.} Only the third and fourth
items require a proof. We write $a$ instead of $a(\eps)$.

\noindent \textbf{(3)} Let $(x_{j},0)$, $j=0,\dots,k$, be the singular
points of \eqref{eq_X}. The ratios of eigenvalues at each singular
point is an analytic invariant under orbital equivalence. When $\eps\in\Sigma_{0}$
these are given by \begin{eqnarray}
\nu_{j} & = & \frac{1+ax_{j}^{k}}{P_{\eps}'(x_{j})}.\label{eq_nu}\end{eqnarray}
 Then \begin{equation}
a=\sum_{j=0}^{k}\frac{1}{\nu_{j}}.\end{equation}
 The quantity $a$ remains bounded when two singular points collide
as \begin{equation}
a=\frac{1}{2\pi i}\int_{r\mathbb{S}^{1}}\frac{1+az^{k}}{P_{\eps}(z)}dz\label{eq_a}\end{equation}
 where $r\mathbb{S}^{1}$ is a circle in $x$-space surrounding $x_{0},\dots,x_{k}$.

We suppose that there is an equivalence $(x,y,\eps)\mapsto\left(\Psi_{\eps}(x,y),h(\eps)\right)$
between the prepared families $\left(X_{\eps}\right)$ and $\left(\tilde{X}_{h(\eps)}\right)$.
This yields an equivalence between $Y_{\varepsilon}:=X_{\eps}/(1+ax^{k})$
and $\tilde{X}_{h(\eps)}/(1+\tilde{a}\left(h\left(\varepsilon\right)\right)\tilde{x}^{k})$.
The map $\Psi_{\eps}(x,y)=\left(H_{1,\varepsilon}(x,y),H_{2,\varepsilon}(x,y)\right)$
sends the singular points to the singular points. Let us first show
that it is possible to construct an equivalence $\Theta_{\varepsilon}=(K_{1,\varepsilon},K_{2,\varepsilon})$
in which the first coordinate $K_{1,\varepsilon}$ depends on $x$
alone. The ideas come from \cite{MR1}.

The map $H_{1,\varepsilon}(x,y)=h_{1,\varepsilon}(x)+r_{\varepsilon}(x,y)$
with $r_{\eps}(x,y)=O(y)$, as a fibration $H_{1,\varepsilon}\,:\,(\mathbb{C}^{2},0)\rightarrow(\mathbb{C},0)$,
is transverse to all fibers, except the fibers through the singular
points. In particular the gradient of $H_{1,\varepsilon}$ is orthogonal
to the vector field along the fibers through the singular points.
This yields that $r_{\varepsilon}(x,y)=P_{\eps}(x)r_{1,\varepsilon}(x,y)$,
with $r_{1,\cdot}$ analytic in $\left(x,y,\varepsilon\right)$. As
in \cite{MR1} we can construct an analytic change of coordinates
$\Theta_{\varepsilon}$ such that $H_{1,\varepsilon}\circ\Theta_{\varepsilon}=h_{1,\varepsilon}$
(see Lemma~\ref{lem:MR} below) which is an equivalence between $Y_{\eps}$
and itself.

So we can suppose that there exists an equivalence $(x,y,\eps)\mapsto\left(\Psi_{\eps}(x,y),h(\eps)\right)$
between the two families, in which $H_{1,\varepsilon}$ depends on
$x$ alone. Then the map $H_{1,\varepsilon}$ is a conjugacy between
$W_{\eps}\left(x\right)=P_{\eps}(x)/(1+ax^{k})\pp{x}$ and $V_{h(\eps)}(\tilde{x})\tilde{W}_{h(\eps)}$,
where $\tilde{W}_{h(\eps)}=\tilde{P}_{h(\eps)}/(1+\tilde{a}\left(h\left(\varepsilon\right)\right)\tilde{x}^{k})\frac{\partial}{\partial\tilde{x}}$
and $\left(\tilde{x},\varepsilon\right)\mapsto V_{h\left(\varepsilon\right)}\left(\tilde{x}\right)$
is an analytic map. By Kostov's Theorem, there exists a germ of conjugacy
$x\mapsto K_{1,\varepsilon}(x)$ between $V_{h(\eps)}\tilde{W}_{h(\eps)}$
and $\tilde{W}_{h(\eps)}$. Let $\Theta_{\varepsilon}(\tilde{x},\tilde{y}):=(K_{1,\varepsilon}(\tilde{x}),\tilde{y})$.
Then $\Theta_{\varepsilon}\circ\Psi_{\varepsilon}$ is an equivalence
between $X_{\eps}$ and $\Theta_{\varepsilon}^{*}\left(\tilde{X}_{h\left(\varepsilon\right)}\right)$
while $L_{\varepsilon}:=K_{1,\varepsilon}\circ H_{1,\varepsilon}$
is a conjugacy between $W_{\eps}$ and $\tilde{W}_{h(\eps)}$. \medskip{}

\textbf{The case $k=1$.} Then $\eps_{0}$ is an analytic invariant
since \[
\frac{1}{\sqrt{-\eps_{0}}}=\frac{1}{\nu_{0}}-\frac{1}{\nu_{1}},\]
 where $x_{0}=\sqrt{-\eps_{0}}$ and $x_{1}=-\sqrt{-\eps_{0}}$. This
comes from \cite{IY} (see also \cite{R2}). \medskip{}

\textbf{The case $k>1$.} It is done in Theorem~\ref{parametres_inv}
below.

\textbf{(4)} We have shown in (3) that $[\left(\eps,a\right)]/\cong$
is an analytic invariant of $X_{\eps}$ and that the equivalence relation
$\cong$ yields an equivalence relation on $x$ given by $x\cong\tilde{x}=\exp(2\pi i\, m/k)x$.
This yields that the set of singular points $\left[\left\{ x_{0},\dots,x_{k}\right\} \right]/\cong$
is an analytic invariant for a prepared family. The eigenvalues of
the linearized vector field of $Z_{\eps}$ at $(x_{j},0)$ are analytic
invariants of the system. They are given by \begin{equation}
(\lambda_{j},\mu_{j})=(P_{\eps}'(x_{j})Q_{\eps}(x_{j}),Q_{\eps}(x_{j})).\end{equation}
 The coefficients $C_{j}(\eps)$ of $Q_{\eps}$ are uniquely determined
from the $x_{j}$ using $Q_{\eps}(x_{j})=\mu_{j}$. \hfill{}$\Box$

\begin{lem}
\label{lem:MR} We consider the vector field \eqref{eq_X}. Let $H_{1,\eps}(x,y)=h_{1,\eps}(x)+P_{\eps}(x)r_{1,\eps}(x,y)$
with $r_{1,\eps}(x,y)=O(y)$, be a family of analytic maps defined
in a neighborhood of the origin in $\mathbb{C}^{2}$ such that $\frac{\partial H_{1,\eps}}{\partial x}$
never vanishes on that neighborhood. There exists a family of analytic
diffeomorphisms $K_{\eps}$ defined over a neighborhood of the origin
in $\mathbb{C}^{2}$ which is an equivalence between \eqref{eq_X}
and itself (\emph{i.e.} an orbital symmetry of \eqref{eq_X}) and
such that $H_{1,\eps}\circ K_{\eps}=h_{1,\eps}$. 
\end{lem}
\begin{proof}
The proof is an adaptation of Lemma (2.2.2), Chapter II of \cite{MR1}.
It is done by the homotopy method. Let \begin{equation}
H_{\eps}(t,x,y)=H_{t,\eps}(x,y)=h_{1,\eps}(x)+tP_{\eps}(x)r_{1,\eps}(x,y)\end{equation}
 and $\omega_{\eps}$ be a 1-form dual to $X_{\varepsilon}$. We look
for a one-parameter family of analytic vector fields $\Xi_{t,\eps}(x,y)$
such that \begin{equation}
\omega_{\eps}(\Xi_{t,\eps})=0.\label{eq_1}\end{equation}
 and \begin{equation}
\Xi_{t,\eps}\cdot H_{1,\eps}=-\frac{\partial H_{t,\eps}}{\partial t}=-P_{\eps}r_{1,\eps}.\label{eq_2}\end{equation}

Then \eqref{eq_1} yields $\Xi_{t,\eps}=g_{\eps}X_{\eps}$ for some
arbitrary family of functions $g_{\eps}$. As $X_{\eps}=P_{\eps}\pp{x}+A_{\eps}\pp{y}$,
\eqref{eq_2} yields \begin{equation}
g_{\eps}\left(x,y\right)\left(\frac{\partial H_{t,\eps}}{\partial x}\left(x,y\right)+tA(x,y)\frac{\partial r_{1,\eps}}{\partial y}\left(x,y\right)\right)=-r_{1,\eps}(x,y),\end{equation}
 which has an analytic solution $g$ over a neighborhood of the origin
since $\frac{\partial H_{t,\eps}}{\partial x}\neq0$ for $y$ sufficiently
small.

The time-$t$ flow of $\Xi_{t,\eps}$ is a diffeomorphism $K_{t,\eps}$
which is an orbital symmetry of $X_{\eps}$ and such that $H_{1,\eps}\circ K_{1,\eps}=h_{1,\eps}$. 
\end{proof}
\bigskip{}

\begin{thm}
\label{parametres_inv} We consider a germ of an analytic change of
coordinates $\Psi\,:\,\left(x,\varepsilon\right)=\left(x,\varepsilon_{0},\ldots,\varepsilon_{k-1}\right)\mapsto\left(\varphi_{\varepsilon}\left(x\right),h_{0}\left(\varepsilon\right),\ldots,h_{k-1}\left(\varepsilon\right)\right)=\left(z,h\right)$
at $\left(0,0,\cdots,0\right)\in\mathbb{C}^{1+k}$. The following
assertions are equivalent : 
\begin{enumerate}
\item the families $\left(\frac{P_{\varepsilon}\left(x\right)}{1+a\left(\varepsilon\right)x^{k}}\frac{\partial}{\partial x}\right)_{\varepsilon}$
and $\left(\frac{P_{h}\left(z\right)}{1+\tilde{a}\left(h\right)z^{k}}\frac{\partial}{\partial z}\right)_{h}$
are conjugate under $\Psi$, 
\item there exist $\lambda$ with $\lambda^{k}=1$ and $T\in\germ{\varepsilon}$
such that

\begin{itemize}
\item $\varphi_{\varepsilon}\left(x\right)=\Phi_{X_{\eps}}^{T(\eps)}\circ R_{\lambda}\left(x\right)$
where $R_{\lambda}(x)=\lambda x$, 
\item $\varepsilon_{j}=\lambda^{j-1}h_{j}\left(\varepsilon\right)$, 
\item $a\left(\varepsilon\right)=\tilde{a}\left(h\left(\varepsilon\right)\right)$. 
\end{itemize}
\end{enumerate}
\end{thm}
\begin{proof}
(2)$\Rightarrow$(1) is trivial so we only consider (1)$\Rightarrow$(2).
We may moreover assume that $k>1$, the case $k=1$ being recalled
in Theorem~\ref{pro:eqv_eps}.

The result is easily shown for $\eps=0$. But let us discuss some
details which will be important in the proof. Indeed the flow $\Phi_{X_{0}}^{t}$
has the form $x(1+g_{t}(x^{k}))$, with $g_{t}(0)=0$. Moreover if
$\varphi_{0}'(0)=\lambda_{0}$ we need have $\lambda_{0}^{k}=1$ in
order to preserve the form of $X_{0}$. So we can compose $\Psi(x,\eps)$
with $R_{\lambda_{0}}$ and the corresponding change of parameters
$\varepsilon_{j}=\lambda_{0}^{j-1}h_{j}\left(\varepsilon\right)$
and only discuss the composed family. Hence we can suppose that $\Psi(x,\eps)$
is such that $\varphi_{0}'(0)=1$. We now need to prove that $h_{j}(\eps)\equiv\eps_{j}$.

It is easy to check that the only changes of coordinates tangent to
the identity which preserve $X_{0}$ are the maps $\Phi_{X_{0}}^{t}$:
using power series, it is easily verified that such changes of coordinates
have the form $x(1+m_{t}(x^{k}))$ with $m_{t}(0)=0$, where $m_{t}$
is completely determined by $m_{t}'(0)=t$. This is exactly the form
of the family $\Phi_{X_{0}}^{t}$. Indeed let $\Phi_{X_{0}}^{t}(x)=b_{t}(x)$.
The function $b_{t}(x)=xd_{t}(x)$ with $d_{t}(0)=1$ is solution
of \[
-\frac{1}{kb_{t}^{k}(x)}+\frac{1}{kx^{k}}+a\ln(b_{t}(x))-a\ln(x)=t\]
 \emph{i.e.} \[
d_{t}^{k}(x)-1+akx^{k}d_{t}^{k}(x)\ln(d_{t}(x))=ktx^{k}d_{t}^{k}(x).\]
 Substituting an unknown power series $d_{t}(x)=1+\sum_{n\geq1}c_{n}x^{n}$
yields the result.

Let \[
G(x,t,\eps):=\Phi_{X_{\eps}}^{t}\circ\varphi_{\eps}\left(x\right)\,,\]
 \[
H(x,t,\eps):=\frac{\partial^{k+1}G}{\partial x^{k+1}}\left(x,t,\varepsilon\right)\]
 and \[
K(t,\eps):=H(0,t,\eps).\]
 $K$ is an analytic map and we have \[
\frac{\partial K}{\partial t}(0,0)=(k+1)!\neq0.\]
 Moreover, let $t_{0}$ be such that $K(t_{0},0)=0$ (in the study
for $\eps=0$ we have shown the existence of $t_{0}$). By the implicit
function theorem there exists a unique function $t(\eps)$ such that
$t(0)=t_{0}$ and $K(t(\eps),\eps)\equiv0$. Composing $\varphi_{\eps}$
with $\Phi_{X_{\eps}}^{t(\eps)}$ we can suppose that the original
family $\Psi$ is such that $\frac{\partial^{k+1}\varphi_{\eps}}{\partial x^{k+1}}(0)=0.$

Under this reduction we will now show that $\varphi_{\eps}=id$. The
argument will be done with an infinite descent. We introduce the ideal
\[
I=\langle\eps_{0},\dots,\eps_{k-1}\rangle.\]
 With our preparation we know that $\varphi_{0}=id$ so we write \[
\varphi_{\varepsilon}\left(x\right):=x+\sum_{n\geq0}f_{n}\left(\varepsilon\right)x^{n}\]
 where each $f_{n}\in I$ and $f_{k+1}\equiv0$.

The conjugacy condition can be written as \begin{equation}
\begin{array}{l}
\left(1+a\left(\varepsilon\right)x^{k}\right)\left(\varphi_{\varepsilon}^{k+1}\left(x\right)+h_{k-1}\left(\varepsilon\right)\varphi_{\varepsilon}^{k-1}\left(x\right)+\cdots+h_{0}\left(\varepsilon\right)\right)\\
\qquad-\left(1+\tilde{a}\left(h\left(\varepsilon\right)\right)\varphi_{\varepsilon}^{k}\left(x\right)\right)\left(x^{k+1}+\varepsilon_{k-1}x^{k-1}+\cdots+\varepsilon_{0}\right)\varphi_{\varepsilon}'\left(x\right)=0.\end{array}\label{eq:formal_reccu}\end{equation}
 It is then clear that $h_{j}(\eps)\in I$. For the sake of simplicity
we simply write $h_{j}$ instead of $h_{j}(\eps)$. Let $g_{j}x^{j}$
be the term of degree $j$ in \eqref{eq:formal_reccu}. We will play
with the infinite set of equations $g_{j}=0$, $j\geq0$.

The equations $g_{j}=0$ with $0\leq j\leq k-1$ yield \[
h_{j}-\eps_{j}\in I^{2},\]
 since all other terms in the expression of $g_{j}$ belong to $I^{2}$.

The equation $g_{k+j}=0$ with $0\leq j\leq k$ yields $f_{j}\in I^{2}$
since the only terms of degree $1$ are $a(h_{j}-\eps_{j})+(k+1-j)f_{j}$
when $j<k$ and $af_{0}+f_{k}$ for $j=k$. Also, our hypothesis is
that $f_{k+1}\equiv0$. Looking at the linear terms in the equations
$g_{\ell}=0$ with $\ell>2k+1$ yields $f_{\ell-k}\in I^{2}$ since
the only linear terms in $g_{\ell}$ are $-(\ell-2k-1)[f_{\ell-k}+af_{\ell-2k}]$.

So we have that $f_{j}\in I^{2}$ for all $j$.

To show the conclusion we will shown by induction that, for any $n$,
$h_{j}-\eps_{j}\in I^{n}$ when $0\leq j\leq k-1$ and $f_{j}\in I^{n}$
whenever $j\geq0$. The conclusion is valid for $n=1,2$. We now suppose
that it is valid for $n$ and we show it for $n+1$.

To show that $h_{j}-\eps_{j}\in I^{n+1}$ for $0\leq j\leq k-1$ we
consider again the corresponding equations $g_{j}=0$, where the only
linear terms are $h_{j}-\eps_{j}$. Hence all other terms of the equation
belong to $I^{n+1}$ yielding $h_{j}-\eps_{j}\in I^{n+1}$.

For the same reason the equation $g_{k+j}=0$ with $0\leq j\leq k$
yields $f_{j}\in I^{n+1}$ and the equations $g_{\ell}=0$ with $\ell>2k+1$
yields $f_{\ell-k}\in I^{n+1}$. 
\end{proof}
\begin{defn}
\label{def:canonic_param} The parameter $\eps=(\eps_{0},\dots,\eps_{k-1})$
is called the \textbf{\textit{\emph{canonical}}} \textit{\emph{(multi-)}}\textbf{\textit{\emph{parameter}}}
of the family \eqref{unfold}. 
\end{defn}
\begin{cor}
\label{cor:facto_prepared} An orbital equivalence or a conjugacy
between two prepared families is the composition of a map which preserves
the canonical parameters with a map $(x,y,\eps)\mapsto(\tilde{x},y,\tilde{\eps})$
where \begin{equation}
\begin{cases}
\tilde{x}=\exp(2\pi i\, m/k)x\\
\tilde{\eps}_{j}=\exp(-2\pi i\, m(j-1)/k)\eps_{j} & \,\,\,\,\, j=0,\dots,k-1\end{cases}\label{sym_prep}\end{equation}
 for some $m\in\ww{Z}/k$. 
\end{cor}

\section{\label{sec:model-family}The model family}

We compare a prepared family of vector fields \eqref{unfold} with
multi-parameter $\eps=(\eps_{0},\dots,\eps_{k-1})\in\ww{C}^{k}$ to
the model family given by \eqref{model} with the same singular points,
and hence same parameters $\eps=(\eps_{0},\dots,\eps_{k-1})$ and
formal invariant $a(\eps)$ given by \eqref{eq_a}. The coefficients
of $Q_{\eps}$ are chosen so that $Z_{\eps}$ has the same eigenvalues
as the model at the singular points $(x_{j},0)$ for $j=0,\dots,k$
as \eqref{unfold}.

\bigskip{}

\subsection{\label{sub:generic_eps}The parameter space $\Sigma_{0}$}

We define \begin{equation}
||\eps||:=\max\left(|\eps_{k-1}|^{1/2},\dots,|\eps_{1}|^{1/k},|\eps_{0}|^{1/(k+1)}\right)\end{equation}
 The parameter space $\mathcal{W}=\left\{ \eps\,:\,||\eps||\leq\rho\right\} $
of $\Xi_{\eps}$ is stratified. The generic stratum $\Sigma_{0}$
is the set of $\eps$ for which the discriminant of $P_{\eps}$ does
not vanish: \begin{equation}
\Sigma_{0}:=\left\{ \eps\in\mathcal{W}\,:\,\mbox{disc}(P_{\eps})\neq0\right\} .\label{Sigma_0}\end{equation}
 The singular part, where the discriminant vanishes, is of codimension
one. Then, as soon as we define analytic and bounded functions on
$\Sigma_{0}$ they can be extended to $\mathcal{W}$ by the theorem
of removable singularities. For these reasons we limit ourselves to
parameter values in $\Sigma_{0}$.

\begin{lem}
\label{lem:mod_rac} All roots of $P_{\varepsilon}$ are contained
in a closed disk of radius at most $\sqrt{k}\left|\left|\varepsilon\right|\right|$. 
\end{lem}
\begin{proof}
Let $\eta:=\frac{1}{\sqrt{k}||\eps||}$. If $\left|x\right|>\frac{1}{\eta}$
then \begin{eqnarray*}
\left|\frac{P_{\varepsilon}\left(x\right)}{x^{k+1}}-1\right| & <\eta^{2}\left|\left|\varepsilon\right|\right|^{2}+\ldots+\eta^{k+1}\left|\left|\varepsilon\right|\right|^{k+1}\leq & 1,\end{eqnarray*}
 since each term is less than $\frac{1}{k}$. 
\end{proof}

\subsection{\label{first_int_model}First integral of the model family}

For $\eps\in\Sigma_{0}$ the model $X_{\eps}^{M}$ has a (multi-valued)
first integral \begin{equation}
H_{\eps}^{M}(x,y)=y\prod_{j=0}^{k}(x-x_{j})^{-\frac{1}{\nu_{j}}}\label{int_prem}\end{equation}
 where $\nu_{j}$ are defined by \eqref{eq_nu}. A first integral
of a vector field $X$ is a function $H$ such that $X\cdot H=0$
or, equivalently, which is constant on integral curves of $X$. If
$H$ is not constant then the connected components of the level sets
of $H$ coincide with the integral curves of $X$.

Below we will describe more precisely the foliation of $X_{\eps}^{M}$
over adequate sectors, but preliminary work is needed to describe
them. These sectors will correspond to domains over which $H_{\varepsilon}^{M}$
is univalued and takes all values in $\ww{C}$.

\subsection{The global and semi-local phase portrait of $P_{\varepsilon}\pp{x}$ }

The following trivial lemma will be used to define an equivalence
relation on the parameter space.

\begin{lem}
\label{P_1} The vector field \begin{equation}
\Xi_{\eps}:=P_{\varepsilon}\pp{x}\label{eq_Xi}\end{equation}
 is transformed into $\eta^{-(k+1)}\Xi_{\eps}$ under \begin{equation}
\begin{cases}
\eps=\left(\eps_{0},\dots,\eps_{k-1}\right)\mapsto\left(\eta^{k+1}\eps_{0},\eta^{k}\eps_{1},\dots,\eta^{2}\eps_{k-1}\right)\\
x\mapsto x/\eta,\end{cases}\label{rescaling}\end{equation}
 where $\eta\in\mathbb{R}_{>0}$. Hence the bifurcation diagram for
the phase portrait of $\Xi_{\eps}$ has a conic structure and is completely
determined on the surface $||\eps||=\rho$. 
\end{lem}
\begin{defn}
We define the following equivalence relation on the set of $\eps$:
\begin{equation}
\eps=\left(\eps_{0},\dots,\eps_{k-1}\right)\simeq\eps'=\left(\eps_{0}',\dots,\eps_{k-1}'\right)\Longleftrightarrow\exists\eta\in\mathbb{R}_{>0}\,:\:\eps_{j}'=\eta^{k+1-j}\eps_{j}.\label{equiv_eps}\end{equation}

\end{defn}
The global phase portrait of $\Xi_{\varepsilon}$ is studied by Douady
and Sentenac in \cite{DS}. They show how the attracting and repelling
separatrices of the saddle point at infinity separates the phase plane
in simply connected regions. Among the different $\Xi_{\varepsilon}$
they make a special discussion of the generic $P_{\varepsilon}\pp{x}$,
which have the property that there is no homoclinic trajectory, \emph{i.e.}
no connection between an attracting and a repelling separatrix at
infinity.

\begin{defn}
The vector field $\Xi_{\varepsilon}$ is generic in the sense of Douady
and Sentenac if all singular points are distinct and there are no
homoclinic trajectories. 
\end{defn}
Douady and Sentenac show that the eigenvalues of the singular points
of a generic vector field all have a nonzero real part and then that
the singular points are nodes or foci.

A homoclinic trajectory $\gamma$ goes in finite real time $T$ from
infinity to infinity since infinity is either a regular point ($k=1$)
or a pole ($k>1$). The close loop $\gamma$ on $\mathbb{C}\mathbb{P}^{1}$
necessarily contains some singular points $x_{j_{1}},\dots,x_{j_{s}}$
in its \lq\lq interior'' The value of $T$ can be calculated by
the residue theorem: \[
T=\int_{\gamma}\frac{dx}{P_{\eps}(x)}=2\pi i\sum_{\ell=1}^{s}\frac{1}{P_{\eps}'(x_{j_{\ell}})}.\]
 (Even if the \lq\lq interior'' of $\gamma$ is not well defined
$T$ is well defined since $\sum_{j=0}^{k}\frac{1}{P_{\eps}'(x_{j})}=0$).
Moreover $T$ cannot vanish since $\gamma$ is non contractible. We
will recall below a lower bound for $T$ calculated in \cite{DS}.
Douady and Sentenac show that the generic $\Xi_{\varepsilon}$ are
dense and also that, given any $\Xi_{\varepsilon}$ with $\eps\in\Sigma_{0}$,
there exists an angle $\theta$ such that $\exp(i\theta)\Xi_{\varepsilon}$
is generic.

We will derive below an adaptation of their result coming from the
fact that we are only interested in $\Xi_{\varepsilon}$ over a disk
$r\mathbb{D}$.

\begin{figure}
\subfigure[$\eps=0$]{\includegraphics[width=5cm]{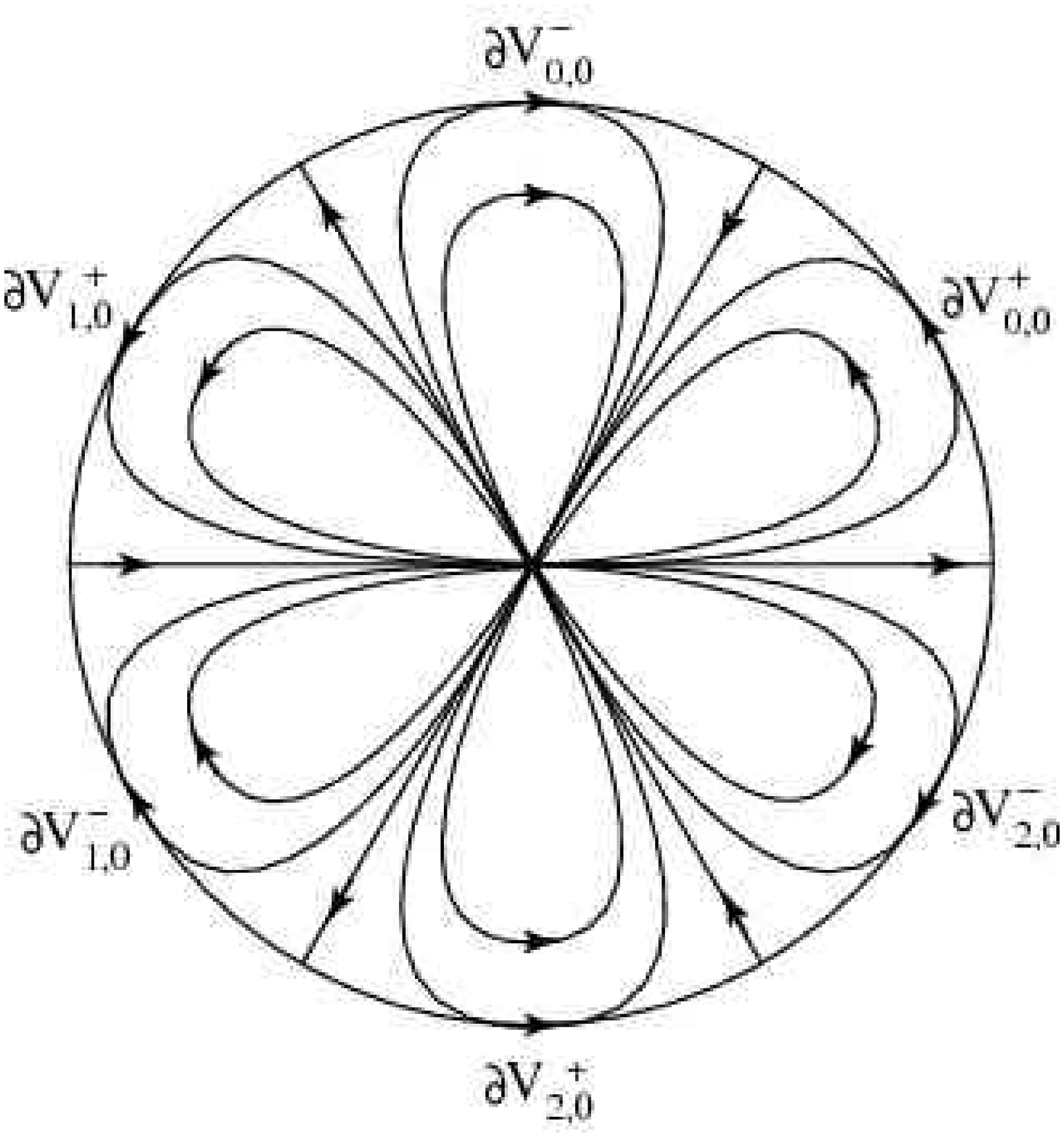}} \hfill{}\subfigure[general
$\eps$]{\includegraphics[width=5cm]{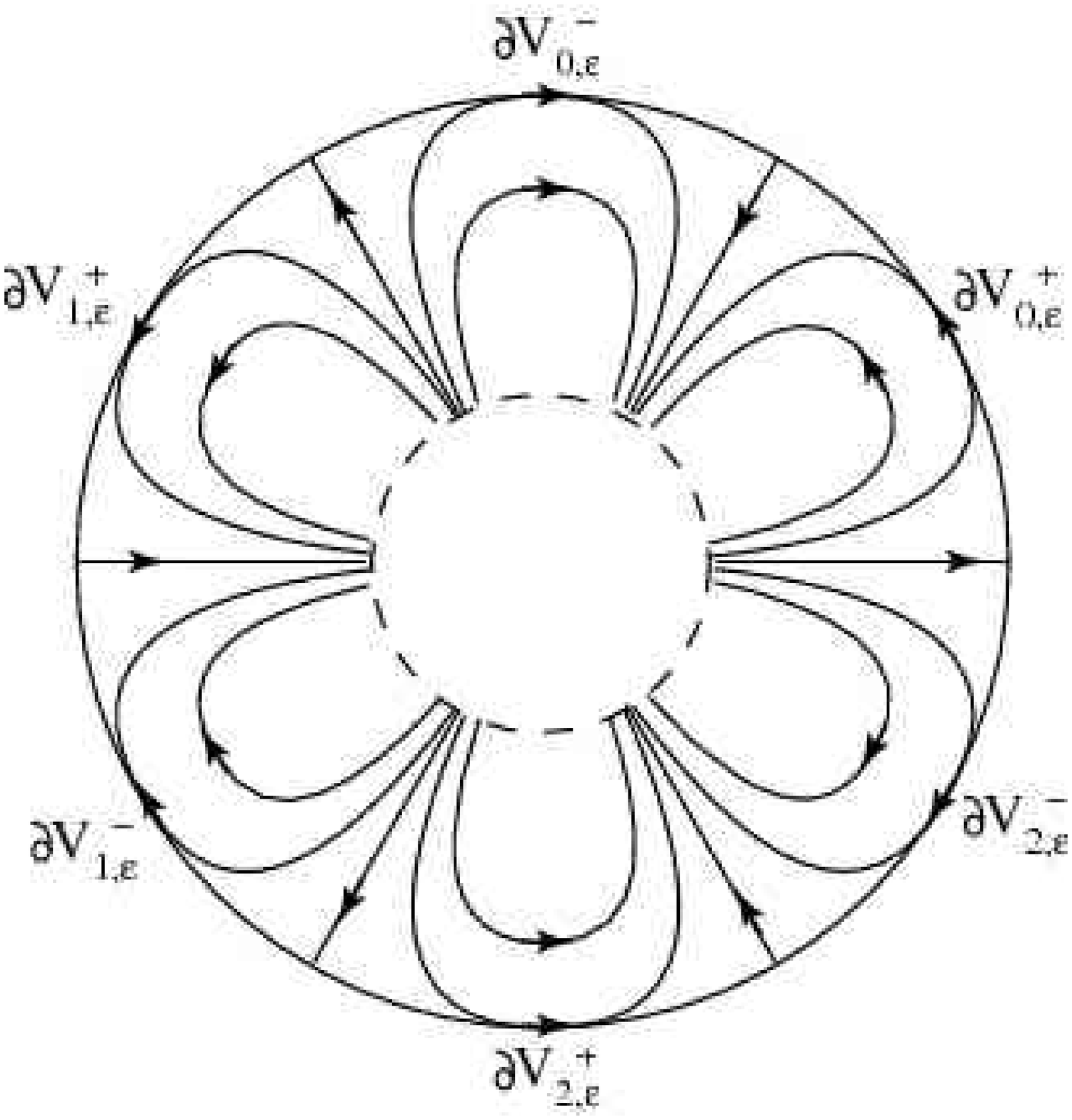}}

\caption{\label{fig_fleur} The dynamics of $\Xi_{\eps}$ near $|x|=r$.} 
\end{figure}

\begin{lem}
\label{lem:mino_deriv} Let $K=2^{k}-1$ be the number of partitions
of $\{x_{0},\dots,x_{k}\}$ into the union of two disjoint non empty
subsets. Let \[
\delta=\frac{\pi}{16K+2}\]
 and \begin{equation}
\Xi_{\eps}(\theta):=\exp(i\theta)\Xi_{\eps}.\label{eq:theta_theta}\end{equation}
 For any $\eps\in\Sigma_{0}$ (\emph{i.e.} such that all roots of
$P_{\eps}$ are distinct) there exists $\theta=\theta(\eps)\in(-\pi/4,\pi/4)$
such that $\Xi_{\eps}(\theta)$ is generic in the sense of Douady
and Sentenac. More precisely, for any partition $\{x_{0},\dots,x_{k}\}=I_{1}\cup I_{2}$
with $I_{1},I_{2}\neq\emptyset$ \begin{eqnarray*}
\left|\arg\left(\exp(i\theta)\sum_{x_{j}\in I_{1}}\frac{1}{P'(x_{j})}\right)\right|-\frac{\pi}{2} & \notin & (-\delta,\delta)\,.\end{eqnarray*}
 Moreover $\theta(\eps)$ can be chosen constant on a neighborhood
of a given $\tilde{\eps}$ and can also be chosen constant on the
equivalence class of $\eps$ under \eqref{equiv_eps}. It is possible
to cover $\Sigma_{0}$ with $m=4K-1$ connected open sets $W_{i}$
on which $\theta(\eps)$ can be chosen constant. 
\end{lem}
\begin{proof}
Let $\theta_{\ell}=\frac{\pi\ell}{8K+1}$ for $\ell\in\{-(2K-1),\dots,0,\dots,2K-1\}$.
We consider the set $J$ of partitions $\{x_{0},\dots,x_{k}\}=I_{\ell_{1}}\cup I_{\ell_{2}}$
with disjoint $I_{\ell_{1}},I_{\ell_{2}}\neq\emptyset$.

We let \[
W_{\ell}=\left\{ \eps\;:\;\left|\arg\sum_{j\in I_{\ell_{1}}}\frac{1}{P_{\eps}'(x_{j})}\right|+\theta_{\ell}-\frac{\pi}{2}\notin[-\delta,\delta]\,\,,\,\,(I_{\ell_{1}},I_{\ell_{2}})\in J\right\} .\]
 We need to show that $\{W_{\ell}\}_{\ell\in\{-(2K-1),\dots,0,\dots,2K-1\}}$
is an open covering of $\Sigma_{0}$ and that the $W_{\ell}$ are
connected.

For this purpose we suppose that $\arg\sum_{j\in I_{\ell_{1}}}\frac{1}{P_{\eps}'(x_{j})}\in[-\pi,\pi]$.
Let $M_{I_{\ell_{1}}}=\left|\arg\sum_{j\in I_{\ell_{1}}}\frac{1}{P_{\eps}'(x_{j})}\right|-\frac{\pi}{2}\in[-\frac{\pi}{2},\frac{\pi}{2}]$.
We divide $[-\frac{\pi}{2},\frac{\pi}{2}]$ in $8K+1$ equal closed
intervals. $4K+1$ of these sub-intervals cover $[-\frac{\pi}{4},\frac{\pi}{4}]$.
Among this group of $4K+1$ intervals there is at least one group
of three consecutive intervals whose union contains no $M_{I_{\ell_{1}}}$
in its interior. We apply one of the rotations $\theta_{\ell}$ to
send this group of intervals to the center interval, namely $[-3\delta,3\delta]$.
The $m=4K-1$ open sets correspond to the $4K-1$ ways to choose three
consecutive intervals (from the $4K+1$ intervals) covering $[-\frac{\pi}{4},\frac{\pi}{4}]$.
The fact that the $W_{\ell}$ are connected comes from the fact that
$\Sigma_{0}$ does not separate the $\eps$-space. 
\end{proof}
\begin{defn}
Let $\rho>0$ be given. 
\begin{enumerate}
\item We consider a neighborhood $\mathcal{W}=\left\{ \varepsilon\,:\,||\eps||\leq\rho\right\} $
of $\eps=0$. An open sector $W$ of $\mathcal{W}\setminus\{0\}$
is an \textit{adequate sector} if it is a union of equivalence classes
of \eqref{equiv_eps} inside $\mathcal{W}$. 
\item Let $\Sigma_{0}$ defined in \eqref{Sigma_0}. An open covering $\left\{ W_{i}\right\} _{i\in I}$
of $\Sigma_{0}$, where $W_{i}\subset\Sigma_{0}$, is an \textit{adequate
covering} of $\Sigma_{0}$ if each $W_{i}$ is an adequate sector. 
\item Given $\eps\in\Sigma_{0}$ we associate to it an angle $\theta(\eps)$.
The angle $\theta(\eps)$ is \textit{adequate} if it satisfies Lemma~\ref{lem:mino_deriv}
and if it can be chosen constant on the equivalence class of $\eps$. 
\end{enumerate}
\end{defn}
The 1-dimensional vector field $\Xi_{0}=P_{0}\frac{\partial}{\partial x}$
is given on $\left\{ x\,:\,|x|\leq r\right\} $ in Figure~\ref{fig_fleur}(a).
For $\eps$ sufficiently small the phase portrait near $\left\{ |x|=r\right\} $
is similar to that of $\Xi_{0}$ (Figure~\ref{fig_fleur}(b)). In
particular the boundary $\left\{ |x|=r\right\} $ has $k$ sub-sectors
on which the vector field goes inwards and $k$ sub-sectors on which
it goes outwards. %

\bigskip{}

To study properly the vector field $\Xi_{\eps}$ it is useful to change
the $x$-coordinate to the complex time-coordinate (the generalized
Fatou coordinate). This is the point of view described by Douady and
Sentenac in \cite{DS}. The following lemmas summarize some properties
described in \cite{DS}, so most of them will be given without proof.
They are also equivalent to some properties described by Oudkerk in
\cite{O}.

\begin{figure}
\subfigure[$z(r\mathbb{D})$ (the
holes are removed)]{\includegraphics[width=7cm]{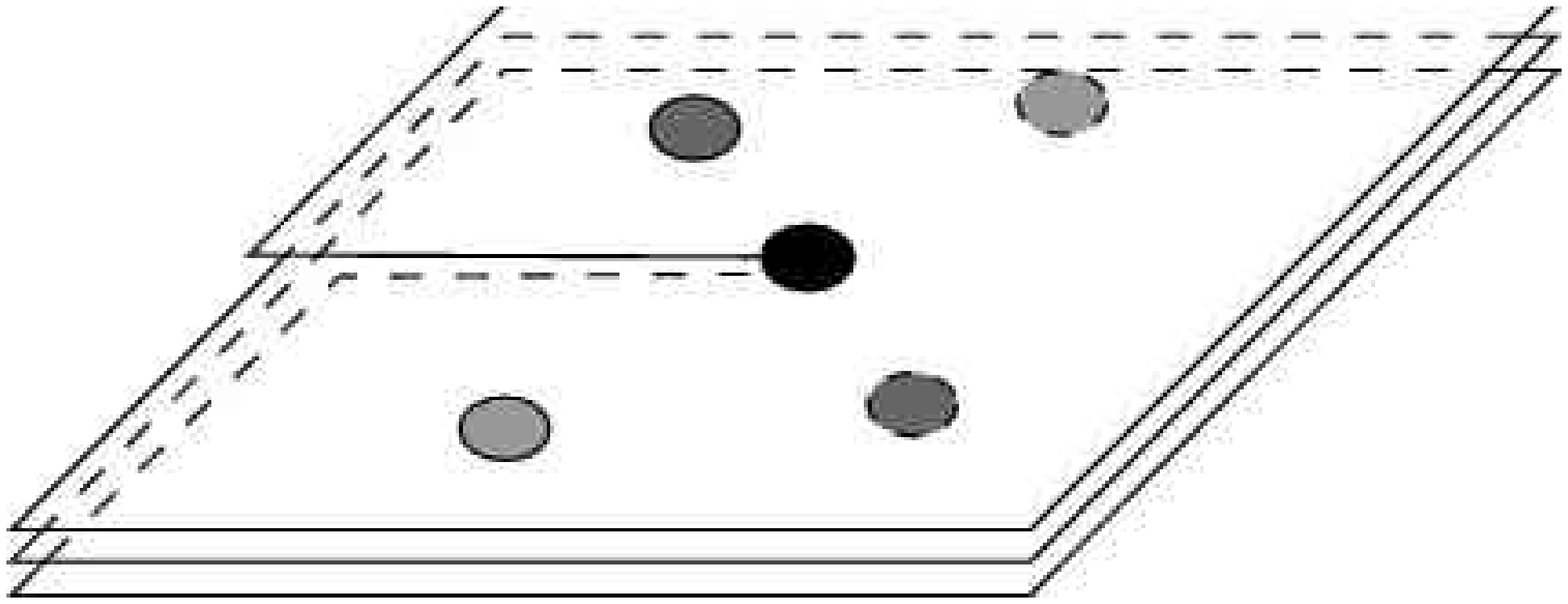}} \hfill{}\subfigure[A horizontal strip in
$z$-space]{\includegraphics[width=7cm]{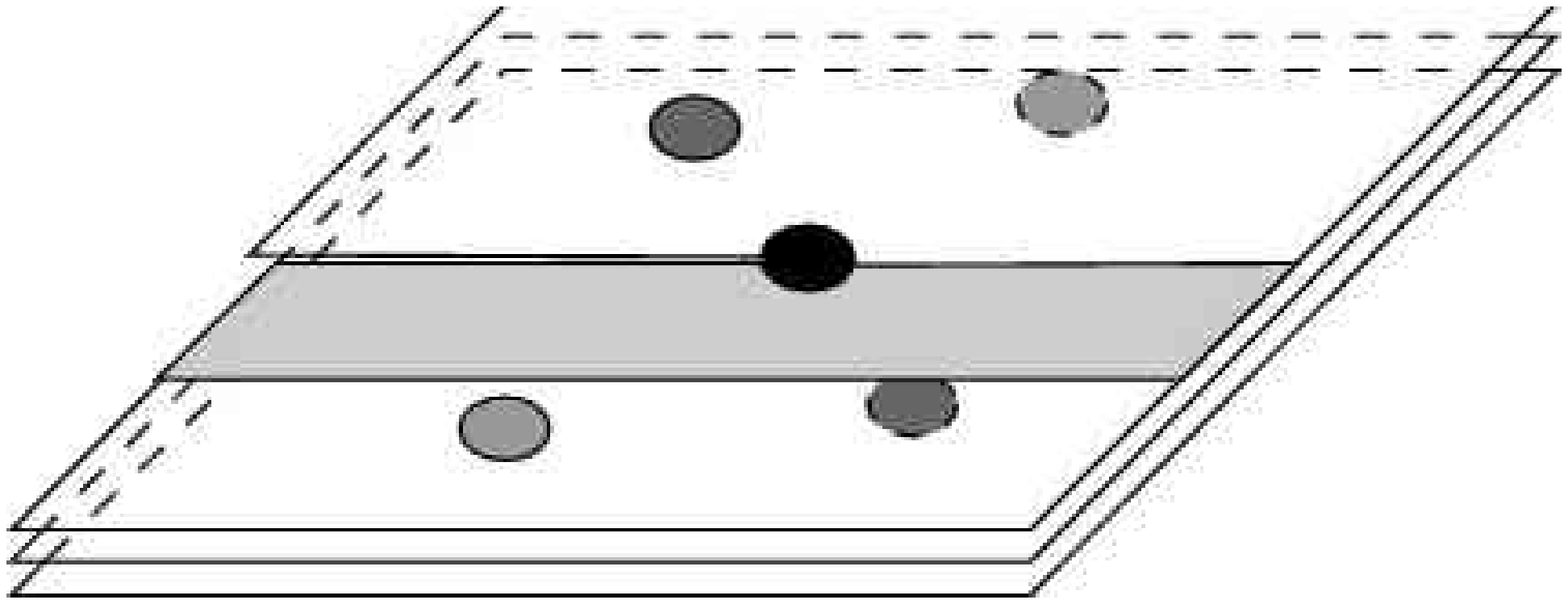}}

\caption{The image of $r\mathbb{D}$ in $z$-space. There are ramifications
at each hole.}

\label{fig2} 
\end{figure}

\begin{lem}
\label{P_2} The change of coordinate \begin{equation}
z=z(x)=\int_{\infty}^{x}\frac{dx}{P_{\eps}(x)}\label{coord_z}\end{equation}
 is a multivalued function defined on $\mathbb{C}\mathbb{P}\setminus\{x_{0},\dots,x_{k}\}$.
It is a $k$-sheeted covering $\mathcal{S}$ of a neighborhood of
$\infty$. The image of the circle $r\mathbb{S}^{1}$ is a $k$-covering
of a closed curve which is approximately a circle of radius $\frac{1}{kr^{k}}$:
see Figure~\ref{fig2}(a). We call $B_{0}$ the interior of the projection
of this curve. 
\end{lem}
In the sequel we consider values of the parameter $\varepsilon\neq0$
for which the multivalueness of $z$ is of logarithmic type at each
branch point. Hence the ball $B_{0}$ will have infinitely many images
distributed in the $z$-plane. A more precise describtion is done
further below.

\begin{figure}
\includegraphics[width=60mm]{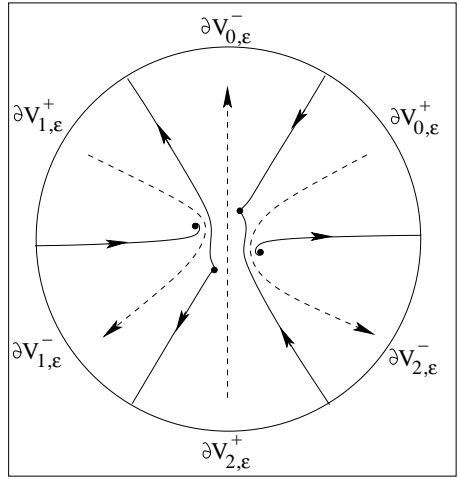}

\caption{\label{fig:fig_sigma}A separating graph $\Gamma$ and associated
correspondence $\sigma$ for $k=3$.} 
\end{figure}

\begin{defn}
The union $\Gamma$ of the separatrices of a generic vector field
in the sense of Douady and Sentenac is called the \textbf{\textit{\emph{separating
graph}}} which divides $\mathbb{C}\setminus\Gamma$ in $k$ simply
connected components. 
\end{defn}
A generic separating graph divides the disk $\left\{ |x|\leq r\right\} $
into $k$ simply-connected domains. At this level of the construction
we will only consider the corresponding $2k$ segments on the circle
$\left\{ \left|x\right|=r\right\} $, which we conveniently denote
by $\partial V_{j,\varepsilon}^{\pm}$ for $j\in\left\{ 0,1,\ldots,k-1\right\} $
(see Figures~\ref{fig_fleur} and \ref{fig:fig_sigma}). A reason
for this notation will be found in Section~\ref{sub:Squid-sectors}
where we define the squid sectors $\sect{j,\varepsilon}{\pm}$.

\begin{lem}
\label{lemma_DS} \cite{DS} We consider a generic vector field in
the sense of Douady and Sentenac. Each connected component of $\mathbb{C}\setminus\Gamma$,
where $\Gamma$ is the separating graph, intersects $r\mathbb{S}^{1}$
in exactly one sector $\partial V_{j,\varepsilon}^{+}$ and one sector
$\partial V_{\ell,\varepsilon}^{-}$ (see Figures~\ref{fig_fleur}
and \ref{fig:fig_sigma}). This yields a correspondence \begin{equation}
\sigma:\{0,\dots,k-1\}\rightarrow\{0,\dots,k-1\},\qquad\quad j\mapsto\ell\label{def_sigma}\end{equation}
 between the sectors $\partial V_{j,\varepsilon}^{+}$ and $\partial V_{\ell,\varepsilon}^{-}$. 
\end{lem}

\begin{proof}
For $\eps\in\Sigma_{0}$ we have that $\frac{dx}{P_{\eps}(x)}=\sum_{j=0}^{k}\frac{dx}{P_{\eps}'(x_{j})(x-x_{j})}$.
Since the logarithm function is multivalued, this yields other images
of the circle $r\mathbb{D}$ as drawn in Figure~\ref{fig2}(a)).
The interior of these curves (which we will also call balls) correspond
to images of the exterior of $r\mathbb{D}$. A straight line of slope
$\theta$ joining two such balls corresponds to a trajectory in real
time of $\Xi_{\eps}(\theta)$ joining a boundary sector of $r\mathbb{S}^{1}$
to a boundary sector of $r\mathbb{S}^{1}$. 
\end{proof}
\begin{lem}
\label{P_3} We consider one of the angles $\theta_{\ell}$ of the
proof of Lemma~\ref{lem:mino_deriv} and the open sector $W_{\ell}$
of values of $\eps$ for which this angle is adequate. There exists
$\rho>0$ sufficiently small so that for any $||\eps||\leq\rho$ all
trajectories of $\Xi_{\eps}(\theta)$ in real time starting on $|x|=r$
end in a singular point $x_{j}$. 
\end{lem}
\begin{proof}
We change to the $z$ coordinate. Each of these sectors corresponds
to either the upper half-circle or the lower half-circle of the $k$-sheeted
covering of the boundary of $B_{0}$. Following the trajectory of
$\Xi_{\eps}(\theta)$ in real time is the same as following a trajectory
of $\Xi_{\eps}$ in time $e^{i\theta}\mathbb{R}$, hence following
a line of slope $\theta$ in $z$ space.

There exists a trajectory $\gamma_{j}$ of $\Xi_{\eps}(\theta)$ in
$r\mathbb{D}$ in time $\exp(i\theta')\mathbb{R}$ for some $\theta'$
starting on each $\partial V_{j,\varepsilon}^{+}$ and ending in $\partial V_{\sigma(j),\varepsilon}^{-}$
and not crossing the separating graph: we can consider this trajectory
oriented from $\partial V_{j,\varepsilon}^{+}$ to $\partial V_{\sigma(j),\varepsilon}^{-}$.
Let $I_{1}$ be the set of singularities to the right of $\gamma_{j}$
and $M(\gamma_{j})$ given by \begin{equation}
M(\gamma_{j})=\exp(-i\theta)\sum_{x_{\ell}\in I_{1}}\frac{1}{P_{\eps}'(x_{\ell})}.\label{eq_mu_j}\end{equation}
 If we consider the time function $z$ defined for $\partial V_{j,\varepsilon}^{+}$,
then it has a period $2\pi iM(\gamma_{j})$, so we should visualize
the $z$ plane with holes which would be the translates of $B_{0}$
by $2\pi iM(\gamma_{j})$. The same $M(\gamma_{j})$ is also valid
for initial conditions on $\partial V_{\sigma(j),\varepsilon}^{-}$.
The trajectory $\gamma_{j}$ can be visualized in $z$-space as a
line from a hole to an adjacent hole. Hence its length is of the order
of $M(\gamma_{j})-\frac{2}{kr^{k}}$ and corresponds to the modulus
of the time to travel along $\gamma_{j}$. We want to show that this
quantity remains large when $\eps$ is small. Indeed the trajectory
$\gamma_{j}$ separates the singular points inside $r\mathbb{D}$
into two non empty sets. As all $x_{j}$ lie in the disk of radius
$\sqrt{k}||\eps||$, then $\gamma_{j}$ has points inside that disk.
Hence the time to travel along $\gamma_{j}$ is bounded in modulus
by twice the modulus of the time to travel from $r\mathbb{S}^{1}$
to $\sqrt{k}||\eps||\mathbb{S}^{1}$.

Instead of evaluating this time we will use the results of \cite{DS}.
Indeed a slanted line of slope $\theta'$ joining two holes in the
$z$-coordinate for $\Xi_{\eps}(\theta)$ corresponds to a horizontal
line joining two holes in $\Xi_{\eps}(\theta-\theta')$. Such a line
(if it joins the center of the two holes) is a homoclinic trajectory.
So we need to estimate the traveling time $M'(\gamma_{j})$ of a homoclinic
trajectory of $\Xi_{\eps}(\theta-\theta')$. In \cite{DS} (Corollary
I.2.2.1) we find the following estimate \begin{equation}
|M'(\gamma_{j})|>\frac{1}{2^{k(k+4)/2}\max(|x_{0}|,\dots,|x_{k}|)}>\frac{1}{2^{k(k+4)/2}\sqrt{k}||\eps||}\label{eq:estim_bandwidth}\end{equation}
 since all $x_{j}$ satisfy $|x_{j}|<\sqrt{k}||\eps||$. Moreover
it is clear that $|M(\gamma_{j})|$ is approximately $|M'(\gamma_{j})|$
minus twice the time to travel from $\infty$ to $|x|=r$. Hence $|M(\gamma_{j})|\sim|M'(\gamma_{j})|-\frac{2}{k|r|^{k}}$.

To finish the proof we know that $\arg M(\gamma_{j})\notin[-\delta,\delta]$.
Hence the horizontal line of $\Xi_{\eps}(\theta)$ will not encounter
any hole if $|M(\gamma_{j})|\sin\delta>\frac{2}{kr^{k}}$. From the
estimate above, this is clearly satisfied as soon as $||\eps||$ is
sufficiently small. 
\end{proof}
\begin{lem}
\label{P_4} If $\eps\in\Sigma_{0}$ and $\theta(\eps)$ is adequate,
a horizontal strip as in Figure~\ref{fig2}(b) will start in a singular
point $x_{n}$ and end in a singular point $x_{s}$ such that $Re\left(e^{i\theta}P_{\varepsilon}'\left(x_{n}\right)\right)>0$
and $Re\left(e^{i\theta}P_{\varepsilon}'\left(x_{s}\right)\right)<0$.
The same holds for an infinite strip with parallel slanted ends as
in Figure~\ref{figsect_eps}. 
\end{lem}
\begin{thm}
\label{covering} There exists a finite open covering $\left\{ W_{i}\right\} _{i\in I}$
of $\Sigma_{0}$, where the $W_{i}$ are adequate sectors and, for
each $W_{i}$, there exists a constant adequate angle $\theta_{i}(\eps)=:\theta_{i}$
such that the conclusion of Lemma~\ref{P_3} holds. 
\end{thm}
\begin{proof}
The proof is immediate if one works in the time coordinate \eqref{coord_z}. 
\end{proof}
\begin{defn}
\label{good_covering} An open covering $\left\{ W_{i}\right\} _{i\in I}$
of $\Sigma_{0}$ as in Theorem~\ref{covering} is called a \textbf{\textit{\emph{good
covering}}} of $\Sigma_{0}$ and the angles $\theta_{i}$ are called
\textbf{\textit{\emph{good angles}}}. Each $W_{i}$ with this property
is called a \textbf{\textit{\emph{good sector}}}. 
\end{defn}
\begin{rem}
To give a good covering of $\Sigma_{0}$ it is sufficient to give
a good covering of $\Sigma_{0}\cap\left\{ \eps\,:\,||\eps||=\rho\right\} $.
Then the good covering is given by the equivalence classes of the
elements of the good covering of $\Sigma_{0}\cap\left\{ \eps\,:\,||\eps||=\rho\right\} $.

In the case $k=1$ there exists a good covering with only two sectors,
for instance $\arg(\eps)\in(-\eta,\pi+\eta)$ and $\arg(\eps)\in(-\pi-\eta,\eta)$,
with $\eta\in(0,\pi)$. The smaller $\eta$, the less spiraling in
the drawing of the sectors. 
\end{rem}

\subsection{\label{sub:Squid-sectors}Squid sectors}

The first integral of the model family is ramified in the $x$-variable.
For each value of $\eps$ in a small neighborhood of the origin we
will define $2k$ sectors in $x$-space (called adapted sectors),
above which the first integral \eqref{int_prem} is univalued and
there is a one-to-one correspondence between the set of  level curves
of the first integral and $\mathbb{C}$.

\begin{figure}
\subfigure[sector in
$x$-space]{\includegraphics[width=5cm]{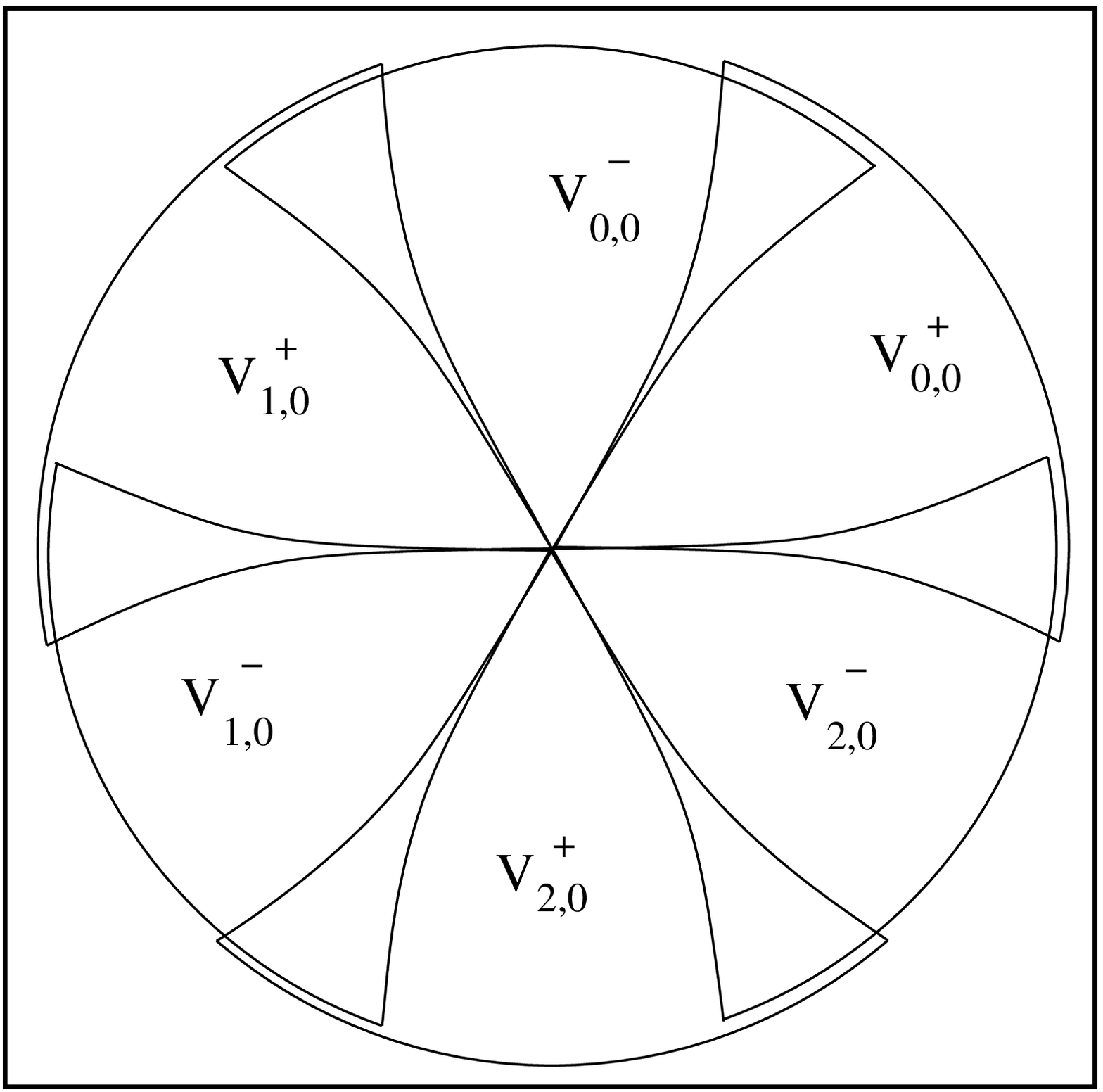}} \hfill{}\subfigure[sector in
$z$-space]{\includegraphics[width=8cm]{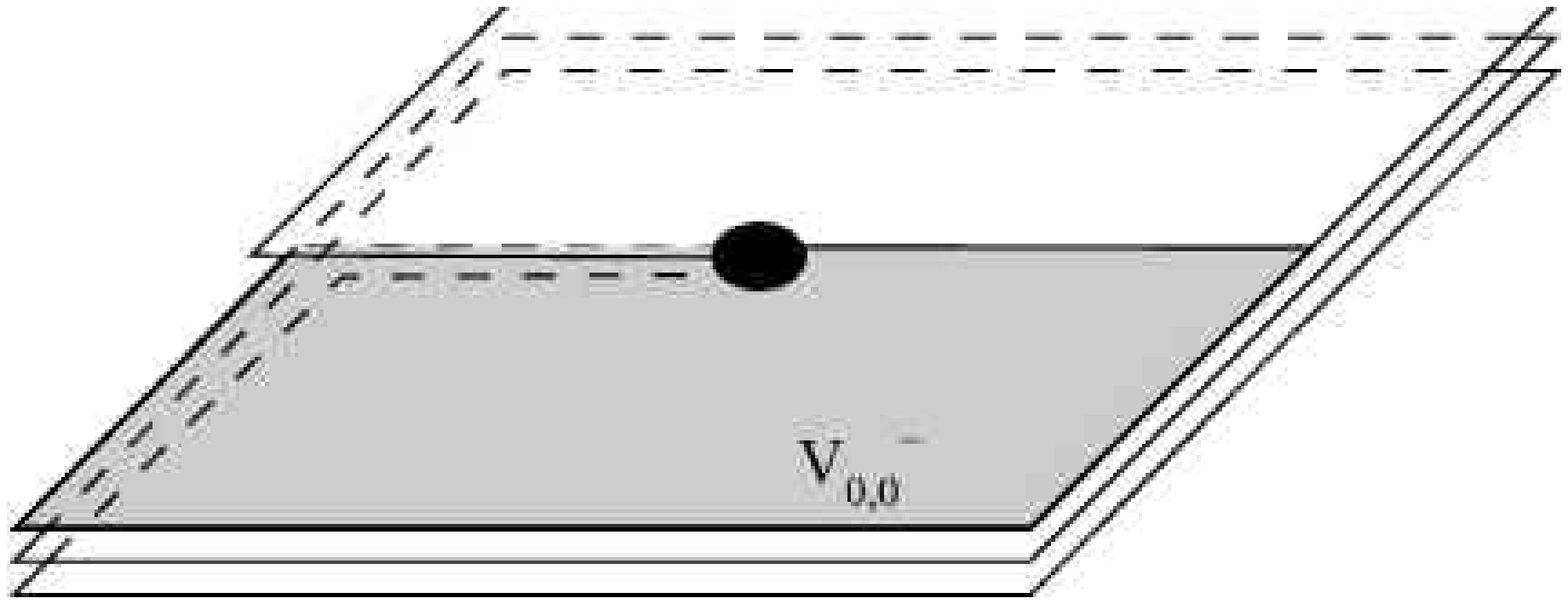}}

\caption{The sectors $\sect{j,0}{\pm}$ for $\eps=0$ }

\label{figsect_0} 
\end{figure}

For $\eps=0$ we define a unique set of $2k$ sectors $\sect{j,0}{\pm}$,
$j=0,\dots,k-1$ (Figure~\ref{figsect_0}). These sectors define
sectors $\partial V_{j,0}^{\pm}$ on the boundary $r\ww{S}^{1}=\partial\left(r\ww{D}\right)$
: the sectors defined for $\eps\neq0$ will be associated to the same
boundary sectors $\partial V_{j,\eps}^{\pm}$ of $r\ww{S}^{1}$. For
a given $\eps\neq0$ the $2k$ sectors may not be uniquely defined
and for each $\eps\neq0$ belonging to several $W_{i}$ there will
be several non-equivalent sets of $2k$ adapted sectors $\sect{j,\varepsilon}{\pm}$
for $j=0,\dots,k-1$, with same boundary sectors $\partial V_{j,\eps}^{\pm}$
(one for each $W_{i}$). In particular, in the generic case, each
sector will be adherent to two singular points and non-equivalent
sectors may be adherent to different pairs of singular points. However
we will limit ourselves to definitions of sectors valid on equivalence
classes of $\eps$ (under the equivalence relation \eqref{equiv_eps}).
When $\eps\to0$ inside an equivalence class, any set of $2k$ sectors
will have the same limit: $\sect{j,\varepsilon}{\pm}\to\sect{j,0}{\pm}$.

\begin{figure}
\includegraphics[width=8cm]{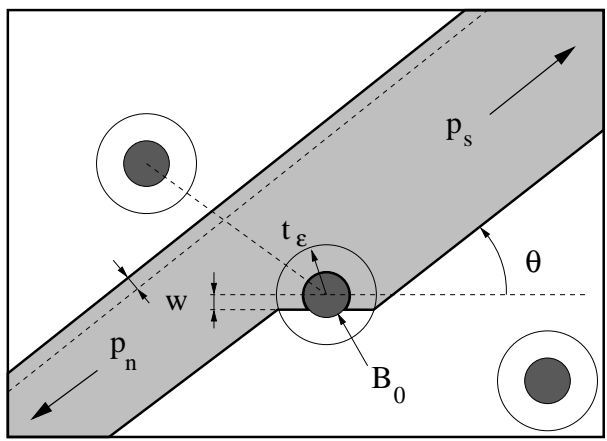}

\caption{A squid sector in $z$-coordinate for $\eps$ general}

\label{fig:z_squid} 
\end{figure}

\begin{figure}
\subfigure{\includegraphics[width=6cm]{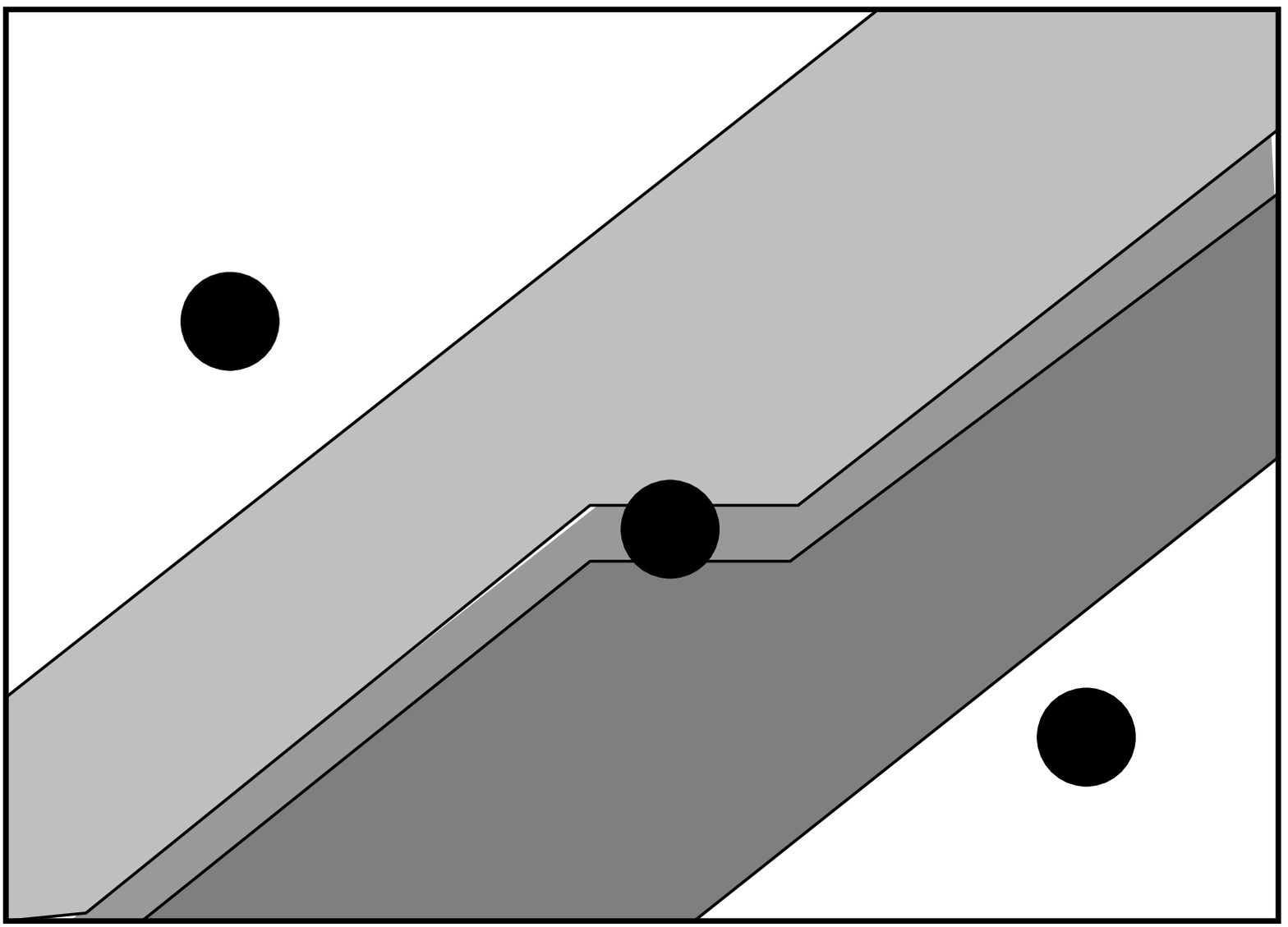}}\qquad{}\qquad{}\qquad{}\subfigure{\includegraphics[width=6cm]{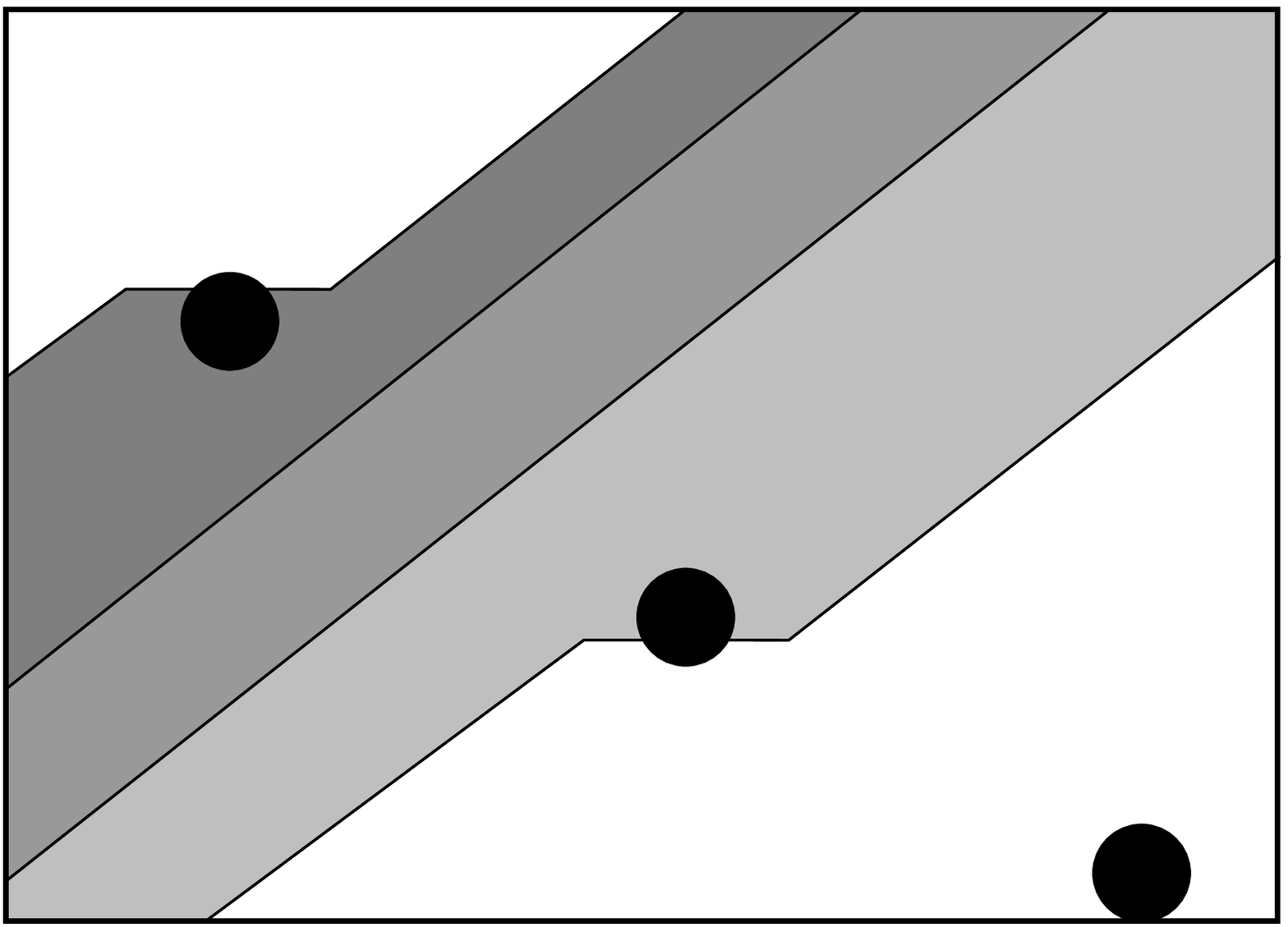}}\\
 \subfigure{\includegraphics[width=6cm]{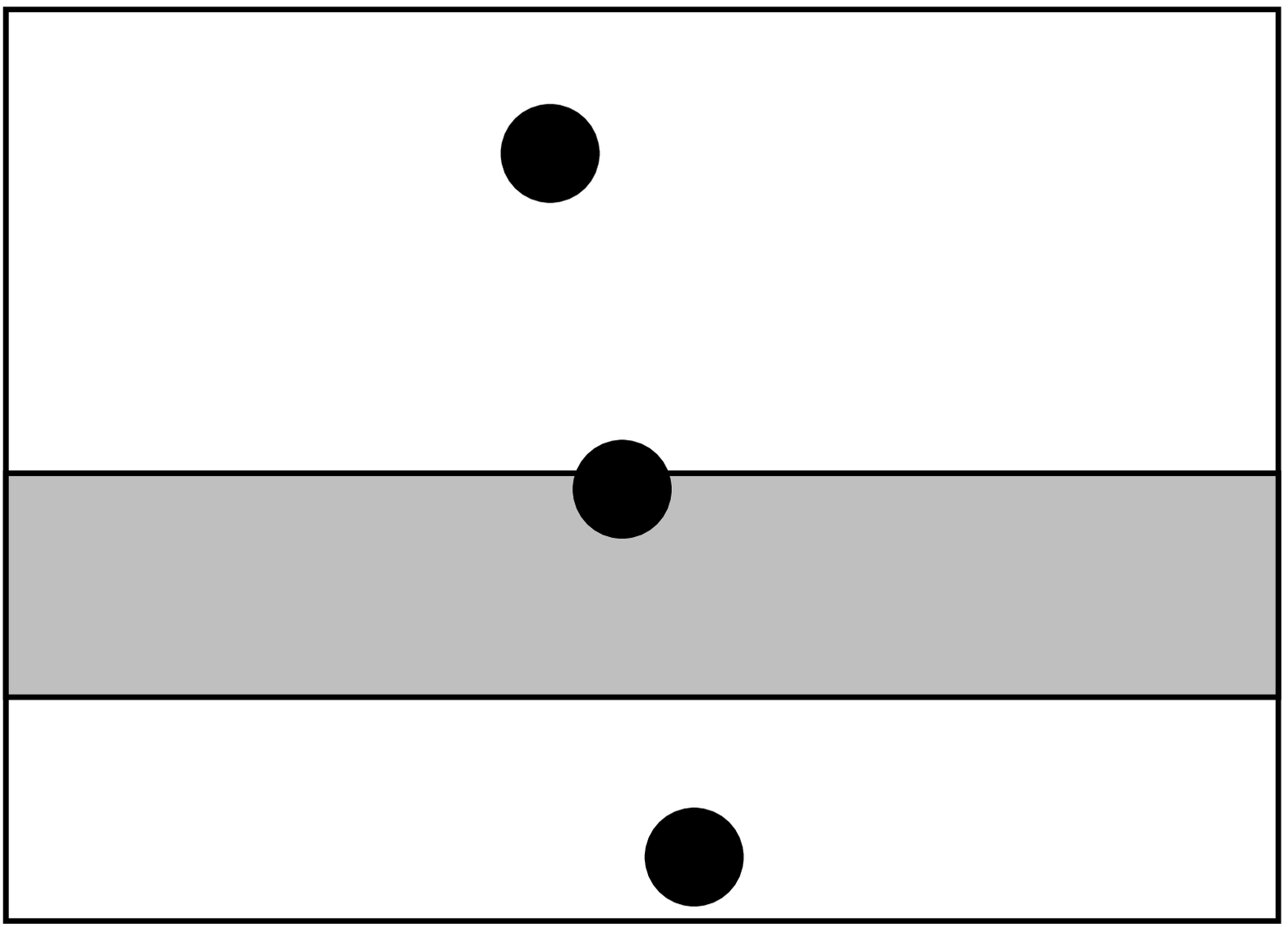}}\qquad{}\qquad{}\qquad{}\subfigure{\includegraphics[width=6cm]{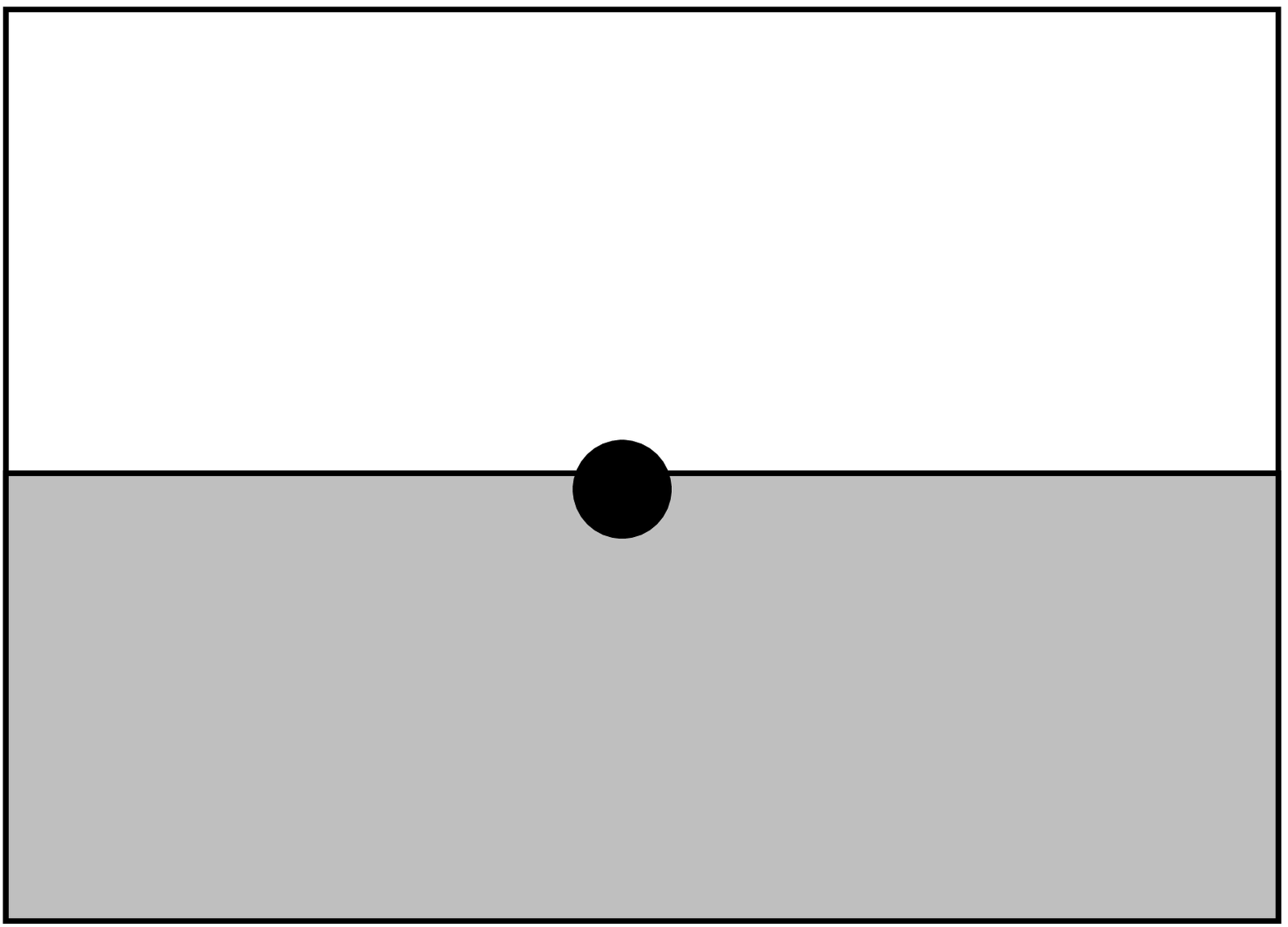}}

\caption{Different squid sectors in $z$-space and their intersections. }

\label{figsect_eps} 
\end{figure}

\begin{figure}
\subfigure{\includegraphics[width=6cm]{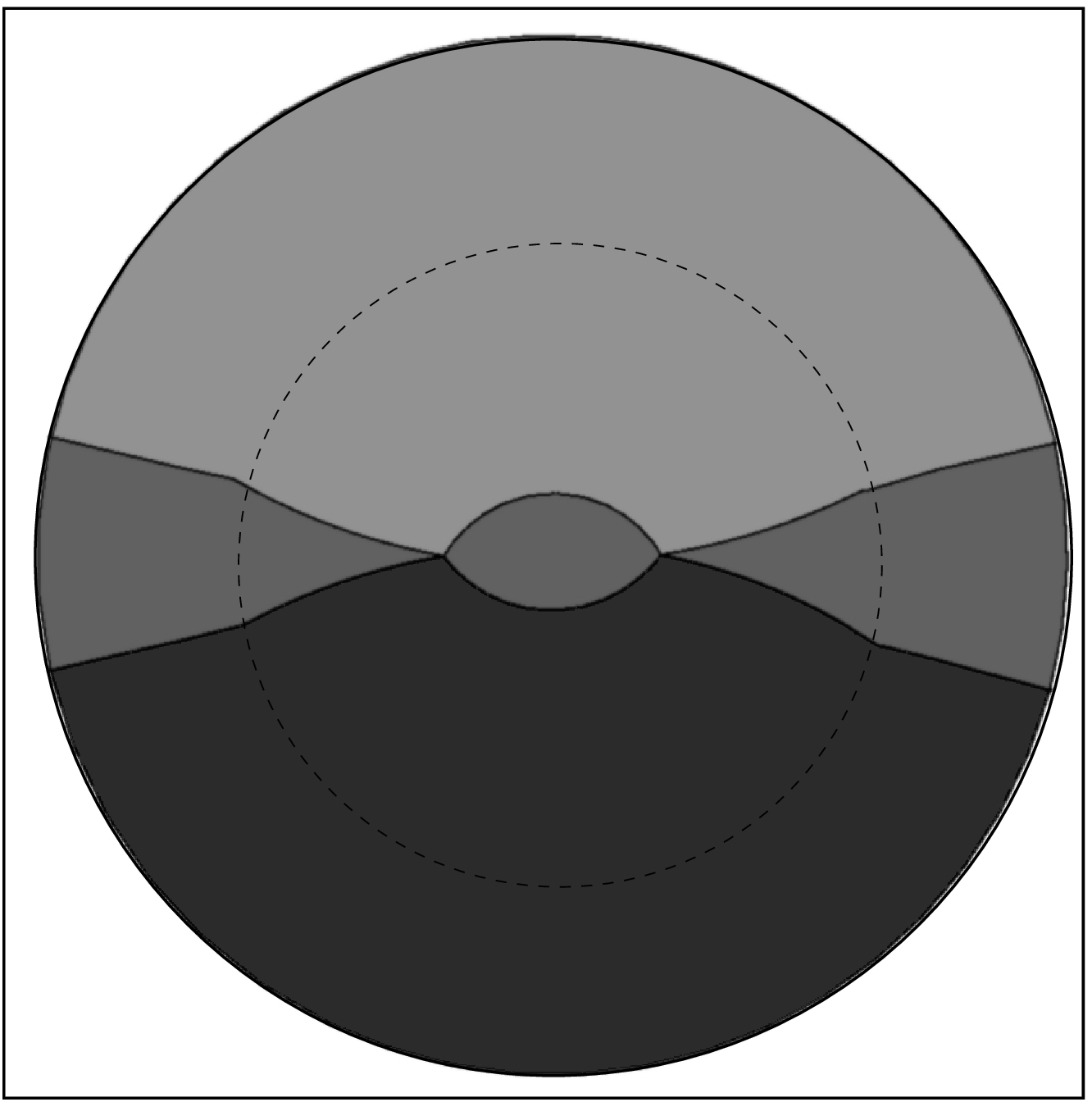}}\hfill{}\subfigure{\includegraphics[width=6cm]{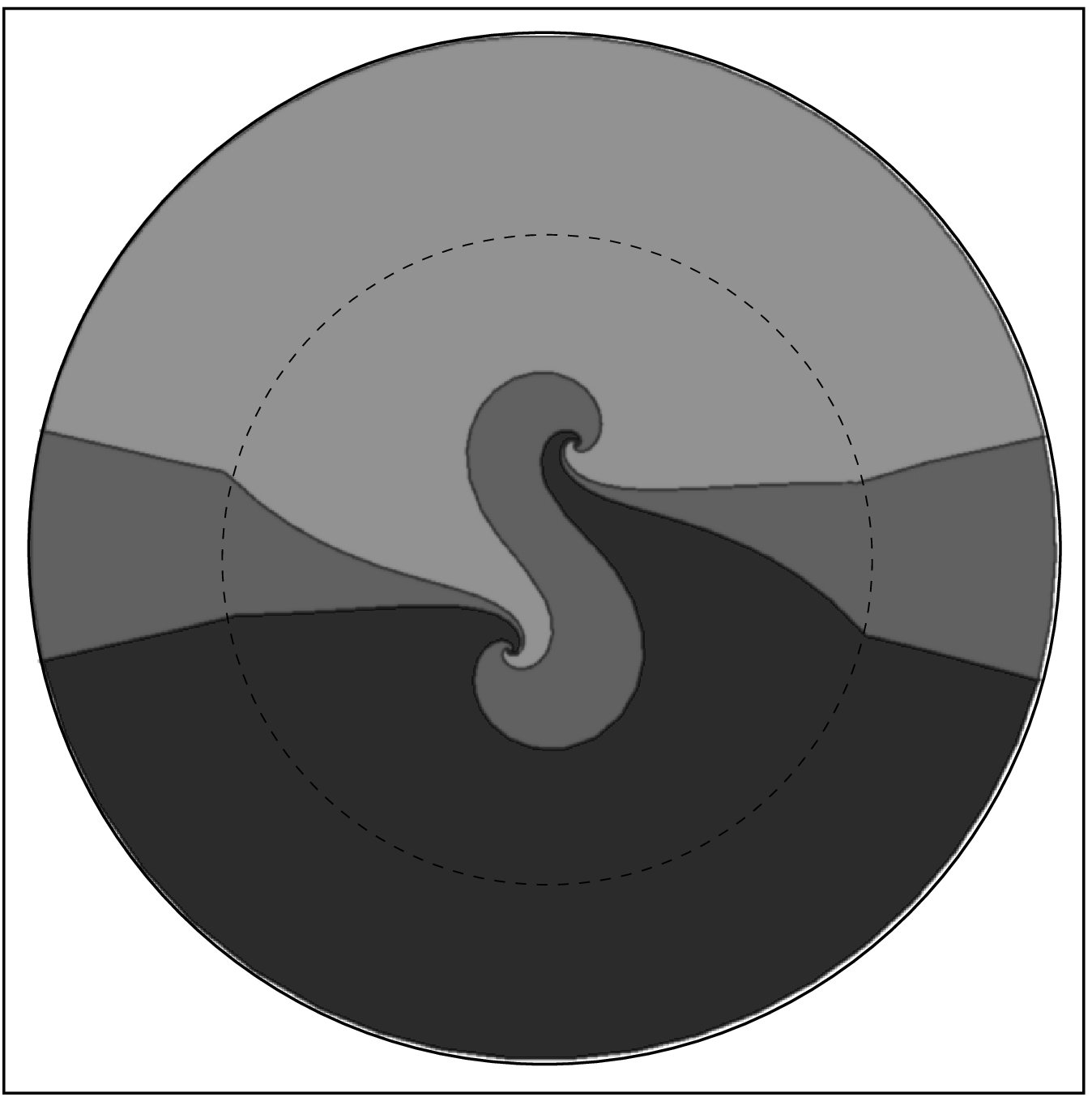}}\\
 \hfill{}\subfigure{\includegraphics[width=6cm]{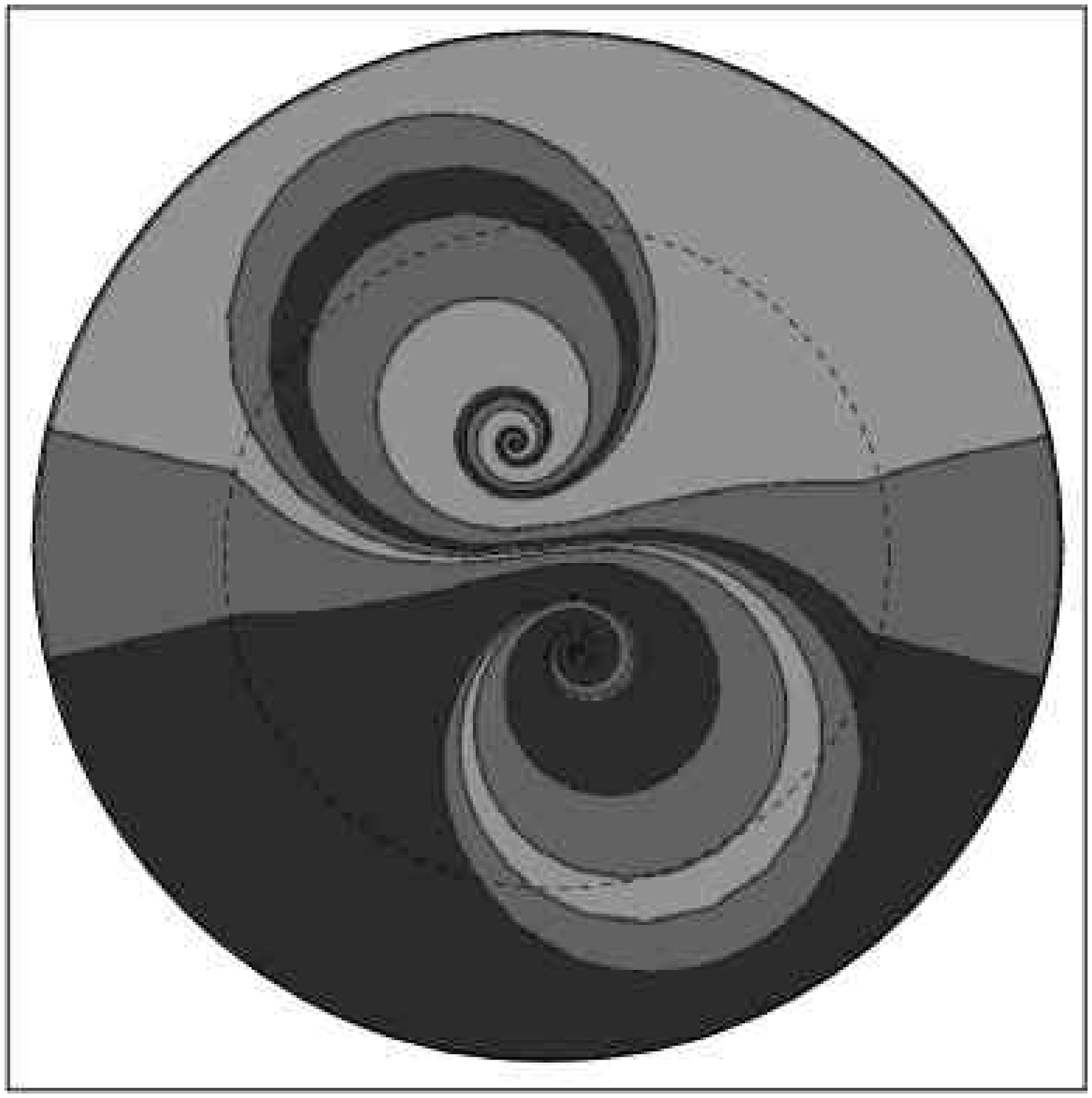}}\hfill{}

\caption{The canonical squid sectors when $k=1$ for different values of $\varepsilon$. }

\label{fig_squid} 
\end{figure}

\begin{defn}
\label{def:squid-sectors} We consider a good sector $W\subset\Sigma_{0}$
on which we fix a good angle $\theta$, and $\rho>0$ so that the
conclusions of Lemmas~\ref{lem:mino_deriv} and \ref{P_3} hold. 
\begin{enumerate}
\item We first build the squid sectors in $z$-coordinate around the ball
$B_{0}$ of center $0$. The other sectors are deduced by translations
and changes of sheet. These sectors are somewhat wider than the $\partial V_{j,\eps}^{\pm}$
of Figure~\ref{fig_fleur}, so as to give an open covering of $r\mathbb{D}\setminus\{x_{0},\dots,x_{k}\}$.
For a given $\varepsilon\in W$ we define\begin{eqnarray*}
t_{\varepsilon} & := & \kappa\left|\left|\varepsilon\right|\right|^{-1}\end{eqnarray*}
 where $\kappa>\frac{1}{k}r^{-k}$ is sufficiently small so that $t_{\varepsilon}\leq\frac{1}{3M(\gamma_{j})}$
and $z^{-1}\left(t_{\varepsilon}\ww{D}\backslash B_{0}\right)$ does
not meet the disc $\sqrt{k}\left|\left|\varepsilon\right|\right|\ww{D}$
containing the roots of $P_{\varepsilon}$ (the number $M\left(\gamma_{j}\right)$
is defined in (\ref{eq_mu_j})). This choice of $\kappa$ can be made
independently on $\varepsilon$ according to the estimate (\ref{eq:estim_bandwidth})
of Lemma~\ref{P_3}. The following construction corresponds to Figure~\ref{fig:z_squid}:

\begin{enumerate}
\item inside $t_{\varepsilon}\ww{D}\backslash B_{0}$ the domain is a horizontal
strip; the distance $w$ between the horizontal boundary of the strip
and the parallel diameter of $B_{0}$ is fixed once and for all satisfying
$0<w<\frac{2}{3k}r^{-k}$. We call it the \textbf{width} of the squid
sectors. 
\item outside the disk $t_{\varepsilon}\ww{D}$ the domain is a slanted
strip comprised between straight lines making an angle $\theta$ with
the horizontal. The distance between the outermost lines and the center
of $B_{0}$ is $\frac{M(\gamma_{j})}{2}+w$. 
\end{enumerate}
The domain in $z$-space is taken so that no two points project on
the same $x$-point, \emph{i.e.} differ by a period of $P_{\eps}(x)$.

\item The projection in $x$-space of such a domain is called a \textbf{squid
sector} and denoted $\sect{j,\varepsilon}{\pm}$. See Figure~\ref{fig_squid}. 
\item For $\varepsilon=0$ we do the same construction with $t_{0}:=+\infty$
and $\theta:=0$, which corresponds to the half-plane with holes $\left\{ Im\left(z\right)<w\right\} \backslash B_{0}$.
See Figure~\ref{figsect_0} and Figure~\ref{figsect_eps}(d). 
\end{enumerate}
\end{defn}

\begin{figure}
\includegraphics[width=80mm]{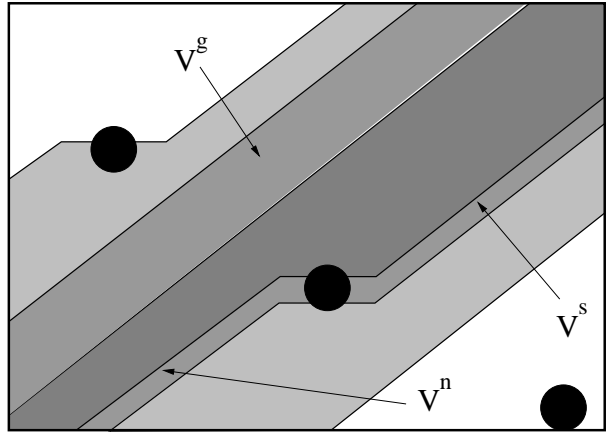}

\caption{\label{fig:intersection} In the case $k=1$ the intersections of
the two sectors is formed of $\sect{\varepsilon}{n}$, $V_{\varepsilon}^{s}$
and $\sect{\varepsilon}{g}$ . To visualize $\sect{\varepsilon}{g}$,
we need to take a translate of one of the sectors by a period.} 
\end{figure}

\begin{figure}
\hfill{}\subfigure[ Here $p_0$ and $p_1$ are of saddle type while $p_2$ is of node type.]{\includegraphics[width=7cm]{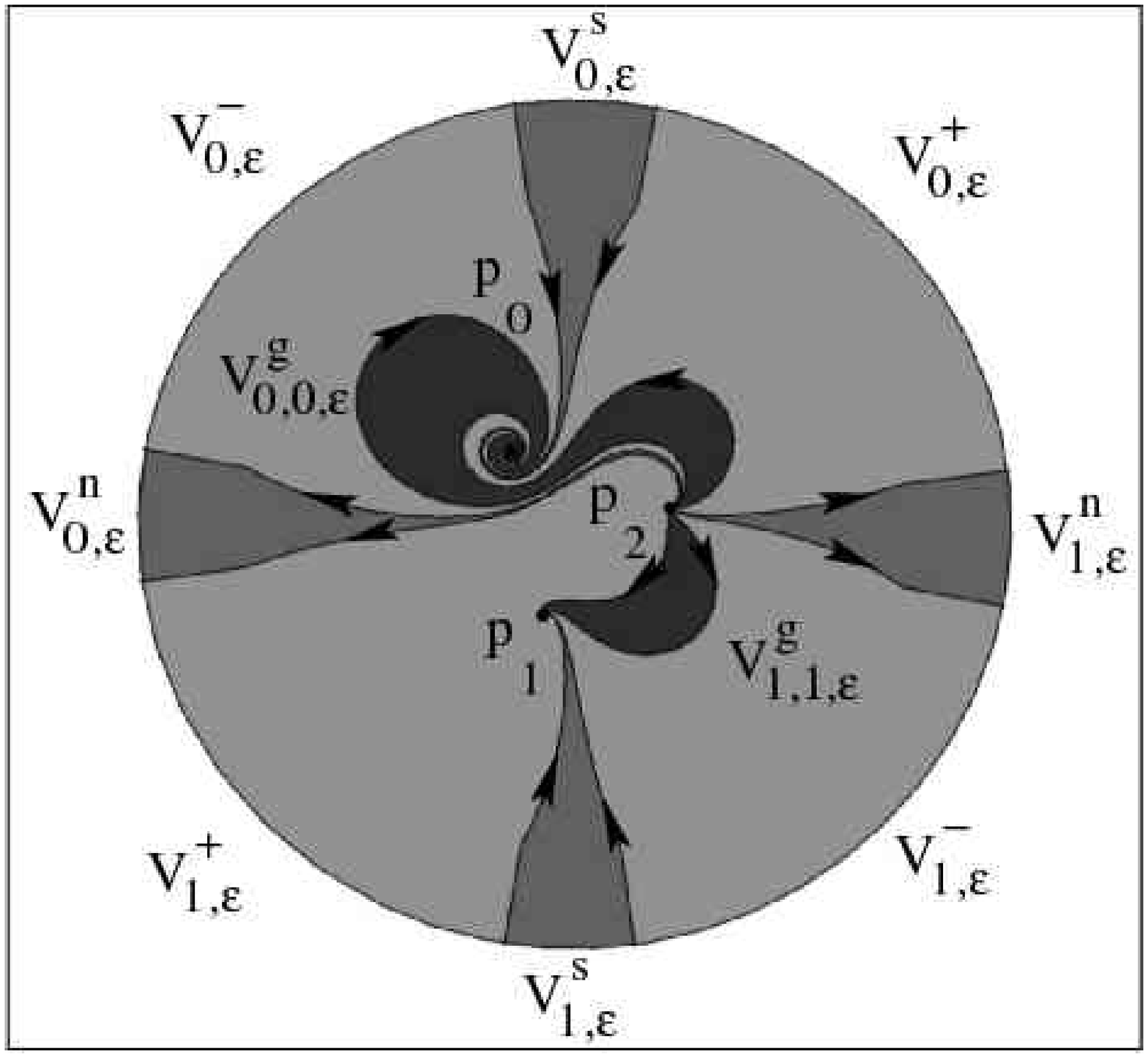}}\hfill{}\subfigure[ Here $p_1$ is of saddle type whereas $p_0$ and $p_2$ are of node type.]{\includegraphics[width=7cm]{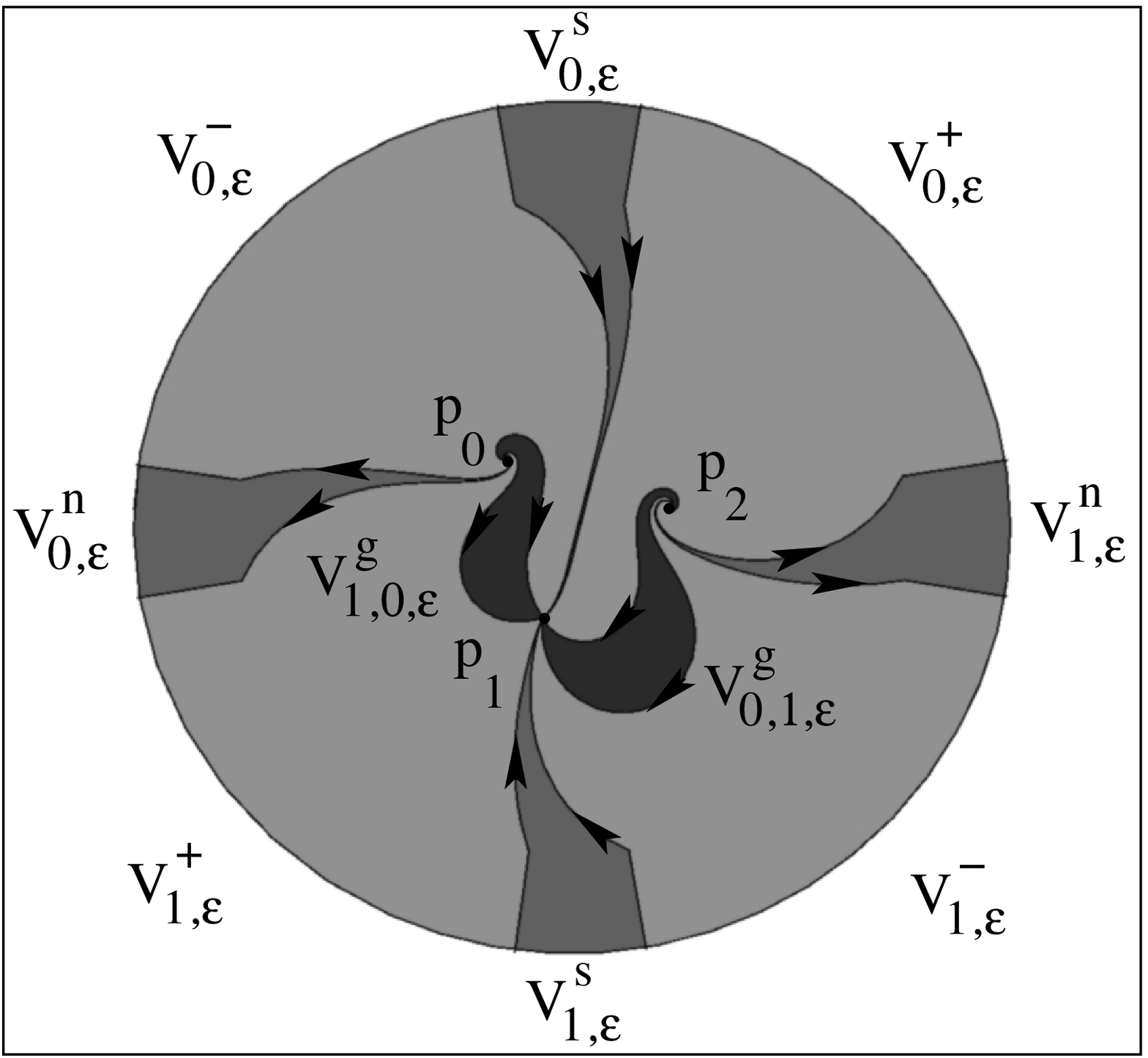}}\hfill{}

\caption{Examples of (non-equivalent) squid sectors in the case $k=2$ for
the same value of $\eps$ and different choices of $\theta$.}

\label{fig_sectors} 
\end{figure}

\begin{lem}
Any compact subset of $\sect{j,0}{\pm}$ is contained in a compact
set of $\sect{j,\varepsilon}{\pm}$ for $\eps$ sufficiently small. 
\end{lem}
\begin{proof}
This is obvious from the Figures \ref{figsect_0} and \ref{fig:z_squid}. 
\end{proof}
\begin{lem}
\label{lem:intersect}\emph{\cite{DS}, \cite{O}.} In the neighborhood
of a generic $\eps\in\Sigma_{0}$, any squid sector is adherent to
two singular points $x_{s}$ and $x_{n}$, one being an attractor
and the other being a repeller for $\Xi_{\eps}(\theta)$ given in
(\ref{eq:theta_theta}). 
\begin{enumerate}
\item In the case $k=1$ the intersection of the two squid sectors is formed
by three sectors $\sect{\varepsilon}{s}$, $\sect{\varepsilon}{n}$
and $\sect{\varepsilon}{g}$, see Figure~\ref{fig_squid} and Figure~\ref{fig:intersection}.
The upper-indices $s$ (\emph{resp.} $n$, $g$) refer to {}``saddle-like''
(\emph{resp.} {}``node-like'' and {}``gate''). The gate structure
was introduced by Oudkerk \cite{O}. $\sect{\varepsilon}{s}$ is adherent
to an attracting point $x_{s}$ for $\Xi_{\eps}(\theta)$, $\sect{\varepsilon}{n}$
is adherent to a repelling point $x_{n}$ for $\Xi_{\eps}(\theta)$
and $V_{\eps}^{g}$ is adherent to both. 
\item In the case $k>1$ the intersection of two consecutive squid sectors
is given by one or two sectors, namely

\begin{itemize}
\item in the case of $\sect{j,\varepsilon}{+}\cap\sect{j,\varepsilon}{-}$
a sector $\sect{j,\varepsilon}{s}$, and an additional sector $\sect{j,\sigma(j),\varepsilon}{g}$
if and only if $\sigma(j)=j$. 
\item in the case of $\sect{j+1,\varepsilon}{+}\cap\sect{j,\varepsilon}{-}$
a sector $\sect{j,\varepsilon}{n}$, and an additional sector $\sect{j,\sigma(j),\varepsilon}{g}$
if and only if $\sigma(j+1)=j$. 
\end{itemize}
$\sect{j,\varepsilon}{s}$ is adherent to an attracting point for
$\Xi_{\eps}(\theta)$, $\sect{j,\varepsilon}{n}$ is adherent to a
repelling point for $\Xi_{\eps}(\theta)$ and $\sect{j,\sigma(j),\varepsilon}{g}$
exists if and only if the two sectors share the same singular points,
in which case it is adherent to both.

In order to be able to give definitions valid for all sectors we will
often use the notation \begin{equation}
\sect{j,\varepsilon}{g}:=\sect{j,\sigma(j),\varepsilon}{g}.\label{sector_gate}\end{equation}

\item Two non consecutive squid sectors $V_{j,\eps}^{+}$ and $V_{\ell,\eps}^{-}$
intersect along a gate sector $V_{j,\ell,\eps}^{g}$ if and only if
$\ell=\sigma(j)$ (see for instance Figure~\ref{fig_sectors}), \emph{i.e.}
in the case where they are adherent to the same singular points. 
\end{enumerate}
\end{lem}

\begin{rem}
The boundary of a squid sector $\sect{j,\varepsilon}{\pm}$ can be
split into four parts for generic $\varepsilon$, which we explicit
in $x$-coordinate since its structure in $z$-coordinate is simpler
to recover :
\end{rem}
\begin{itemize}
\item starting from $x_{n}$ we follow a \lq\lq S''-shaped spiraling
curve to reach $x_{s}$
\item then we reach the circle $r\ww{S}^{1}$ by following a spiral
\item we follow an arc of circle, corresponding to $\partial\sect{j,\varepsilon}{\pm}$
\item we close the loop by reaching $x_{n}$ along another spiral.
\end{itemize}
\begin{lem}
Let $W\subset\Sigma_{0}$ be a good sector. The squid sectors can
be taken depending analytically on $\eps\in W$ and continuously on
$\varepsilon\in W\cup\left\{ 0\right\} $.
\end{lem}

\begin{rem}
Simply taking a horizontal strip in $z$-coordinate would not work
for our construction. Indeed we need to discriminate between singular
points of saddle-type $x_{s}$ and node-type $x_{n}$ for the foliation.
This goal is achieved by taking slanted strips near infinity (corresponding
to spirals in $x$-coordinates) which forces the leaves respectively
either not to or to accumulate on the point when following the sector.
We discuss this property in the upcoming Section~\ref{sub:Study-model}.

Moreover simply taking a slanted strip in $z$-coordinate would not
ensure the convergence $\sect{j,\varepsilon}{\pm}\to\sect{j,0}{\pm}$,
as $\theta$ is fixed when $\varepsilon$ belongs to a good sector
$W_{i}$ and $\theta$ is taken to be $0$ when $\varepsilon=0$.
This is the reason why we need a mixed construction for the squid
sectors.
\end{rem}

\subsection{\label{sub:Study-model}Study of the foliations of the model family}

We consider the first integral of the model \eqref{int_prem} over
fibered squid sectors \begin{equation}
\ssect{j,\varepsilon}{\pm}:=\sect{j,\varepsilon}{\pm}\times\ww{C}\end{equation}
 constructed with an adapted set of squid sectors $\sect{j,\varepsilon}{\pm}$.

\begin{figure}
\hfill{}\subfigure{\includegraphics[width=7cm]{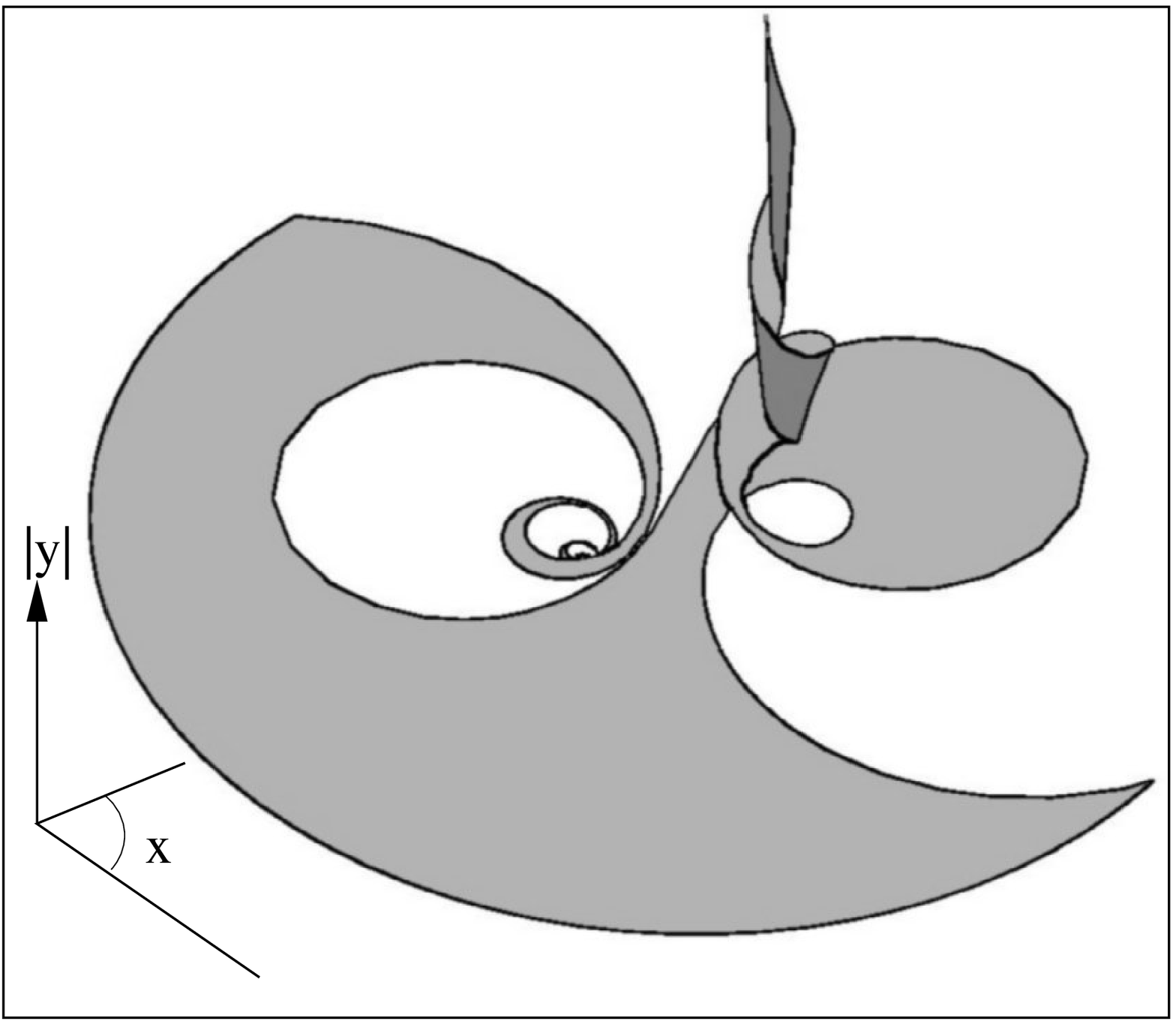}}\hfill{}\subfigure{\includegraphics[width=7cm]{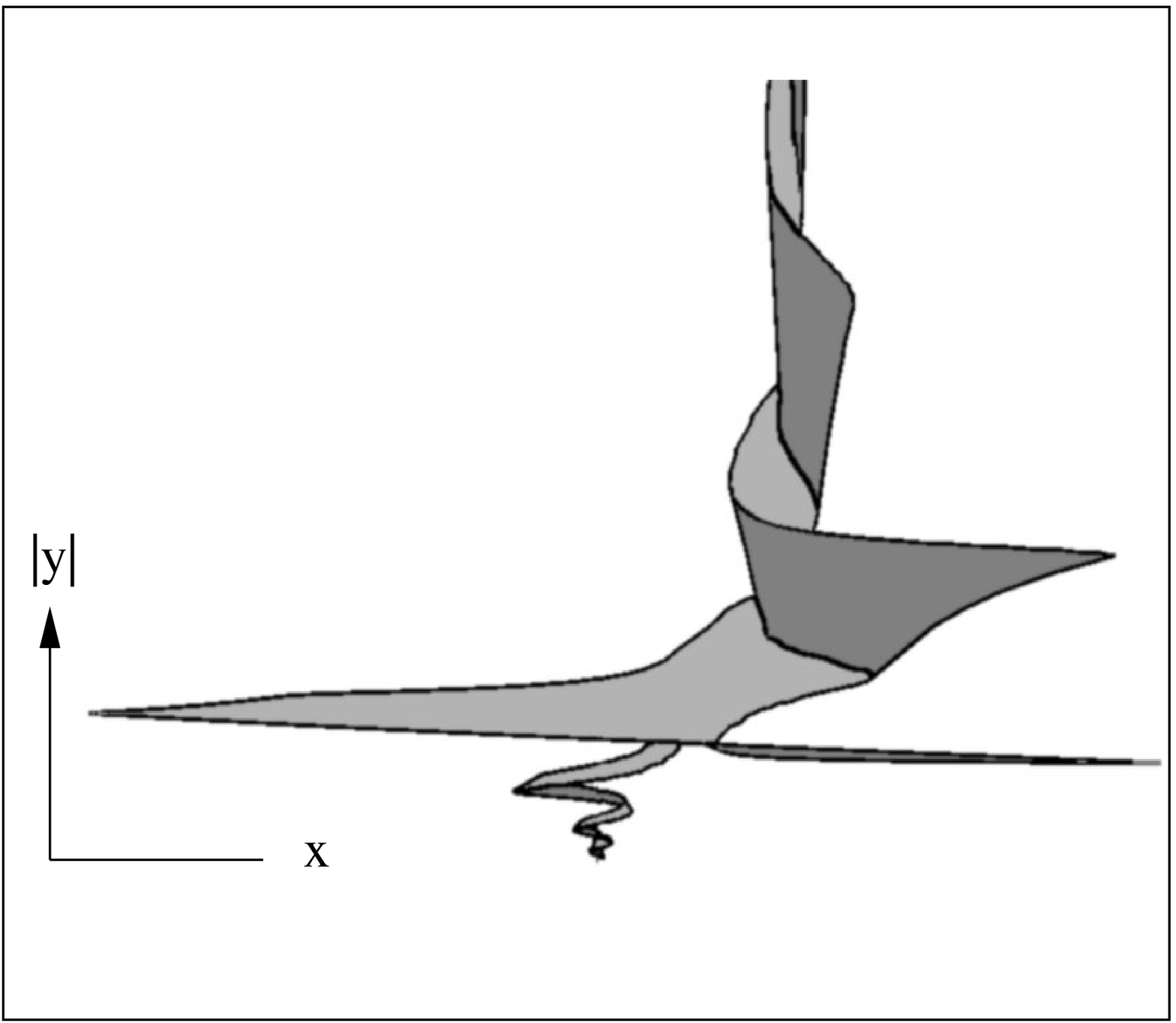}}\hfill{}

\caption{\label{fig:foliation}Modulus of a leaf of the model foliation over
a squid sector ($k=1$). These drawings justify the terms {}``node
type'' (on the left of each figure) and {}``saddle type'' (on the
right) qualifying the singular points.} 
\end{figure}

\begin{defn}
\label{def:pt_sing}~ 
\begin{enumerate}
\item Let $\varepsilon\in\Sigma_{0}$. Each fibered squid sector is adherent
to two distinct singular points $(x_{n},0)$ and $(x_{s},0)$ of $X_{\eps}^{M}$.
The point $(x_{n},0)$ (\emph{resp.} $(x_{s},0)$) is said to be of
\textbf{node type} (\emph{resp.} \textbf{saddle type}) if $Re\left(\exp(i\theta)P_{\eps}'(x_{n})\right)>0$
(\emph{resp.} $Re\left(\exp(i\theta)P_{\eps}'(x_{s})\right)<0$).
We note $p_{n}:=(x_{n},0)$ and $p_{s}:=(x_{s},0)$. Depending on
the context we may also use the notation $p_{j,n}=(x_{j,n},0)$ or
$p_{j,n}^{\pm}=(x_{j,n}^{\pm},0)$ (similarly $p_{j,s}$ or $p_{j,s}^{\pm}$)
to emphasize that we consider the singular points of node and saddle
type associated to the fibered squid sector $\ssect{j,\varepsilon}{\pm}$. 
\item If $\varepsilon=0$ we set $p_{j,n}:=p_{j,s}:=\left(0,0\right)$. 
\end{enumerate}
\end{defn}
For each fibered squid sector we fix the principal holomorphic determination
$H_{j,\varepsilon}^{M,\pm}$ of the first integral $H_{\varepsilon}^{M}$
(defined in \eqref{int_prem}) on $\ssect{j,\varepsilon}{\pm}$. This
is defined starting on the boundary in sector $\partial V_{0,\eps}^{+}$
of Figure~\ref{fig_fleur}, turning on $|x|=r$ in the positive direction
and then extending to the interior of the sectors.

\begin{figure}
\hfill{}\subfigure[in $x$-coordinate]{\includegraphics[width=7cm]{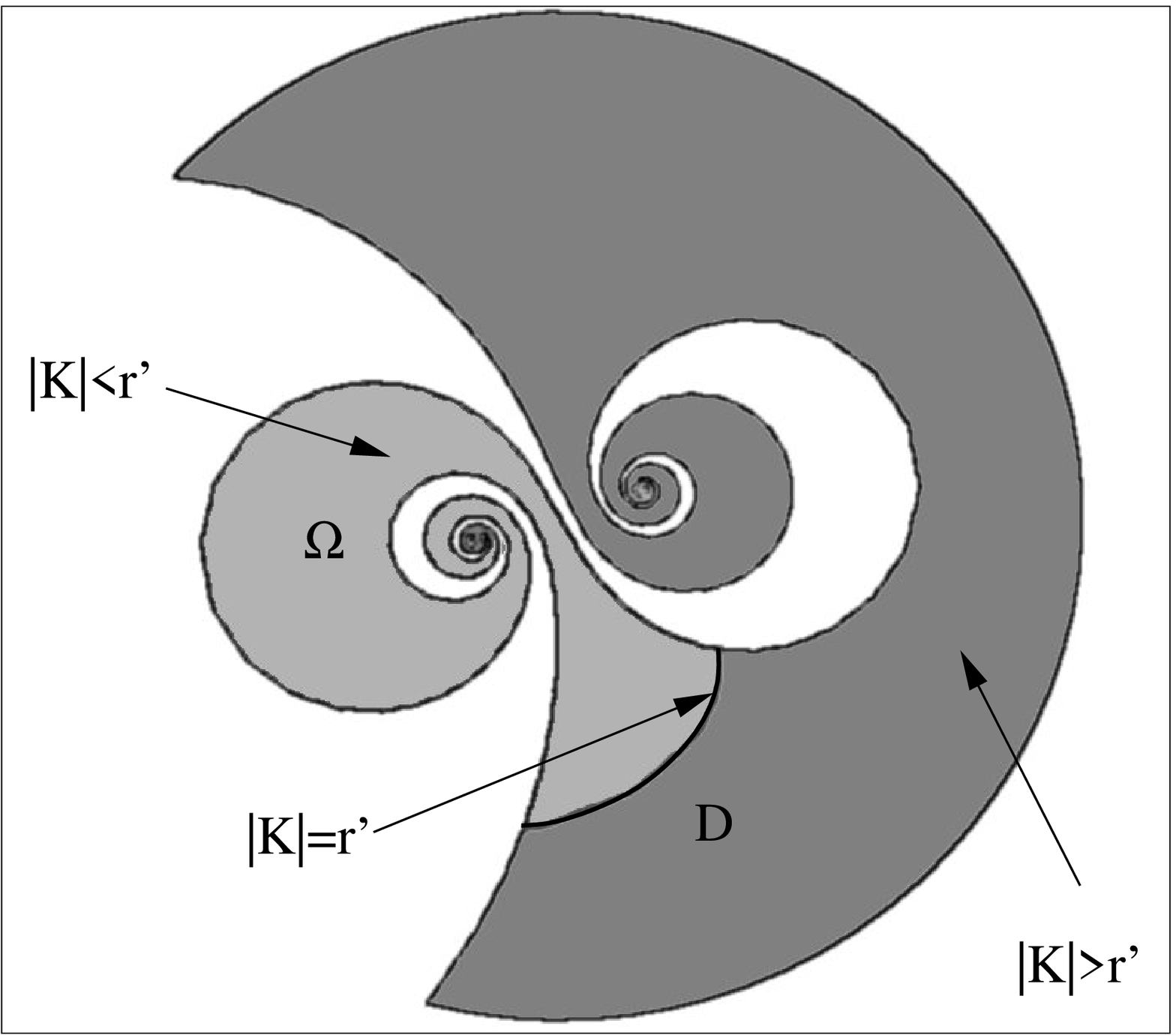}}\hfill{}\subfigure[in $z$-coordinate]{\includegraphics[width=7cm]{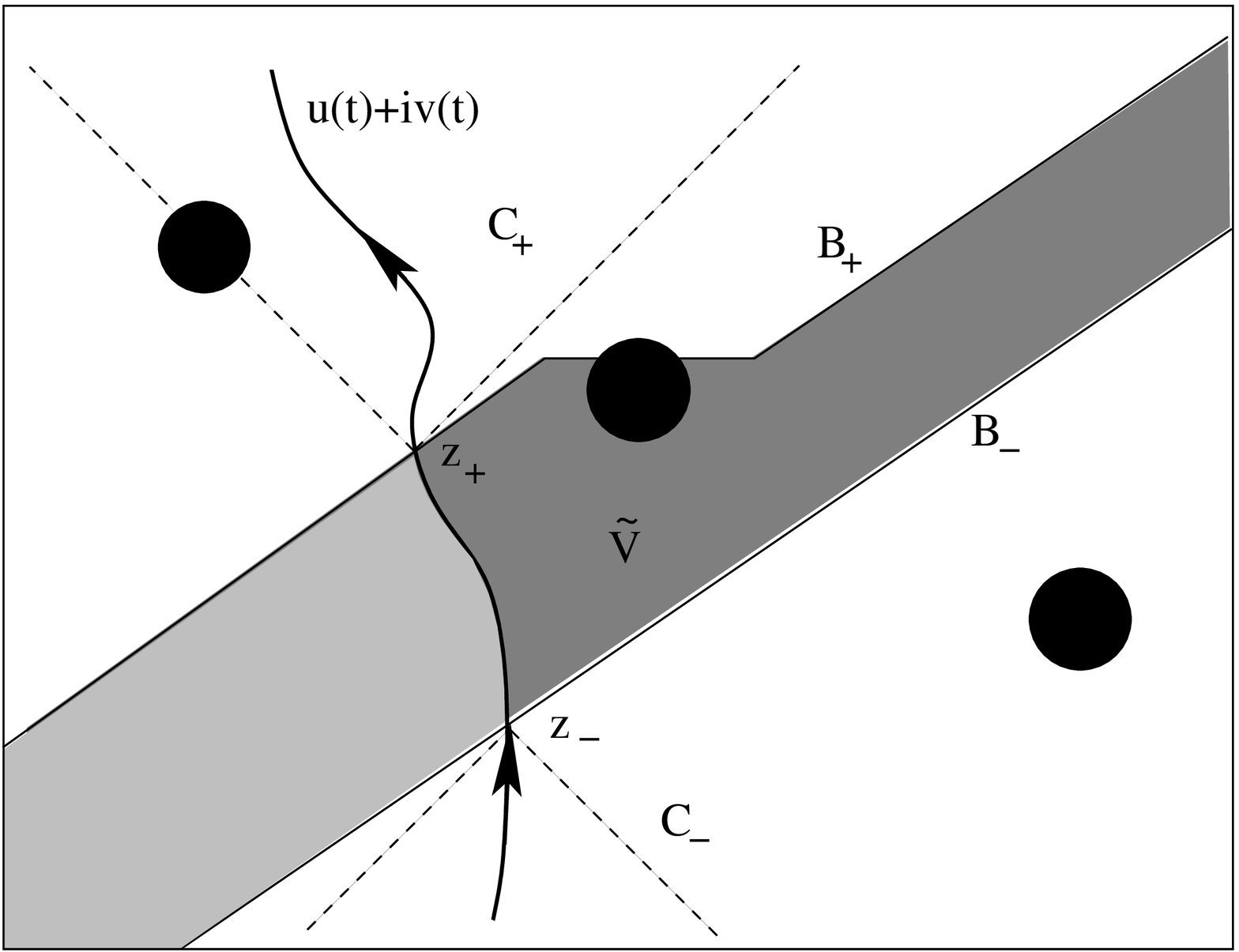}}\hfill{}

\caption{\label{fig:squid_modulus}The trace $D:=\left\{ x\,:\,\left|K_{j,\varepsilon}^{\pm}\left(x,h\right)\right|=r'\right\} $
of the leaf on $\sect{j,\varepsilon}{\pm}\times r'\ww{S}^{1}$.} 
\end{figure}

\begin{lem}
The convergence $H_{j,\varepsilon}^{\pm,M}\to H_{j,0}^{\pm,M}$ is
uniform on any compact set of $V_{j,0}^{\pm}$. 
\end{lem}
\begin{prop}
\label{pro:first-integ-model} Let $r>0$ and $\varepsilon\in\Sigma_{0}\cup\left\{ 0\right\} $
be given. The foliation $\fol{j,\varepsilon}^{\pm}$ induced by $X_{\varepsilon}^{M}$
on $\ssect{j,\varepsilon}{\pm}$ satisfies the following properties
: 
\begin{enumerate}
\item For each leaf $\mathcal{L}$ of $\fol{j,\varepsilon}^{\pm}$ there
exists $h\in\ww{C}$ such that $\mathcal{L}=\left(H_{j,\varepsilon}^{M,\pm}\right)^{-1}\left(h\right)$.
On the other hand, for any $h\in\ww{C}$ the set $\left(H_{j,\varepsilon}^{M,\pm}\right)^{-1}\left(h\right)$
is a leaf of $\fol{j,\varepsilon}^{\pm}$. 
\item There exists a holomorphic function $K\,:\,\sect{j,\varepsilon}{\pm}\times\ww{C}\to\ww{C}$
such that the leaf of $\fol{j,\varepsilon}^{\pm}$ corresponding to
$h\in\ww{C}$ coincides with the graph of $x\mapsto K\left(x,h\right)$. 
\item There exists $r,\rho>0$ such that for any $r'>0$, any $\varepsilon\in\Sigma_{0}\cup\left\{ 0\right\} $
with $\left|\left|\varepsilon\right|\right|\leq\rho$ and any $\left(\overline{x},\overline{y}\right)\in\sect{j,\varepsilon}{\pm}\times r'\ww{D}\backslash\left\{ 0\right\} $
the closure of $\left[K\left(\cdot,H_{j,\varepsilon}^{\pm,M}\left(\overline{x},\overline{y}\right)\right)\right]^{-1}\left(r'\ww{S}^{1}\right)$
is a (connected) real analytic curve $D$ which separates $\sect{j,\varepsilon}{\pm}$
into two connected components, and crosses transversely the boundary
$\partial\sect{j,\varepsilon}{\pm}$ in exactly two points (see Figure~\ref{fig:squid_modulus}(a)).
One connected component of $\sect{j,\varepsilon}{\pm}\backslash D$
accumulates on $x_{j,n}$ while the other accumulates on $x_{j,s}$.
Moreover $D\to\left\{ x_{j,n}\right\} $ as $r'\to0$. 
\item Let $r',r>0$ be given. Then :

\begin{enumerate}
\item $H_{j,\varepsilon}^{M,\pm}\left(\sect{j,\varepsilon}{\#}\times r'\ww{D}\right)=\ww{C}$
for $\#\in\left\{ \pm,n,g\right\} $ 
\item $H_{j,\varepsilon}^{M,\pm}\left(\sect{j,\varepsilon}{s}\times r'\ww{D}\right)=\eta(r')\ww{D}$
with $\eta\left(r'\right)=r'O\left(1\right)$ uniformly in $\varepsilon$
belonging to a good sector. 
\end{enumerate}
\item There exists a unique distinguished leaf, the zero level curve of
$H_{j,\varepsilon}^{M,\pm}$, which is adherent to both singular points
$p_{j,s}$ and $p_{j,n}$. The distinguished leaves over the different
sectors glue in a global leaf in $r\mathbb{D}\times\mathbb{C}$, which
actually is $\left(r\ww{D}\backslash\left\{ x_{0},\dots,x_{k}\right\} \right)\times\left\{ 0\right\} $. 
\item Assume $\varepsilon\neq0$. In each sector $\ssect{j,\varepsilon}{\pm}$
all leaves, except the distinguished leaf, are adherent to exactly
one of the singular points, namely $p_{j,n}$. 
\end{enumerate}
\end{prop}
\begin{proof}
We drop all indices so as to enlighten the ideas. Because the squid
sectors are simply connected the first integral $H$ is univalued.
We note $\mathcal{L}_{h}:=H^{-1}\left(h\right)$ the curve of level
$h$ of $H$.

\textbf{(1)} Firstly the relation $X\cdot H=0$ yields that the function
$H$ is constant on each leaf of $\fol{}$. Thus each $\mathcal{L}_{h}$
is a union of leaves. For a fixed $\overline{x}\in V$ the map $H_{\overline{x}}\,:\, y\mapsto H\left(\overline{x},y\right)$
is linear and invertible. Hence, given $\left(\overline{x},\overline{y}\right)\in\ssect{}{}$,
there exists $h:=H\left(\overline{x},\overline{y}\right)$ such that
$\left(\overline{x},\overline{y}\right)\in\mathcal{L}_{h}$; in particular
every leaf of $\fol{}$ is contained in some $\mathcal{L}_{h}$. On
the other hand the injectivity of $H_{\overline{x}}$ implies that
if $\left(\overline{x},\overline{y}\right)$ and $\left(\overline{x},\tilde{y}\right)$
lie in distinct leaves then $H\left(\overline{x},\overline{y}\right)\neq H\left(\overline{x},\tilde{y}\right)$.
The conclusion follows since any $\left(\tilde{x},\tilde{y}\right)$
may be connected within a leaf to some $\left(\overline{x},y\right)$
using the fact that $\fol{}$ is transverse to the lines $\left\{ x=\mbox{cst}\right\} $
on $\ssect{}{}$.

\textbf{(2)} is a direct consequence of (1). Each leaf $\mathcal{L}_{h}$
coincides with the graph of the holomorphic function \begin{eqnarray}
K\left(\cdot,h\right)\,:\, x & \mapsto & h\frac{y}{H\left(x,y\right)}\,.\end{eqnarray}
 (Note that $\frac{y}{H\left(x,y\right)}$ is a function of $x$ alone.)

\textbf{(3)} We work in $z$-coordinate as in Lemma~\ref{P_2} and
assume that $h\neq0$. The map $\tilde{K}\,:\, z\mapsto K\left(x\left(z\right),h\right)$
satisfies the differential equation\begin{eqnarray}
\frac{d\tilde{K}}{dz}\left(z\right) & = & \tilde{K}\left(z\right)\left(1+ax\left(z\right)^{k}\right)\,.\label{eq:z-coord_A}\end{eqnarray}
 Hence \begin{eqnarray}
\tilde{K}\left(z\right) & = & \overline{y}\exp A\left(z\right)\neq0\\
A\left(z\right) & = & \int_{\overline{z}}^{z}\left(1+ax\left(s\right)^{k}\right)ds\nonumber \end{eqnarray}
 where $A$ is holomorphic on a neighborhood of the strip $\tilde{V}$
corresponding in $z$-coordinate to the closure of the squid sector
$V$ (see Lemmas~\ref{P_3} and \ref{P_4}). We let $z=u+iv$. The
level sets $\left\{ \left|\tilde{K}\right|=r'\right\} =\left\{ Re\left(A\right)=\ln\frac{r'}{|\ov{y}|}\right\} $
for fixed $\ov{y}$ and different $r'>0$ define a regular real analytic
foliation of $\tilde{V}$ through the differential system \begin{eqnarray}
\dot{u} & = & -\pp{v}Re\left(A\left(u+iv\right)\right)=Im\left(a(x\left(u+iv\right))^{k}\right)\\
\dot{v} & = & \pp{u}Re\left(A\left(u+iv\right)\right)=1+Re\left(a(x\left(u+iv\right))^{k}\right)\,\,\,\,.\nonumber \end{eqnarray}
 A first observation is that $t\mapsto v\left(t\right)$ is strictly
monotonous provided that $r^{k}\left|a\right|<1$. Indeed we have
$\left|\dot{u}\right|<\left|a\right|r^{k}$ and $\left|\dot{v}-1\right|<\left|a\right|r^{k}$.
We will assume now that $r^{k}\left|a\right|<1$, which can be achieved
for $r$ sufficiently small independently on $\varepsilon$. By integrating
the previous inequalities between $0$ and $t$ we obtain : \begin{eqnarray}
\left|v\left(t\right)-v\left(0\right)-t\right| & \leq & \left|t\right|\left|a\right|r^{k}\label{eq:z-coord_v-estim}\\
\left|\frac{u\left(t\right)-u\left(0\right)}{v\left(t\right)-v\left(0\right)}\right| & \leq & \frac{r^{k}\left|a\right|}{1-r^{k}\left|a\right|}\end{eqnarray}
 We further require that $\eta:=\frac{r^{k}\left|a\right|}{1-r^{k}\left|a\right|}<1$
by potentially diminishing $r$ if necessary. The curve $\tilde{D}\,:\, t\mapsto\left(u\left(t\right),v\left(t\right)\right)$
lies thus in the union of the conic regions $\mathcal{C}_{+}:=\left\{ \left|u-u\left(0\right)\right|\leq\eta\left(v-v\left(0\right)\right)\right\} $
(for $t\geq0$) and $\mathcal{C}_{-}:=\left\{ v-v\left(0\right)\leq-\frac{1}{\eta}\left|u-u\left(0\right)\right|\right\} $
(for $t\leq0$).

Let us write $\partial\tilde{V}=B_{-}\cup B_{+}$ where $B_{+}$ (\emph{resp.}
$B_{-}$) comes from the upper (\emph{resp.} lower boundary) of $\tilde{V}$
on $z$-coordinate (see Figure~\ref{fig:squid_modulus}(b)). Because
of (\ref{eq:z-coord_v-estim}) we derive that the integral curve obtained
for $\overline{z}=u\left(0\right)+iv\left(0\right)$ cuts $\partial\tilde{V}$
in at least two points $z_{+}\in B_{+}$ and $z_{-}\in B_{-}$. Indeed
$\left|\theta\left(\varepsilon\right)\right|<\frac{\pi}{4}$. Hence
if we take any starting point $z=u\left(0\right)+iv\left(0\right)$
in $B_{+}$ (\emph{resp.} $B_{-}$) the set $\mathcal{C}_{+}$ (\emph{resp.}
$\mathcal{C}_{-}$) intersects $\tilde{V}$ only at $z$. This yields
the uniqueness of $z_{\pm}$. See Figure~\ref{fig:squid_modulus}(b).

The fact that one component accumulates of $x_{n}$ and the other
on $x_{s}$ comes clearly from Lemma~\ref{P_4}. If $\varepsilon\in\Sigma_{0}$
the fact that $D$ cannot converge to $\left\{ x_{s}\right\} $ follows
from the construction of $\theta$ since $\cos\theta>0$ and \begin{eqnarray}
\left|K\left(x\left(te^{i\theta}+\overline{z}\right),h\right)\right| & \sim_{t\to\pm\infty} & A\exp\left(t\cos\theta\right)\\
\lim_{t\to\pm\infty}x\left(te^{i\theta}+\overline{z}\right) & = & x_{\#}\end{eqnarray}
 where $\#=n$ (\emph{resp.} $\#=s$) if $Re\left(e^{i\theta}P_{\varepsilon}'\left(x_{n}\right)\right)>0$
(resp. $Re\left(e^{i\theta}P_{\varepsilon}'\left(x_{s}\right)\right)<0$)
and $t\to-\infty$ (\emph{resp.} $t\to+\infty$). More details can
be found in Lemma~\ref{lem:modulus_estimate}.

\textbf{(4)} According to the discussion made just above we have,
for fixed $\overline{y}$,\begin{eqnarray}
\lim_{x\to x_{n}\,,\, x\in V}\left|H\left(x,\overline{y}\right)\right| & = & \infty\\
\lim_{x\to x_{s}\,,\, x\in V}\left|H\left(x,\overline{y}\right)\right| & = & 0\end{eqnarray}
 Because $H$ is linear in $\overline{y}$ and the argument of $\overline{y}$
takes all values, every $h\in\ww{C}$ is reached on any sector accumulating
on $x_{n}$. The same argument shows that $H\left(V^{s}\right)$ is
a disk of radius \begin{eqnarray}
\eta\left(r'\right) & = & \sup_{V^{s}\times r'\ww{D}}\left|H\left(x,y\right)\right|\\
 & = & r'\max_{\partial V^{s}\backslash\left\{ x_{s}\right\} }\left|H\left(x,1\right)\right|\,.\nonumber \end{eqnarray}
 The fact that \begin{eqnarray}
H_{0}^{M}\left(x,y\right) & = & yx^{-a\left(0\right)}\exp\frac{1}{kx^{k}}\end{eqnarray}
 yields $\eta_{0}\left(r'\right)=r'\exp\left(\frac{c}{kr^{k}}\right)$
for some constant $c>1$ depending on the width of $V_{0}^{s}$. This
proves the claim as $\partial\sect{\varepsilon}{s}\backslash\left\{ x_{s}\right\} \to\partial\sect{0}{s}\backslash\left\{ 0\right\} $
and $H_{\varepsilon}\left(\cdot,1\right)\to H_{0}\left(\cdot,1\right)$
uniformly on $\overline{\sect{0}{s}}\backslash\left\{ 0\right\} $.

\textbf{(5)} The line $\left\{ y=0\right\} \backslash\left\{ p_{0},\dots,p_{k}\right\} $
clearly is the curve of level $0$ of $H$, so is a leaf. All the
principal determinations of $H_{\varepsilon}^{M}$ agrees on $\left\{ y=0\right\} $
so that these distinguished leaves glue in a global leaf.

\textbf{(6)} It follows from (3). Indeed, let $\mathcal{L}_{h}$ be
a leaf of $\fol{}$ with $h\neq0$. On the one hand $\mathcal{L}_{h}$
cannot accumulate on $p_{s}$ because $x_{s}$ belongs to the closure
of $\left\{ \left|K\left(\cdot,h\right)\right|>r'\right\} $ for all
$r'>0$. On the other hand $D\to\left\{ x_{n}\right\} $ as $r'\to0$
so that $x_{n}$ lies in $\overline{\mathcal{L}_{h}}$. 
\end{proof}

\section{\label{sec:center}The center manifolds}

This section is purely orbital so we work with a prepared family $X_{\eps}$
of vector fields of the form \eqref{eq_X}. Let us define for $k>1$
the sectors \begin{equation}
\sect{j,\varepsilon}{}:=\sect{j,\varepsilon}{+}\cup\sect{j,\varepsilon}{-}\end{equation}
 and \begin{equation}
\ssect{j,\varepsilon}{}:=\sect{j,\varepsilon}{}\times r'\mathbb{D},\end{equation}
 built from the squid sectors obtained in Definition~\ref{def:squid-sectors}
(see Figure~\ref{fig:center}). If $k=1$ we merge $\sect{0,\varepsilon}{+}$
and $\sect{0,\varepsilon}{-}$ only on the saddle and gate sides.
This yields a sector of opening greater than $2\pi$, which must be
considered in the universal covering of $x$-space punctured at $x_{n}$.

\begin{lem}
\label{attractor} Each sector $\ssect{j,\varepsilon}{}$ contains
a singular point $p_{j,s}=(x_{j,s},0)$ such that the $x$-eigenvalue
(\emph{resp.} $y$-eigenvalue) of the linearized vector field of $\exp(i\theta\left(\varepsilon\right))X_{\eps}$
at $p_{j,s}$ has a negative (\emph{resp.} positive) real part. 
\end{lem}
\bigskip{}
 We choose to study the family \eqref{eq_X} over a fixed polydisk
in $(x,y)$-space, taken as $r\mathbb{D}\times r'\mathbb{D}$. For
$\eps=0$ the vector field has a formal center manifold given by a
(generically divergent) power series $y=\hat{S}(x)=\sum_{n\geq2}a_{n}x^{n}$.
The sum of this series gives $k$ center manifolds as graphs of functions
$\left\{ y=S_{j,0}(x)\right\} $ over the sectors $\sect{j,0}{}$
provided $r$ is sufficiently small with respect to $r'$.

\begin{thm}
\label{center_manifold} We consider a prepared family of the form
\eqref{eq_X}. There exists $\rho>0$ such that for each $\eps$ with
$||\eps||\leq\rho$ and adapted set of sectors $\sect{j,\varepsilon}{}$,
$j=0,\dots,k-1$, there exist $k$ leaves which are center manifolds
defined by graphs $\left\{ y=S_{j,\eps}(x)\right\} $ over $\sect{j,\varepsilon}{}$
and such that $\lim_{x\to x_{j,s}}S_{j,\eps}(x)=0$. In the limit
when $\eps\to0$ inside an equivalence class then $S_{j,\eps}\to S_{j,0}$
uniformly on compact sets of $\ssect{j,0}{}$. The $S_{j,\eps}$ are
unique on $\sect{j,\varepsilon}{}$ and are called \textbf{sectorial
center manifolds}. Let $W\subset\Sigma_{0}$ be a good sector. Then
the $S_{j,\eps}$ depend analytically on $\eps\in W$. Moreover the
$S_{j,\eps}$ are uniformly bounded in $\eps\in W\cup\left\{ 0\right\} $. 
\end{thm}
\begin{rem}
\label{rem:sep_estim_poly}In fact we have $S_{j,\varepsilon}=O\left(P_{\varepsilon}\right)$. 
\end{rem}
\begin{proof}
The proof is adapted from that of \cite{R2}, with ideas borrowed
from Glutsyuk \cite{G}. The idea is that the graph of the function
$S_{j,\eps}(x)$ is the stable manifold of $(x_{j,s},0)$ given in
Lemma~\ref{attractor}. The function $S_{j,\eps}(x)$ of the theorem
must be a solution of the nonlinear differential equation: \begin{equation}
P_{\eps}(x)S_{j,\eps}'(x)=S_{j,\eps}(x)(1+a(\eps)x^{k})+S_{j,\eps}^{2}(x)R_{2,\eps}(x,S_{j,\eps}(x))+P_{\eps}(x)R_{0,\eps}(x),\label{center.1}\end{equation}
 such that $S_{j,\eps}(x_{j})=0$. For $\eps=0$ the solution of \eqref{center.1}
is $k$-summable in all directions except in the directions $\exp(\frac{2\pi i\ell}{k})\mathbb{R}_{\geq0}$
for $\ell\in\ww{Z}/k$, see \cite{MR1}. If $r$ is chosen sufficiently
small the equation \eqref{center.1} with $\eps=0$ has a solution
over $\sect{j,0}{}$ for each $j=0,\dots,k-1$. We can always suppose
that $r$ is sufficiently small so that $|S_{j,0}(x)|<|x|$ for $|x|=r$
(this comes from the fact that $S_{j,0}(x)$ has an asymptotic expansion
of the form $O(x^{k+1})$ near $x=0$).

The equation \eqref{center.1} has an analytic solution defined in
the neighborhood of $x_{j,s}$ and vanishing at $x_{j,s}$ (because
the quotient of eigenvalues is neither zero nor a positive real number).
For $\eps$ sufficiently small in an equivalence class we now need
to extend this solution to $\sect{j,\varepsilon}{}$. For $(x,\eps)$
sufficiently small the inequality $|\dot{y}|>|\dot{x}|$ is satisfied
for $(x,y)$ in the cones: $K_{\ell}(\eps)=\left\{ (x,y)\,:\,|y|>|x-x_{\ell}|\right\} $,
$\ell=0,\dots k-1$. Also leaves of the foliation of \eqref{center.1}
contain trajectories with real time of all systems of the form $X_{\eps}(\theta)=e^{i\theta}X_{\eps}.$

We need to find points $(x',S_{j,\eps}(x'))$, with $|x'|=r$, which
{}``should'' belong to the center manifold and are located under
the cones $K_{\ell}(\eps)$. The extension of their trajectories under
the different $v_{\eps}(\theta)$ will yield the full center manifold.
The details are as follows.

We let $x'=r\exp(\frac{\pi i(2j+1)}{k})$. Let $\Phi_{0}^{t}$ be
the flow of $X_{0}$. Then for all $(x'',S_{j,0}(x''))$ with $|x''|=r$
and $x''\in\overline{\sect{j,0}{}}$ there exists $t(x'')\in\mathbb{C}$
such that $(x'',S_{j,0}(x''))=\Phi_{0}^{t(x'')}(x',S_{j,0}(x'))$.

Let $\eta>0$ small. The trajectories with real time starting at $(x',y)$
with $|y-S_{j,0}(x')|=\eta$, \emph{i.e} on a circle $B$, cross the
cylinder $C$ given by $|y|=r'$ along a non-contractible loop $\gamma$:
this yields a continuous map $\Pi_{0}$ from the circle $B$ to the
cylinder $C$.

We limit ourselves to values of $\eps$ with $||\eps||\leq\rho$ where
$\rho$ is sufficiently small so that $x_{j}$ remain inside $|x|<r$.
For small $\eps$ the map $\Pi_{0}$ is deformed to a continuous map
$\Pi_{\eps}$ from the circle $B$ to the cylinder $C$. Hence there
is a topological obstruction to the continuous extension of $\Pi_{\eps}$
to the disk $D=\left\{ \left(x,y\right)\,:\, x=x',\,|y-S_{j,0}(x)|\leq\eta\right\} $
given by the interior of $B$ inside the section $\{x=x'\}$, yielding
that the orbit of at least one point $(x',y_{\eps}')$ of $D$ does
not meet the cylinder.

Then the forward trajectory of $(x',y_{\eps}')$ {}``remains under''
the cones $K_{\ell}$, and in particular lies in the region $|y|<|x-x_{j,s}|$.

For all $x''$ with $|x''|=r$ and $x''\in\ov V_{j,\eps}$ there exists
$t_{\eps}(x'')$ such that $\Phi_{\eps}^{t_{\eps}(x'')}(x',y_{\eps}')=(x'',y_{\eps}'')$.
We let $S_{j,\eps}(x'')=y_{\eps}''$. When $\eps$ is small the map
$S_{j,\eps}$ is close to $S_{j,0}$ on $\{|x|=r\}\cap\ov V_{j,0}$.
In particular, if $\eps$ is sufficiently small we have $|y_{\eps}''|<|x-x_{\ell}|$
for all $\ell$.

We limit ourselves to values of $\eps$ in a good sector. Hence, if
$\theta$ is a good angle, all trajectories of $\exp(i\theta)X_{\eps}$
starting at points $(x'',y_{\eps}'')$ belong to the invariant manifold
of $(x_{j,s},0)$, \emph{i.e.} give an extension of $S_{j,\eps}(x)$.

\begin{figure}
\includegraphics[width=7cm]{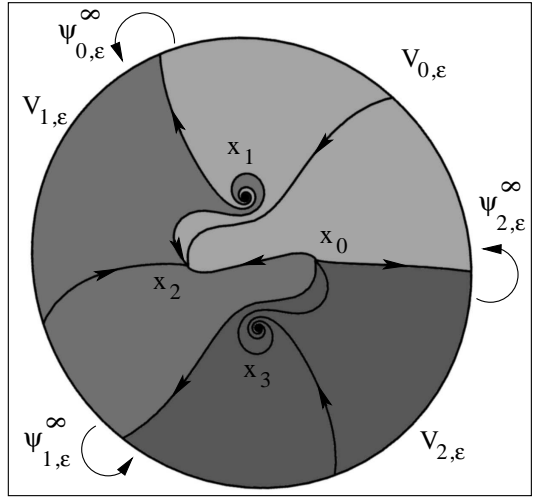}

\caption{\label{fig:center}An example of the sectors $\sect{j,\varepsilon}{}$
together with the node part of the modulus $\psi_{j,\varepsilon}^{\infty}$
when $k=3$.} 
\end{figure}

The uniform boundedness of the $S_{j,\eps}$ comes from the fact that
all graphs of functions $S_{j,\eps}$ over $V_{j,\eps}$ are located
below the cones $K_{j,s}(\eps)$. 
\end{proof}
\begin{rem}
Although the $k$ center manifolds seem to be attached to the attracting
parts of $\partial r\mathbb{D}$, they are only unique when an adapted
set of squid sectors is given. Different center manifolds attached
to different sets of adapted squid sectors may not coincide near $\partial r\mathbb{D}$
(see Figure~\ref{fig_sectors}). 
\end{rem}
The Martinet-Ramis modulus for analytic classification is a $2k$-tuple
of germs of analytic maps $\mathcal{N}_{0}=\left(\psi_{0}^{\infty},\ldots,\psi_{k-1}^{\infty},\phi_{0}^{0},\dots,\phi_{k-1}^{0}\right)$,
the $\psi_{j}^{\infty}$ being affine maps. These $k$ affine maps
will unfold as affine maps $\psi_{j,\eps}^{\infty}$ which will measure
the shift between the $k$ center manifolds: to do this we will need
to introduce adequate coordinates on which to define the $\psi_{j,\eps}^{\infty}$.
These coordinates will parameterize the space of leaves over the neighborhoods
$\sect{j,\varepsilon}{\pm}$. In particular the $\psi_{j,\eps}^{\infty}$
will all be linear when the $k$ sectorial center manifolds glue together
as a global invariant manifold. The relative position of the $k$
center manifolds can be read precisely from the $\psi_{j,\eps}^{\infty}$.

\begin{example}
\label{ex_center_man} Let us interpret the $\psi_{j,\eps}^{\infty}$
in Figure~\ref{fig:center}. The points $(x_{2},0)$ and $(x_{3},0)$
have stable manifolds. The points $(x_{0},0)$ and $(x_{1},0)$ may
have weak invariant manifolds if the quotient of their eigenvalues
is not in $1/\ww{N}_{\neq0}$. 
\begin{enumerate}
\item $\psi_{0,\eps}^{\infty}$ measures if the stable manifold of $(x_{2},0)$
is ramified at $(x_{1},0)$: it is indeed the case if $\psi_{0,\eps}^{\infty}$
is nonlinear. In that case, if $(x_{1},0)$ has a weak invariant manifold,
then necessarily it does not coincide with the stable manifold of
$(x_{2},0)$. 
\item $\psi_{1,\eps}^{\infty}$ measures if the stable manifolds of $(x_{2},0)$
and $(x_{3},0)$ coincide or not: they coincide precisely if $\psi_{1,\eps}^{\infty}$
is linear. 
\item $\psi_{2,\eps}^{\infty}$ measures if the stable manifolds of $(x_{2},0)$
and $(x_{3},0)$ coincide on the other side of $(x_{0},0)$. 
\item From (1) and (2) it is possible to decide if the stable manifold of
$(x_{2},0)$ coincides with the weak invariant manifold of $(x_{0},0)$.
Indeed if $\psi_{1,\eps}^{\infty}$ (\emph{resp.} $\psi_{2,\eps}^{\infty}$)
is linear and $\psi_{2,\eps}^{\infty}$ (\emph{resp.} $\psi_{1,\eps}^{\infty}$)
is nonlinear, then necessarily the stable manifold of $(x_{2},0)$
does not coincide with the weak invariant manifold of $(x_{0},0)$.
In the particular case where $(x_{0},0)$ would be a resonant node
this would imply that it would have no weak invariant manifold (this
is the \textbf{parametric resurgence phenomenon} described in \cite{R2}). 
\item It is also possible to decide directly if the stable manifold of $(x_{2},0)$
coincides with the weak invariant manifold of $(x_{0},0)$ even if
both $\psi_{0,\eps}^{\infty}$ and $\psi_{1,\eps}^{\infty}$ are nonlinear,
but this is more involved and requires two additional tools: the Lavaurs
maps and the other part of the modulus, namely the functions $\phi_{j,\eps}^{0}$.
Indeed we need a characterization of the weak invariant manifold of
$(x_{0},0)$: it is a leaf which is fixed when one turns around $(x_{0},0)$.
Following the leaves when one turns around a singular point requires
the transition maps between the space of leaves associated to the
different $\sect{j,\varepsilon}{\pm}$ over the sectors $\sect{j,\sigma(j),\varepsilon}{g}$.
These transition maps, called Lavaurs maps, are global linear maps.
The maps $\phi_{j,\eps}^{0}$ are related to transition maps between
the space of leaves associated to the different $\sect{j,\varepsilon}{\pm}$
over the sectors $\sect{j,\varepsilon}{s}$. We come back to this
in Section~\ref{subsec:reading}, Example~\ref{ex:11_1}. 
\end{enumerate}
\end{example}
\begin{cor}
\label{cor:red_center_man} We consider a prepared family \eqref{eq_X}
and $W\subset\Sigma_{0}$ a good sector on which the conclusion of
Theorem~\ref{center_manifold} is satisfied. The family of changes
of coordinates $\left(x,y\right)\mapsto\left(x,y-S_{j,\eps}(x)\right)$
transforms the family $\left(X_{\varepsilon}\right)_{\varepsilon}$
into \begin{equation}
\left(X_{j,\eps}\right)_{\varepsilon}:=\left(P_{\eps}(x)\pp{x}+y\left(1+a(\eps)x^{k}+\tilde{R}_{j,\eps}(x,y)\right)\pp{y}\right)_{\varepsilon}\label{unfold_sect}\end{equation}
 with $\tilde{R}_{j,\eps}=O(y)$ over $\sect{j,\varepsilon}{}$ and
uniformly in $\varepsilon\in W\cup\left\{ 0\right\} $. 
\end{cor}

\section{\label{sec: asymptotic}Asymptotic paths}

We want to show that the foliation $\fol{j,\varepsilon}^{\pm}$ induced
by $Z_{\varepsilon}$ over {}``canonical sectors'' is {}``trivial''
in some way, the triviality being expressed in terms of asymptotic
cycles. This property will ensure that $Z_{\varepsilon}$ is analytically
conjugate to the model over these sectors, as explained in Section~\ref{sec:cohomological}.

In order to build the canonical sectors we first give some definitions.

\subsection{Basic definitions }

Throughout this section $\overline{\Omega}$ stands for the topological
closure of $\Omega$.

\begin{defn}
\label{def:homol-asy}Let $Z$ be a vector field with components holomorphic
on a neighborhood of $\overline{\Omega}$ for some open set $\Omega\subset\ww{C}^{2}$
and consider the foliation $\fol{}$ induced by $Z$ on $\Omega$. 
\begin{enumerate}
\item A piecewise-$C^{1}$ map $\gamma\,:\,\ww{R}\to\overline{\Omega}$
satisfying

\begin{enumerate}
\item there exists a leaf $\mathcal{L}$ or a singular point $\mathcal{L}:=\left\{ S\right\} $
such that: $\left(\forall t\in\ww{R}\right)\,\,\gamma\left(t\right)\in\mathcal{L}$ 
\item $\lim_{t\to\pm\infty}\gamma\left(t\right)=p_{\pm}\in\overline{\Omega}$ 
\end{enumerate}
is called an \textbf{asymptotic path}, linking $p_{-}$ to $p_{+}$
within $\fol{\mathcal{}}$. These points need not belong to the same
leaf (they can be singularities of $\fol{}$) and are called the \textbf{endpoints}
of $\gamma$. They will be referred to as $\gamma\left(\pm\infty\right)$.
The map $\gamma$ will be called an \textbf{asymptotic cycle} if $p_{-}=p_{+}$.

\item A piecewise-$C^{1}$ map $h\,:\,\ww{R}\times\ww{R}\to\overline{\Omega}$
such that

\begin{enumerate}
\item $\left(\forall t\in\ww{R}\cup\left\{ \pm\infty\right\} \right)\,\, h\left(t,\cdot\right),\, h\left(\cdot,t\right)$
are asymptotic paths 
\item the family $\left(h\left(t,\cdot\right)\right)_{t\in\ww{R}}$ (\emph{resp.}
$\left(h\left(\cdot,t\right)\right)_{t\in\ww{R}}$) converges uniformly
to $h\left(\pm\infty,\cdot\right)$ (\emph{resp.} $h\left(\cdot,\pm\infty\right)$)
as $t\to\pm\infty$ 
\end{enumerate}
is called an \textbf{asymptotic homology} between $\gamma_{-\infty}:=h\left(-\infty,\cdot\right)$
and $\gamma_{+\infty}:=h\left(+\infty,\cdot\right)$.

\end{enumerate}
\end{defn}
This notion will be useful to express the triviality of $\fol{}$
over the canonical sectors. In fact one could give a definition of
what may be called {}``asymptotic homology of $\fol{}$ over $\Omega$''
by considering the complex of $\ww{Z}$-modules generated by points
(0-chains), asymptotic paths (1-chains) and asymptotic homologies
(2-chains) endowed with boundary operators. We will not need these
refinements in our present study but that is what is at work here.

\begin{defn}
\label{def:connex-asy}~ 
\begin{enumerate}
\item Let $p\in\overline{\Omega}$. We define the \textbf{connected component}
of $p$ in $\fol{}$ as the set\begin{equation}
\left\{ q\in\overline{\Omega}\,:\,\left(\exists\gamma\mbox{ an asymptotic path}\right)\,\gamma\left(-\infty\right)=p\mbox{ and }\gamma\left(+\infty\right)=q\right\} \,.\end{equation}
 We say that $\fol{}$ is \textbf{connected} when there exists a point
$p\in\ov{\Omega}$ such that all points of $\overline{\Omega}$ belong
to the connected component of $p$. 
\item We say that $\fol{}$ is \textbf{simply connected} when each asymptotic
cycle is asymptotically homologous to a constant path. 
\item A foliation $\fol{}$ both connected and simply connected will be
called (\textbf{asymptotically}) \textbf{trivial.} 
\end{enumerate}
\end{defn}
\begin{rem}
We will show below that the foliation over each squid sector is connected.
The point $p$ we will choose will be the point of node type in the
closure of the squid sector. Remark that the connected component $p'$
of an interior point of the squid sector will only be the closure
of the leaf through that point. 
\end{rem}

\subsection{Canonical sectors}

After these preparations we shall prove the following:

\begin{thm}
\label{thm:secto-trivial} We consider an adapted set of squid sectors
$\sect{j,\varepsilon}{\pm}$ covering $r\mathbb{D}$ in $x$-space
where $\varepsilon$ belongs to some good sector $W\subset\Sigma_{0}$.
Each $\sect{j,\varepsilon}{\pm}$ is adherent to two points $x_{j,s}$
and $x_{j,n}$. Let $\ssect{j,\varepsilon}{\pm}$ be the interior
of the connected component of $p_{j,n}=(x_{j,n},0)$ in the foliation
induced by $Z_{\varepsilon}$ over $\sect{j,\varepsilon}{\pm}\times r'\ww{D}$.
There exist $r,r',\rho>0$ such that the following assertions hold
for $\varepsilon\in W\cup\left\{ 0\right\} $ : 
\begin{enumerate}
\item For each $p\in\mathcal{V}_{j,\varepsilon}^{\pm}$ there exists an
asymptotic path $\gamma_{j,\varepsilon}^{\pm}\left(p\right)$ within
$\fol{j,\varepsilon}^{\pm}$ such that $\gamma_{j,\varepsilon}^{\pm}\left(-\infty\right)=p_{j,n}$
and $\gamma_{j,\varepsilon}^{\pm}\left(t\right)=p$ for all $t\geq0$. 
\item The domain $\mathcal{V}_{j,\varepsilon}^{\pm}$ contains a fibered
squid sector $\sect{j,\varepsilon}{\pm}\times r''\ww{D}$. We denote
$\fol{j,\varepsilon}^{\pm}$ the foliation induced by $Z_{\varepsilon}$
over $\ssect{j,\varepsilon}{\pm}$. 
\item There exists a unique leaf $\mathcal{S}_{j,\varepsilon}^{\pm}$ of
$\fol{j,\varepsilon}^{\pm}$ accumulating on both $p_{j,n}$ and $p_{j,s}$
corresponding to the sectorial center manifold of $Z_{j,\varepsilon}$.
This leaf is the graph of a holomorphic function \begin{eqnarray}
S_{j,\varepsilon}^{\pm}\,:\,\sect{j,\varepsilon}{\pm} & \to & r'\ww{D}\end{eqnarray}
 which extends as a continuous function on the closure $S_{j,\varepsilon}^{\pm}\left(x_{j,n}\right)=S_{j,\varepsilon}^{\pm}\left(x_{j,s}\right)=0$.
The sectorial central manifolds $\mathcal{S}_{j,\varepsilon}^{\pm}$
glue on $\ssect{j,\varepsilon}{}=\ssect{j,\varepsilon}{+}\cup\ssect{j,\varepsilon}{-}$
and coincide with the graph of $x\mapsto S_{j,\varepsilon}\left(x\right)$
(see Theorem~\ref{center_manifold}). (Of course in the case $k=1$
the gluing only occurs on the saddle and gate sides and $\ssect{j,\varepsilon}{}=\ssect{j,\varepsilon}{+}\cup\ssect{j,\varepsilon}{-}$
is a sector of opening greater than $2\pi$.) 
\item The foliation $\fol{j,\varepsilon}^{\pm}$ is asymptotically trivial. 
\end{enumerate}
\end{thm}
We postpone the proof of Theorem~\ref{thm:secto-trivial} till Section~\ref{proof_4.4}.

\begin{defn}
\label{def:canonical-sector} The $2k$ sectors $\ssect{j,\varepsilon}{\pm}$
are called the \textbf{canonical sectors} associated to $Z_{\varepsilon}$. 
\end{defn}
For a good sector $W\subset\Sigma_{0}$ and $\varepsilon\in W\cup\left\{ 0\right\} $
we let \begin{equation}
\ssect{\varepsilon}{}:=\mbox{int}\left(\mbox{clos}\left(\cup_{j=0}^{k-1}\left(\mathcal{V}_{j,\varepsilon}^{+}\cup\mathcal{V}_{j,\varepsilon}^{-}\right)\right)\right).\label{mathcalV_0}\end{equation}
 This is an open neighborhood of $\left(0,0\right)$ containing a
polydisk $r\ww{D}\times r''\ww{D}$ independent of $\varepsilon$.

\subsection{\label{proof_4.4}Proof of Theorem~\ref{thm:secto-trivial}}

\begin{figure}
\includegraphics[width=8cm]{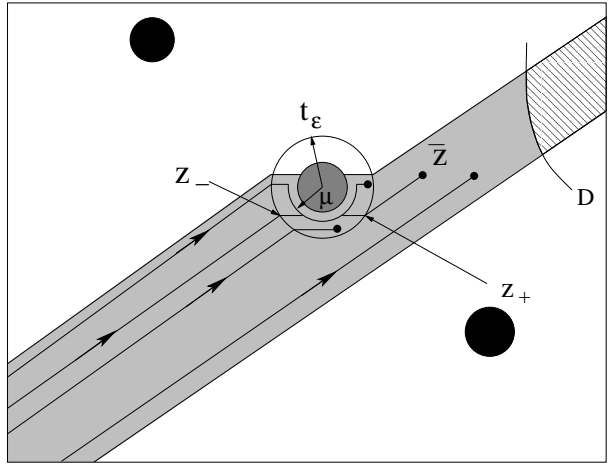}

\caption{\label{fig:asy_path_z_coord}Some asymptotic paths in $z$-coordinate.
The curve $D$ indicates where the $y$-component of the leaf reaches
$r'\ww{S}^{1}$.} 
\end{figure}

We fix a good angle $\theta$ associated to a good sector $W$, see
Lemma~\ref{good_covering}. Note that \textbf{(3)} has been proved
in Theorem~\ref{center_manifold}, so we can apply the change of
coordinates $\left(x,y\right)\mapsto\left(x,y-S_{j,\varepsilon}\left(x\right)\right)$,
where $x\mapsto S_{j,\varepsilon}\left(x\right)$ is the sectorial
central manifold over $\sect{j,\varepsilon}{+}\cup\sect{j,\varepsilon}{-}$,
which sends the foliation $\fol{j,\varepsilon}^{\pm}$ on $\fol{}'$
defined by \begin{eqnarray}
\tilde{X}_{\varepsilon}\left(x,y\right) & := & P_{\varepsilon}\left(x\right)\pp{x}+y\left(1+a(\eps)x^{k}+\tilde{R}_{\varepsilon}\left(x,y\right)\right)\frac{\partial}{\partial y}\end{eqnarray}
 on $\sect{j,\varepsilon}{\pm}\times r'\ww{D}$ (see Corollary \ref{cor:red_center_man}).
We will prove the remaining claims \textbf{(1)}, \textbf{(2)} and
\textbf{(4)} for that foliation which, after possibly decreasing $r,r',\rho>0$,
will still hold back in the original coordinates. The proof will rely
on the following straightforward estimate which we give without proof
:

\begin{lem}
\label{lem:modulus_estimate}Let $\chi\left(t\right):=\left|y\left(t\right)\right|$
in the following non-autonomous system:\begin{eqnarray}
\dot{x}\left(t\right) & = & \exp(i\theta)P_{\varepsilon}\left(x(t)\right)\\
\dot{y}\left(t\right) & = & \exp(i\theta)y(t)\left(1+\tilde{R}_{\varepsilon}\left(x(t),y(t)\right)\right),\label{eq:syst_diff_path}\end{eqnarray}
 so that\begin{eqnarray}
\dot{\chi}\left(t\right) & = & \chi\left(t\right)Re\left(\exp(i\theta)\left(1+\tilde{R}_{\varepsilon}\left(x\left(t\right),y\left(t\right)\right)\right)\right).\label{eq_chi}\end{eqnarray}

\begin{enumerate}
\item Assume that for some $r,r'$ we have \begin{eqnarray}
0<\alpha\leq & Re\left(\exp(i\theta)(1+\tilde{R}\left(x,y\right))\right) & \leq\beta\end{eqnarray}
 for any $\left(x,y\right)\in r\ww{D}\times r'\ww{D}$. Then for any
$t\leq0$ : \begin{equation}
\chi(0)e^{\beta t}\leq\chi(t)\leq\chi(0)e^{\alpha t}\label{eq:module-path}\end{equation}

\item It is possible to find $r,r'$ small enough so that $\alpha$ and
$\beta$ are as close to $\cos\theta$ as we wish, independently on
$\varepsilon$. 
\end{enumerate}
\end{lem}
\textbf{(1)} First we build the path in $z$-coordinate (the function
$z$ is defined in Lemma~\ref{P_2}) and we refer to the notations
given in Figure~\ref{fig:asy_path_z_coord}. We define $t_{\varepsilon}:=\kappa\left|\left|\varepsilon\right|\right|^{-k}$
as in Definition~\ref{def:squid-sectors}.

\begin{itemize}
\item If $\overline{z}:=z\left(\overline{x}\right)$ belongs to the part
of the strip which can be linked to $Im\left(z\right)=-\infty$ in
a straight line of slope $\theta$ we define \begin{eqnarray*}
z\left(t\right) & := & \overline{z}+te^{i\theta}\end{eqnarray*}
 for $t\leq0$. 
\item If $\overline{z}$ belongs to the disk $t_{\varepsilon}\ww{D}$ we
choose a path $t\mapsto z\left(t\right)$ avoiding the central hole
$B_{0}$ and reaching $z_{-}$ on the boundary of the disk. The path
consists of horizontal line(s) and possibly an arc of the circle of
fixed radius $\mu=Ar^{-k}$. We then link $z_{-}$ to $-\infty$ with
a straight line as before. 
\item Otherwise we link $\overline{z}$ to some point $z_{+}$ of the circle
$t_{\varepsilon}\ww{S}^{1}$ with a straight line of slope $\theta$,
then proceed as above starting from $z_{+}$. 
\end{itemize}

We call $\gamma$ the path $t\mapsto x\left(t\right)$ corresponding
to the path in $z$-coordinate build just above. We need to show the
existence of $r''>0$ such that $\gamma$ can be lifted into the foliation
for all $\left(\overline{x},\overline{y}\right)\in\sect{j,\varepsilon}{\pm}\times r''\ww{D}$.
According to the previous lemma, if $\left|z\left(t\right)\right|>t_{\varepsilon}$
or if $z\left(t\right)$ belongs to a horizontal line then $t\mapsto\left|y\left(t\right)\right|$
is decreasing, so that these parts of $\gamma$ can be lifted in $\sect{j,\varepsilon}{\pm}\times r'\ww{D}$
as soon as $\overline{y}\in r'\ww{D}$. In all the other cases the
paths used are of bounded length, independently of $\varepsilon$.
Hence using finitely many flow-boxes we derive the existence of $r''>0$
satisfying the expected property.

\textbf{(2)} The claim is proved through the Lemma~\ref{lem:modulus_estimate}.

\textbf{(3)} This comes from Theorem~\ref{center_manifold}.

\textbf{(4)} So far we have proved that $\fol{}'$ is connected. We
now show that it is simply connected. In fact we prove a slightly
stronger result :

\begin{prop}
\label{pro:leaf-is-graph} There exists $r>0$ independent of small
$\varepsilon$ such that: 
\begin{enumerate}
\item The leaf of $\fol{j,\varepsilon}^{\pm}$ passing through $\left(\overline{x},\overline{y}\right)$
is the graph of a holomorphic function $x\in\Omega_{\ov{x},\ov{y}}\mapsto K_{j,\varepsilon,\ov{x},\ov{y}}^{\pm}\left(x\right)$
and $\Omega_{\ov{x},\ov{y}}\subset V_{j,\eps}^{\pm}$ is simply connected.
Let $\Omega=\bigcup_{(\ov{x},\ov{y})\in\mathcal{V}{j,\eps}^{\pm}}\Omega_{\ov{x},\ov{y}}\times\{(\ov{x},\ov{y})\}\subset V_{j,\eps}^{\pm}\times\ssect{j,\varepsilon}{\pm}$.
Then there exists a holomorphic function $K_{j,\varepsilon}^{\pm}\,:\Omega\to\ww{C}$
such that $K_{j,\eps}^{\pm}(x,\ov{x},\ov{y})=K_{j,\eps,\ov{x},\ov{y}}^{\pm}(x)$
. 
\item The closure of $\Omega_{\ov{x},\ov{y}}$ is also simply connected. 
\item The closure $D$ of $\left[K_{j,\varepsilon}^{\pm}\left(\cdot,\overline{x},\overline{y}\right)\right]^{-1}\left(r'\ww{S}^{1}\right)$
is a (connected) real analytic curve which separates $\sect{j,\varepsilon}{\pm}$
into two connected components, and cuts the boundary $\partial\sect{j,\varepsilon}{\pm}$
in exactly two points (see Figure~\ref{fig:squid_modulus}(a) and
Figure~\ref{fig:squid_intersect_trivial}). One connected component
of $\sect{j,\varepsilon}{\pm}\backslash D$ accumulates on $x_{j,n}$
while the other accumulates on $x_{j,s}$. Besides $D\to\left\{ x_{j,n}\right\} $
as $r'\to0$. 
\end{enumerate}
\end{prop}
\begin{proof}
The proof is done as in Proposition~\ref{pro:first-integ-model}(3),
using the estimates of Lemma~\ref{lem:modulus_estimate}. Indeed
the only obstruction to the analytic continuation of a leaf is the
constraint $\left|y\left(x\right)\right|<r'$ because the foliation
is transverse to the lines $\left\{ x=\mbox{cst}\right\} $. 
\end{proof}
Because $\Omega_{\ov{x},\ov{y}}$ is simply connected the endpoint
of a candidate asymptotic cycle $\gamma$ of $\fol{j,\varepsilon}^{\pm}$
must be a singular point. The leaf cannot accumulate on $p_{j,s}$
so that $\gamma\left(\pm\infty\right)=p_{j,n}$. In that case $\gamma$
is asymptotically homologous to its endpoint as the closure of $\Omega_{\ov{x},\ov{y}}$
is simply connected. Thus the sectorial foliation is simply connected.
\hfill{}$\Box$

\subsection{Asymptotic homology over the intersections of sectors}

\begin{figure}
\subfigure[in
$x$-coordinate]{\includegraphics[width=7cm]{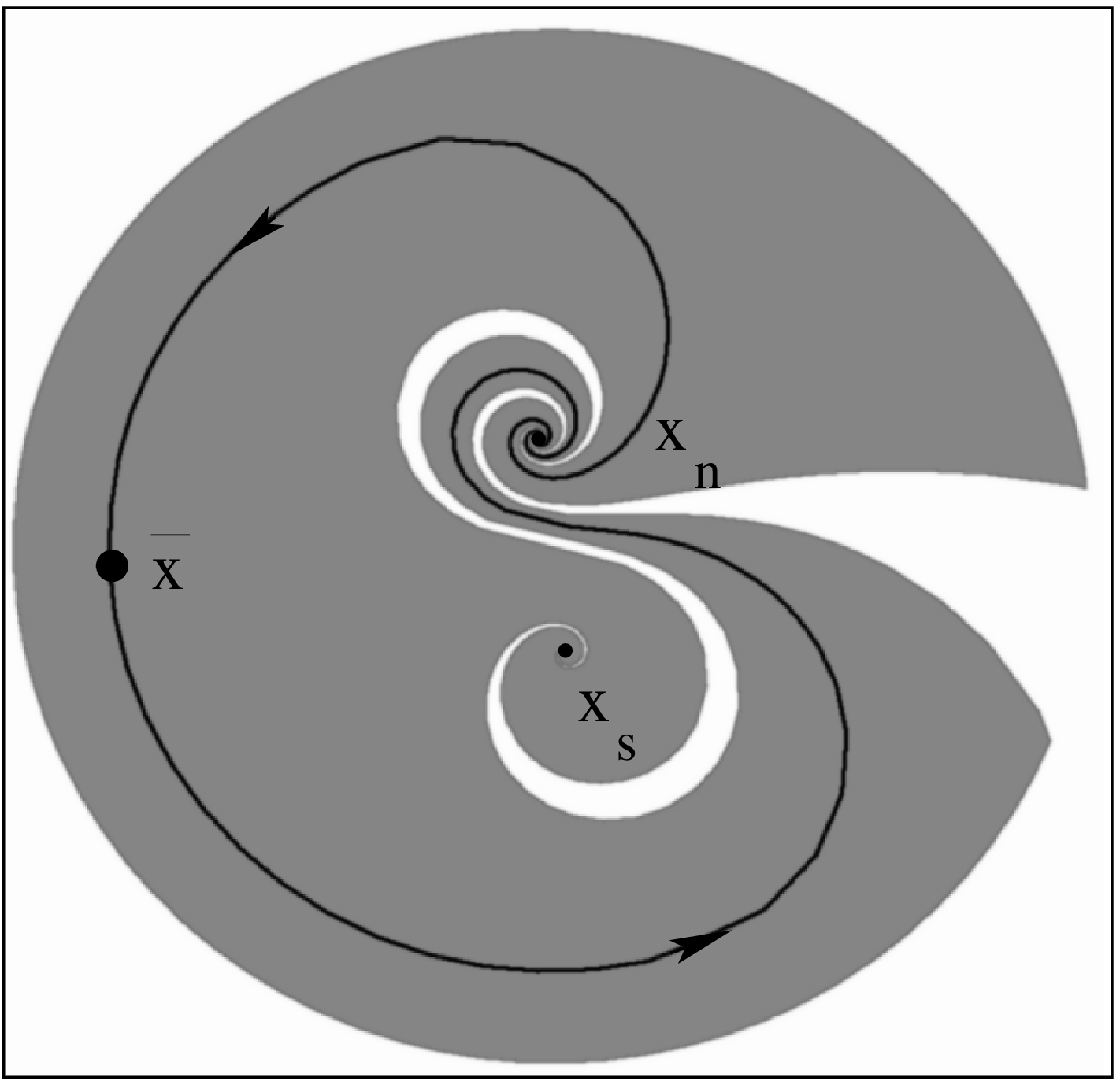}} \hfill{}\subfigure[lifted in the
leaf]{\includegraphics[width=7cm]{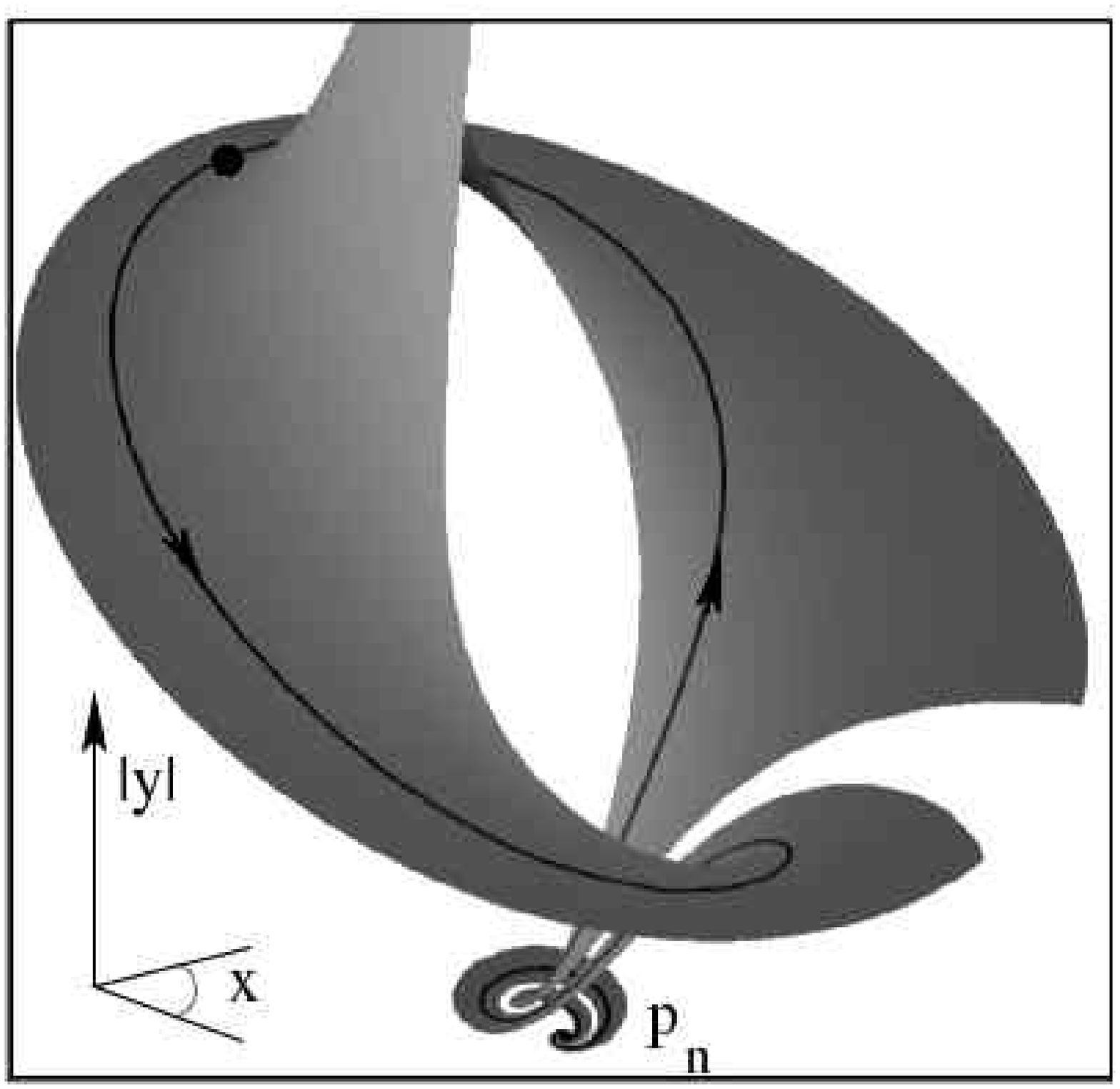}}

\caption{\label{fig:asy_path_2}An asymptotic cycle ($k=1$). } 
\end{figure}

\begin{figure}
\hfill{}\subfigure[$\gamma^n \subset \mathcal{V}_{j,\eps}^n$ is homologous to its endpoints]{\includegraphics[width=7cm]{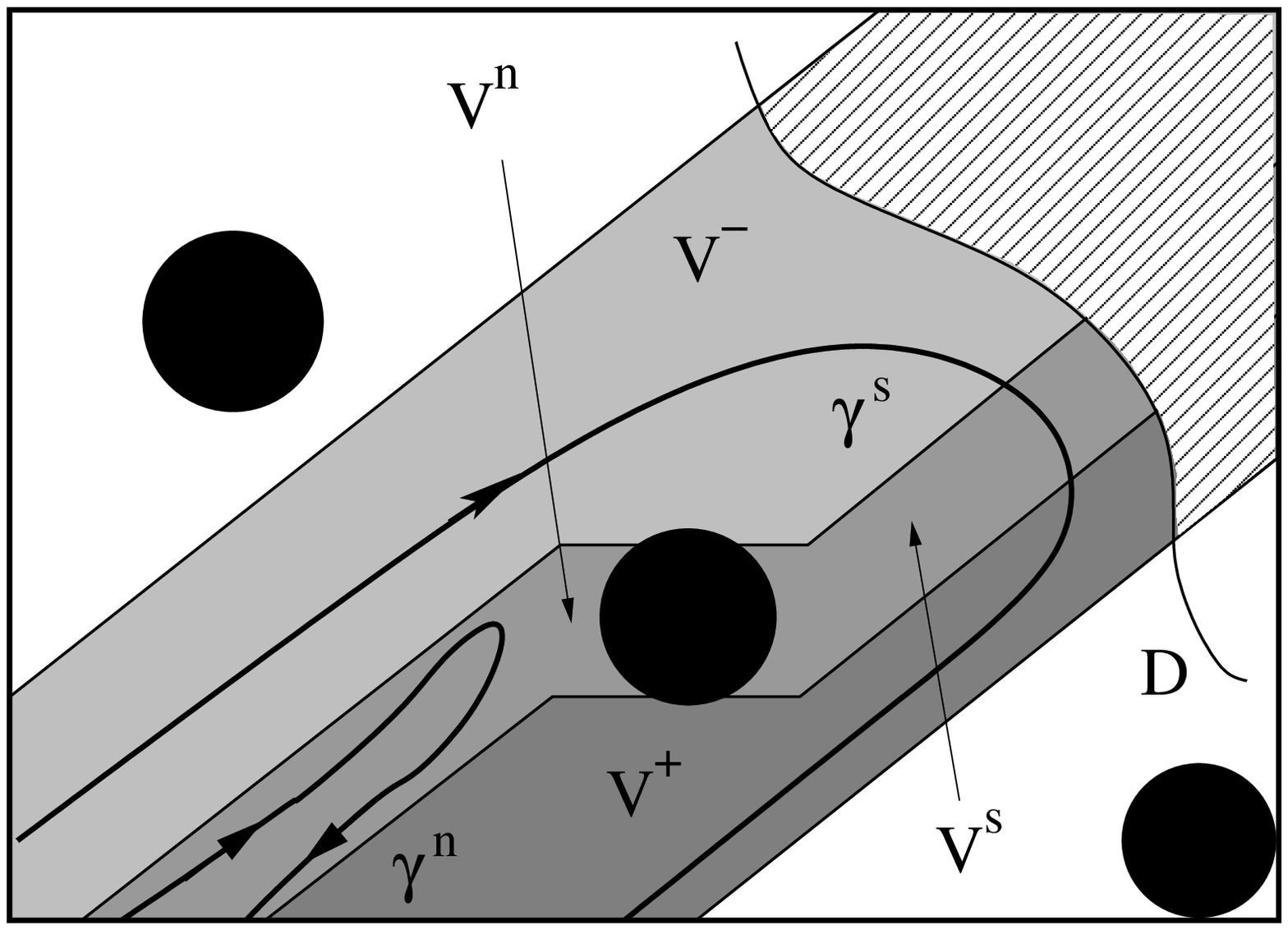}}\hfill{}\subfigure[$\gamma^g \subset \mathcal{V}_{j,\eps}^g$
is homologous to its endpoints]{\includegraphics[width=7cm]{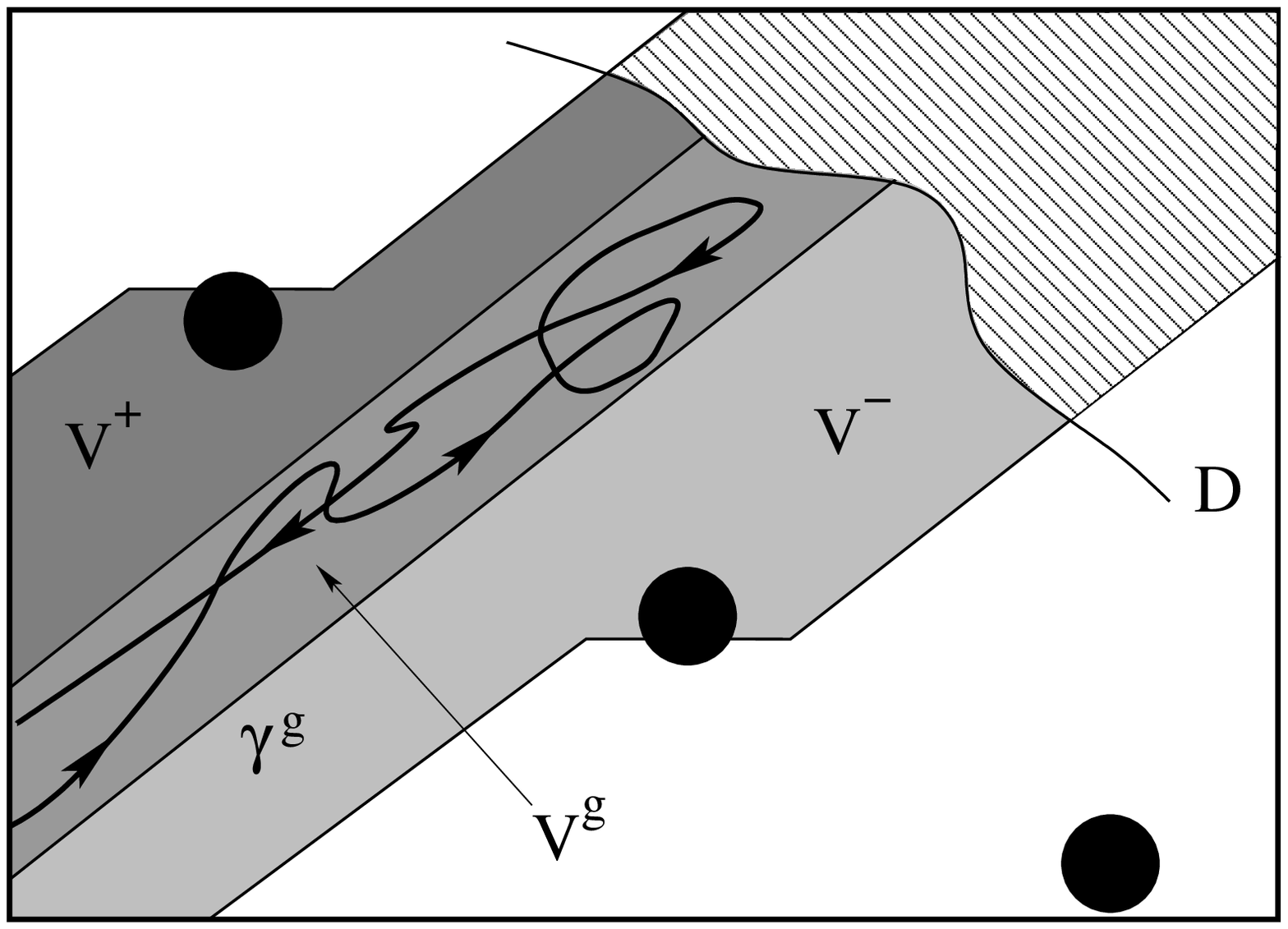}}\hfill{}

\caption{\label{fig:squid_intersect_trivial} The asymptotic cycle $\gamma^{s}=\gamma_{j,\varepsilon}^{s}\left(p\right)$
is not trivial in $\ssect{j,\varepsilon}{+}\cup\ssect{j,\varepsilon}{-}$
when $p\in\ssect{j,\varepsilon}{s}$. We recall that the $z$-coordinate
is a $k$-sheeted covering of $\ww{C}$ minus the black disks so that
$\sect{j,\varepsilon}{s}$ and $\sect{j,\varepsilon}{n}$ may not
actually belong to the same sheet. In particular the two endpoints
of $\gamma^{s}$ may correspond to different points $p_{j,n}^{\pm}$
of node type.} 
\end{figure}

\begin{defn}
\label{def:asy_path} We show below that the intersections of different
$\ssect{j,\varepsilon}{\pm}$ correspond to sectors $\ssect{j,\varepsilon}{s}$,
$\ssect{j,\varepsilon}{n}$ and $\ssect{j,\sigma(j),\varepsilon}{g}$
as in the case of squid sectors (Lemma~\ref{lem:intersect}). Each
$p$ in such an intersection yields an asymptotic path : 
\begin{itemize}
\item $\gamma_{j,\varepsilon}^{s}\left(p\right):=\gamma_{j,\varepsilon}^{+}\left(p\right)-\gamma_{j,\varepsilon}^{-}\left(p\right)$
if $p\in\ssect{j,\varepsilon}{s}$. 
\item $\gamma_{j,\varepsilon}^{n}\left(p\right):=\gamma_{j+1,\varepsilon}^{+}\left(p\right)-\gamma_{j,\varepsilon}^{-}\left(p\right)$
if $p\in\ssect{j,\varepsilon}{n}$. 
\item $\gamma_{j,\sigma(j),\varepsilon}^{g}$$\left(p\right):=\gamma_{j,\varepsilon}^{+}\left(p\right)-\gamma_{\sigma(j),\varepsilon}^{-}\left(p\right)$
if $p\in\ssect{j,\sigma(j),\varepsilon}{g}$. 
\end{itemize}
This path links the corresponding node-type singular points within
the foliation induced on $\mathcal{V}_{j,\varepsilon}^{+}\cup\ssect{j,\varepsilon}{-}$,
$\mathcal{V}_{j+1,\varepsilon}^{+}\cup\ssect{j,\varepsilon}{-}$ or
$\ssect{j,\varepsilon}{+}\cup\ssect{\sigma(j),\varepsilon}{-}$. We
call it the \textbf{canonical asymptotic path} associated to $p$. 
\end{defn}
\begin{prop}
\label{pro:intersect_secto-trivial}Under the same hypothesis as Theorem~\ref{thm:secto-trivial}
the following assertions hold. 
\begin{enumerate}
\item For $k>1$ the intersections $\mathcal{V}_{j,\varepsilon}^{+}\cap\mathcal{V}_{j,\varepsilon}^{-}$
splits into one or two connected components, namely $\ssect{j,\varepsilon}{s}$
(and possibly $\ssect{j,\sigma\left(j\right),\varepsilon}{g}$ if
$\varepsilon\neq0$ and both sectors are adherent to the same singular
points $p_{j,n}$ and $p_{j,s}$). Similarly the intersections $\mathcal{V}_{j,\varepsilon}^{-}\cap\mathcal{V}_{j+1,\varepsilon}^{+}$
splits into $\mathcal{V}_{j,\varepsilon}^{n}$ (and possibly $\mathcal{V}_{j+1,\sigma\left(j+1\right),\varepsilon}^{g}$
if $\varepsilon\neq0$ and both sectors are adherent to the same singular
points $p_{j,n}$ and $p_{j,s}$). Also $\mathcal{V}_{j,\sigma(j),\varepsilon}^{g}=\mathcal{V}_{j,\varepsilon}^{+}\cap\mathcal{V}_{\sigma(j),\varepsilon}^{-}$
if $\sigma(j)\neq j,j+1$.

In the case $k=1$ the intersection of the two sectors splits into
two or three components.

\item The foliation induced by $Z_{\varepsilon}$ over $\ssect{j,\varepsilon}{n}$
or $\ssect{j,\sigma\left(j\right),\varepsilon}{g}$ is trivial. Moreover
the canonical asymptotic path through a point $p$ lying in one of
those sectors is asymptotically homologous to the node-type singular
point. 
\item The foliation induced by $Z_{\varepsilon}$ over $\ssect{j+1,\varepsilon}{+}\cup\ssect{j,\varepsilon}{-}$
is asymptotically trivial. 
\item If $p\in\ssect{j,\varepsilon}{s}\backslash\mathcal{S}_{j,\varepsilon}$
then $\gamma_{j,\varepsilon}^{s}\left(p\right)$ is not homologous
to a constant path. Any other asymptotic path in the same leaf over
$\ssect{j,\varepsilon}{+}\cup\ssect{j,\varepsilon}{-}$, not homologous
to a constant path, is homologous to $\gamma_{j}^{s}\left(p\right)$
in $\ssect{j,\varepsilon}{+}\cup\ssect{j,\varepsilon}{-}$ (up to
reversing of orientation) and has same endpoints. 
\item In $\mathcal{S}_{j,\varepsilon}\cup\left\{ p_{j,s}\right\} $ any
asymptotic path is homologous to $\left\{ p_{j,n}\right\} $. 
\end{enumerate}
\end{prop}
\begin{proof}
This proof is mainly graphical. We refer to Figure~\ref{fig:squid_intersect_trivial}.

\textbf{(1)} According to Proposition~\ref{pro:leaf-is-graph} each
leaf is a graph over a domain $\Omega$ (the complement of the hatched
area in Figure~\ref{fig:squid_intersect_trivial}) and this domain
looks like a half-strip or a strip in $z$-coordinate. Hence $\ssect{j,\varepsilon}{s}$,
$\ssect{j,\varepsilon}{n}$ and $\ssect{j,\sigma(j),\varepsilon}{g}$
are connected. The remaining claims are easy from Lemma~\ref{lem:intersect}.

\textbf{(2)}--\textbf{(5)} are immediate. 
\end{proof}

\section{\label{sec:cohomological}Cohomological equations}

In order to work out the moduli of classification, we need to describe
precisely the global obstructions to solve equations of the form\begin{eqnarray}
X_{\varepsilon}\cdot F_{\varepsilon} & = & G_{\varepsilon}\label{eq:eq-homog}\end{eqnarray}
 where $G_{\varepsilon}$ is given. We will explain why it is so in
the oncoming Section~\ref{orb_classification}. Such equations are
called \textbf{cohomological equations}. A natural intuitive solution
is given by $F_{\eps}=\int G_{\eps}$, where the integral is taken
along trajectories of the vector field. In order to make this formal
we need to choose well a base point. A natural base point is the point
of node type. As it is reached in infinite time we need to define
the notion of asymptotic integrals.

\begin{defn}
Let $\Omega$ be an open set of $\ww{C}^{2}$ and consider a 1-form
$\omega$ with coefficients meromorphic on a neighborhood of the closure
of $\Omega$. Consider an asymptotic path $\gamma$ which avoids the
poles of $\omega$ (except maybe at its end-points). We define the
\textbf{asymptotic integral} of $\omega$ along $\gamma$ as \begin{eqnarray}
\int_{\gamma}\omega & := & \int_{-\infty}^{+\infty}\gamma^{*}\omega\,.\end{eqnarray}

\end{defn}
~

\begin{defn}
\label{def:holo_bound} Let \emph{$\mathcal{A}:=\left(\mathcal{A}_{\varepsilon}\right)_{\varepsilon\in W\cup\left\{ 0\right\} }:=\left(\mathcal{V}_{j,\varepsilon}^{\pm}\right)_{\varepsilon\in W\cup\left\{ 0\right\} }$}
be a family of canonical sectors associated to a good sector $W$
in $\varepsilon$-space and to a boundary sector of $r\ww{D}$ (\emph{i.e.}
$\mathcal{A}$ is one of the $2k$ families $\left(\mathcal{V}_{j,\varepsilon}^{\pm}\right)_{\varepsilon\in W\cup\left\{ 0\right\} }$).
We define the algebra $\mathcal{O}_{b}\left(\mathcal{A},W\right)$
of all functions $\left(x,y,\varepsilon\right)\mapsto G_{\varepsilon}\left(x,y\right)$
holomorphic and bounded on $\bigcup_{\varepsilon\in W}\mathcal{V}_{j,\varepsilon}^{\pm}\times\left\{ \varepsilon\right\} $
with continuous extension to the closure and such that $G_{0}$ be
holomorphic on $\mathcal{A}_{0}$. 
\end{defn}
\begin{thm}
\label{thm:secto-solve} Let $\tau_{\varepsilon}:=\frac{dx}{P_{\eps}(x)}$,
let $\mathcal{A}:=\left(\ssect{j,\varepsilon}{\pm}\right)_{\varepsilon\in W\cup\left\{ 0\right\} }$
be a family of canonical sectors where $W$ is a good sector. We suppose
that for each $\eps\in W$ the singular points of node and saddle
type in $\ssect{j,\varepsilon}{\pm}$ are given by $p_{j,n}=(x_{j,n},0)$
and $p_{j,s}=(x_{j,s},0)$ respectively. Consider a function $G\in\mathcal{O}_{b}\left(\mathcal{A},W\right)$
with \[
G_{\varepsilon}=O\left(P_{\varepsilon}\left(x\right)\right)+O\left(y\right).\]

\begin{enumerate}
\item For $p\in\ssect{j,\varepsilon}{\pm}$ define\begin{eqnarray}
F_{j,\varepsilon}^{\pm}\left(p\right) & := & \int_{\gamma_{j,\varepsilon}^{\pm}\left(p\right)}G_{\varepsilon}\,\tau_{\varepsilon}\label{def_F}\end{eqnarray}
 where $\gamma_{j,\varepsilon}^{\pm}\left(p\right)$ is constructed
in Theorem~\ref{thm:secto-trivial}. This asymptotic integral is
absolutely convergent for all $\varepsilon\in W\cup\left\{ 0\right\} $. 
\item Each function $F_{j,\varepsilon}^{\pm}$ defined in \eqref{def_F}
is holomorphic on $\ssect{j,\varepsilon}{\pm}$ and satisfies \begin{eqnarray}
X_{\varepsilon}\cdot F_{j,\varepsilon}^{\pm} & = & G_{\varepsilon}.\end{eqnarray}

\item The function extends continuously to $\ssect{j,\varepsilon}{\pm}\cup\left\{ p_{s},p_{n}\right\} \times r'\ww{D}$.
In that case the function $y\mapsto F_{j,\varepsilon}^{\pm}\left(x_{\#},y\right)$,
with $\#\in\{n,s\}$, is holomorphic and\begin{eqnarray}
\left|F_{j,\varepsilon}^{\pm}\left(x_{\#},y\right)-F_{j,\varepsilon}^{\pm}\left(x,y\right)\right| & \leq & A\left|x-x_{\#}\right|\end{eqnarray}
 for some $A>0$ independent of $\left(x,y\right)\in\ssect{j,\varepsilon}{\pm}$
and $\varepsilon$. 
\item The function $\left(x,y,\varepsilon\right)\mapsto F_{j,\varepsilon}^{\pm}\left(x,y\right)$
belongs to $\mathcal{O}_{b}\left(\mathcal{A},W\right)$. 
\item If $\gamma$ is an asymptotic path within $\ssect{j,\varepsilon}{\pm}$
with same endpoints and orientation as $\gamma_{j}^{\pm}\left(p\right)$
and lying in the same leaf, then $F_{j,\varepsilon}^{\pm}\left(p\right)=\int_{\gamma}G_{\varepsilon}\,\tau_{\varepsilon}$. 
\item Any other bounded holomorphic solution $F\in\mathcal{O}_{b}\left(\mathcal{A},W\right)$
of (\ref{eq:eq-homog}) differs from $F_{j,\varepsilon}^{\pm}$ by
the addition of a function $f\,:\,\varepsilon\mapsto f\left(\varepsilon\right)$,
$f\in\mathcal{O}_{b}\left(\mathcal{A},W\right)$, which corresponds
to the freedom in the choice of $F_{j,\varepsilon}^{\pm}\left(p_{n}\right)$. 
\end{enumerate}
\end{thm}
\begin{defn}
A function $F_{j,\varepsilon}^{\pm}$ constructed above will be called
a \textbf{sectorial solution} to the equation $X_{\varepsilon}\cdot F=G$. 
\end{defn}
The proof of this theorem is done in Section~\ref{sub:Proof-of-Theorem_secto_solve}.
The basic idea behind this result can nonetheless be shown very simply
: it is the foliated equivalent of the fundamental theorem of calculus.

\begin{lem}
\label{lem:integ-over-tangent} Let $\Omega\subset\mathcal{V}_{\varepsilon}$
be a domain and $F$ be a holomorphic function on $\Omega$. Let $\tau$
be a meromorphic 1-form on $\Omega$ such that $\tau\left(X_{\varepsilon}\right)=1$,
and let $\gamma\,:\,\left[0,1\right]\to\Omega$ be a tangent path
avoiding the poles of $\tau$. Then\begin{eqnarray}
F\left(\gamma\left(1\right)\right)-F\left(\gamma\left(0\right)\right) & = & \int_{\gamma}\left(X_{\varepsilon}\ddt F\right)\,\tau\,.\end{eqnarray}

\end{lem}
\begin{proof}
We set $G:=X_{\varepsilon}\ddt F$ which is holomorphic on $\Omega$
and use the relation :\begin{eqnarray}
\int_{\gamma}G\,\tau & = & \int_{\left[0,1\right]}\left(\gamma^{*}G\right)\left(\gamma^{*}\tau\right)\,.\end{eqnarray}
 Since $\gamma'\left(t\right)=c\left(t\right)X_{\varepsilon}\circ\gamma\left(t\right)$
we deduce $\gamma^{*}\left(\tau\right)=cdt$ and $c\gamma^{*}X_{\varepsilon}=\pp{t}$;
in particular \begin{eqnarray*}
c\gamma^{*}\left(G\right) & = & c\left(\gamma^{*}X_{\varepsilon}\right)\ddt\left(\gamma^{*}F\right)=\pp{t}\left(F\circ\gamma\right)\,\,.\end{eqnarray*}
 Hence \begin{eqnarray}
\int_{\gamma}G\tau & = & \int_{\left[0,1\right]}\pp{t}\left(F\circ\gamma\left(t\right)\right)dt\\
 & = & F\left(\gamma\left(1\right)\right)-F\left(\gamma\left(0\right)\right)\,.\nonumber \end{eqnarray}

\end{proof}

\subsection{Global equations and solutions}

We discuss, for fixed $\varepsilon$, the case where $G$ comes from
a global function $G_{\varepsilon}\in\mathcal{O}\left(\mathcal{V}_{\varepsilon}\right)$,
that is $\mathcal{A}_{\varepsilon}:=\mathcal{V}_{\varepsilon}=\cup_{j}\left(\ssect{j,\varepsilon}{+}\cup\ssect{j,\varepsilon}{-}\right)$
and $G=G_{j,\varepsilon}^{\pm}:=G_{\varepsilon}|_{\ssect{j,\varepsilon}{\pm}}$.
First let us describe how the sectorial solutions $F_{j,\varepsilon}^{\pm}$
glue :

\begin{cor}
\label{cor:solution_over_big_sectors}Let $G_{\varepsilon}\in\mathcal{O}\left(\ssect{\varepsilon}{}\right)$
such that $G_{\varepsilon}=O\left(P_{\eps}\right)+O(y)$. The sectorial
solutions $F_{j+1,\varepsilon}^{+}$ and $F_{j,\varepsilon}^{-}$
to $X_{\varepsilon}\cdot F=G_{\varepsilon}$ coincide on $\ssect{j,\varepsilon}{n}$.
The solutions $F_{j,\varepsilon}^{+}$ and $F_{\sigma\left(j\right),\varepsilon}^{-}$
coincide on $\ssect{j,\sigma(j),\varepsilon}{g}$. 
\end{cor}
\begin{proof}
Assume $p\in\ssect{j,\sigma(j),\varepsilon}{g}$. Then $\ssect{j,\varepsilon}{+}$
and $\ssect{\sigma(j),\varepsilon}{-}$ share the same node point
$p_{j,n}^{+}=p_{\sigma(j),n}^{-}$. The concatenation $\gamma_{j,\eps}^{+}(p)-\gamma_{\sigma(j),\eps}^{-}(p)$
yields an asymptotic cycle $\gamma_{\varepsilon}^{g}\left(p\right)$
through $p_{j,n}$ which is asymptotically homologous to $\left\{ p_{j,n}\right\} $
according to Proposition~\ref{pro:intersect_secto-trivial}. Item
(5) of Theorem~\ref{thm:secto-solve} yields the conclusion. The
same argument applies when $p\in\ssect{j,\varepsilon}{n}$. 
\end{proof}
Hence the obstructions to obtain global holomorphic solutions $F\in\mathcal{O}\left(\mathcal{V}_{\varepsilon}\right)$
lie only in the intersections $\ssect{j,\varepsilon}{s}$.

\begin{cor}
\label{cor:eq-homog} Let $G_{\varepsilon}\in\mathcal{O}\left(\ssect{\varepsilon}{}\right)$
such that $G_{\varepsilon}=O\left(P_{\eps}\right)+O(y)$. Let $p_{j,n}^{\pm}$
be the point of node type associated to $\sect{j,\varepsilon}{\pm}$.
We have that $p_{j,n}^{-}=p_{j+1,n}^{+}$. The asymptotic path $\gamma_{j,\varepsilon}^{s}\left(p\right)$
links the point $p_{j,n}^{-}$ to the point $p_{j,n}^{+}$. Let $I\left(j\right)$
be the value of the integral:\begin{eqnarray*}
I\left(j\right) & := & \int_{\gamma_{j,\varepsilon}^{s}\left(p\right)}G_{\varepsilon}\tau_{\varepsilon}\,,\end{eqnarray*}
 where the integration is done on canonical asymptotic paths given
by Definition~\ref{def:asy_path}.

There exists a holomorphic function $F\in\mathcal{O}\left(\ssect{\varepsilon}{}\right)$
such that $X_{\varepsilon}\cdot F=G_{\varepsilon}$ if, and only if
the following two conditions are satisfied: 
\begin{itemize}
\item the value of $I\left(j\right)$ does not depend on the choice of $p$
in a fixed sector $\mathcal{V}_{j,\eps}^{s}$. 
\item For all $m\geq1$ and all $j,j+1,\dots,j+m$ such that $p_{j,n}^{-}=p_{j+m,n}^{+}$
we have \begin{equation}
I(j)+\dots+I(j+m)=0.\label{equation_I}\end{equation}

\end{itemize}
\end{cor}
\begin{rem}
As will be seen in the last section the second condition is redundant
as the graph whose edges are the $\gamma_{j,\varepsilon}^{s}$ linking
distinct points of node type actually is a tree. 
\end{rem}
\begin{proof}
This derives from the construction of the $F_{j,\varepsilon}^{\pm}$
in (\ref{def_F}). Firstly the conditions are clearly necessary by
continuity of $F$ and Lemma~\ref{lem:integ-over-tangent}. Indeed
for all $p\in\ssect{j,\varepsilon}{s}$ there exists two node-type
singular points $p_{j,n}^{-}$ and $p_{j,n}^{+}$, not necessarily
distinct, linked by $\gamma_{j,\varepsilon}^{s}\left(p\right)$. Hence
:\begin{eqnarray}
I\left(j\right) & = & F\left(p_{j,n}^{+}\right)-F\left(p_{j,n}^{-}\right),\end{eqnarray}
 which in turn implies \eqref{equation_I}.

Let us now look at the converse: since $\mathcal{V}_{\eps}$ is connected
we can build the unique sectorial solutions given by \begin{eqnarray*}
F_{0,\varepsilon}^{-}\left(p_{0,n}^{-}\right) & := & 0\\
F_{0,\varepsilon}^{+}\left(p_{0,n}^{+}\right) & := & I\left(0\right)+F_{0,\varepsilon}^{-}\left(p_{0,n}^{-}\right)\\
F_{1,\varepsilon}^{-}\left(p_{1,n}^{-}\right) & := & F_{0,\varepsilon}^{+}\left(p_{0,n}^{+}\right)\\
\vdots & \vdots & \vdots\\
F_{j,\varepsilon}^{+}\left(p_{j,n}^{+}\right) & := & I\left(j\right)+F_{j,\varepsilon}^{-}\left(p_{j,n}^{-}\right)\\
F_{j+1,\varepsilon}^{-}\left(p_{j+1,n}^{-}\right) & := & F_{j,\varepsilon}^{+}\left(p_{j,n}^{+}\right)\\
\vdots & \vdots & \vdots\end{eqnarray*}
 The conditions precisely ensure that all $F_{j,\varepsilon}^{\pm}$
glue together in a uniform $F_{\eps}$ and that the $F_{\eps}(p_{j,n}^{\pm})$
are well defined. 
\end{proof}

\subsection{\label{sub:Proof-of-Theorem_secto_solve}Proof of Theorem~\ref{thm:secto-solve}}

\subsubsection{Preliminaries}

~

Without loss of generality we can straighten the sectorial center
manifold, since $S_{j,\varepsilon}=O\left(P_{\varepsilon}\right)$
as stated in Remark~\ref{rem:sep_estim_poly}. We drop all indices
$j$ and $\pm$ and write\begin{eqnarray}
X_{\varepsilon}\left(x,y\right) & = & X_{\varepsilon}^{M}\left(x,y\right)+y^{2}R_{\varepsilon}\left(x,y\right)\pp{y}\,.\label{separatrix_straighten}\end{eqnarray}

Let $p=\left(\ov{x},\ov{y}\right)$. We refer to Theorem~\ref{thm:secto-trivial}
and to Section~\ref{proof_4.4} for the construction of $\gamma_{\varepsilon}\left(p\right)\left(t\right)=\left(x\left(t\right),y\left(t\right)\right)$.
This path is solution to the differential system \begin{equation}
\begin{array}{lll}
\dot{x}\left(t\right) & = & \exp(i\theta)P_{\varepsilon}\left(x(t)\right)\\
\dot{y}\left(t\right) & = & \exp(i\theta)y(t)\left(1+a\left(\varepsilon\right)x\left(t\right)^{k}+y\left(t\right)R_{\varepsilon}\left(x(t),y(t)\right)\right),\end{array}\label{eq:syst_diff_path_bis}\end{equation}
 as soon as $\left|x\left(t\right)\right|<\sqrt{k}\left|\left|\varepsilon\right|\right|$.
It satisfies the same system with $\theta=0$ when $\left|x\left(t\right)\right|<\vartheta r$
for some fixed $\vartheta<1$ independent on small $\varepsilon$.

The following proposition is the key to the uniformity with respect
to $\varepsilon$.

\begin{prop}
\label{pro:length_multi_spiral} There exists a constant $C>0$ independent
of $\varepsilon\in W$ such that \begin{eqnarray*}
\int_{-\infty}^{0}\left|P_{\varepsilon}\left(x\left(t\right)\right)\right|dt & \leq & \frac{C}{\sin\delta}\left|x\left(0\right)-x_{n}\right|\end{eqnarray*}
 for any asymptotic path $t\mapsto x\left(t\right)$ landing at $x_{n}$
built in Section~\ref{proof_4.4}. The same estimate is true if we
replace $x_{n}$ by $x_{s}$ provided that we integrate $\left|P_{\varepsilon}\right|$
between $0$ and $+\infty$. 
\end{prop}
\begin{proof}
Notice that $\left|P_{\varepsilon}\left(x\left(t\right)\right)\right|=\left|\dot{x}\left(t\right)\right|$
most of the time, so what we try to achieve here is to bound the growth
of the length of spirals. We will rely on the following trivial computation
: 
\begin{lem}
\label{lem:length_spiral} Take $\overline{x}\in\ww{C}$. We consider
a logarithmic spiral $t\leq0\mapsto x\left(t\right)=\overline{x}e^{\lambda t}$
with $Re\left(\lambda\right)>0$. Then \begin{eqnarray*}
\int_{-\infty}^{0}\left|x\left(t\right)\right|dt & = & \frac{\left|\overline{x}\right|}{Re\left(\lambda\right)}\,.\end{eqnarray*}

\end{lem}
Let us first explain the strategy of the proof. It is done by induction
on the number $k+1$ of singular points enclosed in $r\ww{D}$. We
make essential use of the conic structure of $\Sigma_{0}$ by applying
the change of coordinate $x\mapsto\tilde{x}:=x\left|\left|\varepsilon\right|\right|^{-1}$
which transforms $P_{\varepsilon}\left(x\right)$ into $\left|\left|\varepsilon\right|\right|^{k+1}P_{\tilde{\varepsilon}}\left(\tilde{x}\right)$
with $\left|\left|\tilde{\varepsilon}\right|\right|=1$. Then we show
that we can isolate the singularities inside small disks $D\left(\tilde{x}_{j},\eta\right)$,
each one containing at most $k$ singular points, where $\eta$ is
independent of $\eps$ with $||\eps||=1$. Notice that if $\varepsilon$
belongs to a good sector and $m$ singular points are contained in
a disk \[
A=D\left(\tilde{x}_{j},\eta\right)\]
 then the $m$-dimensional multi-parameter associated to the $m$
singular points contained within $A$ also belongs to a good sector
of this $m$-dimensional parameter space. (The multi-parameter is
formed of the coefficients of the normalization of the monic polynomial
of degree $m$ vanishing at the singular points, the normalization
being done \emph{via} a translation so that the sum of the singular
points vanishes.) Hence we will be able to apply the recursion hypothesis
in each $A$. The last step consists in providing the estimate outside
the disks, which is not difficult.

\medskip{}

We consider \begin{equation}
\Omega_{\eta}:=\cup_{0\leq j\leq k}D\left(\tilde{x}_{j},\eta\right)\label{def_Omega}\end{equation}
 the union of disks centered at $\tilde{x}_{j}$ and of a given radius
$\eta>0$. We claim that there exists $\eta$ small enough and independent
on $\tilde{\varepsilon}$ such that $\Omega_{\eta}$ has at least
two connected components. If indeed it were not true we could find
a decreasing sequence $\left(\eta_{\ell}\right)_{\ell\in\ww{N}}\to0$
and a sequence $\left(\tilde{\varepsilon}_{\ell}\right)_{\ell}$ such
that all roots of $P_{\tilde{\varepsilon}_{n}}$ would be contained
in a domain of diameter at most $2\left(k+1\right)\eta_{\ell}$. Since
$\left\{ \tilde{\varepsilon}\,:\,\left|\left|\tilde{\varepsilon}\right|\right|=1\right\} $
is compact there exists a point of accumulation $\tilde{\varepsilon}_{\infty}$
for which the polynomial $P_{\tilde{\varepsilon}_{\infty}}$ has one
root of multiplicity $k+1$. This necessarily means $\tilde{\varepsilon}_{\infty}=0$
and is impossible. We have thus isolated in each connected component
of $\Omega_{\eta}$ at most $k$ singular points.

We assume that $\tilde{x}$ is bound to remain within some $\zeta\ww{D}$
with $\zeta>k+1$.

\medskip{}

We first deal with the case $k=1$, which contains all the ingredients
for the general case. Consider the disk $A_{n}=D\left(\tilde{x}_{0,n},\eta\right)$
containing the singular point $\tilde{x}_{n}:=\tilde{x}_{0,n}$ of
node type. This singularity is either hyperbolic or a node, in which
case the vector field $e^{i\theta}P_{\varepsilon}\pp{x}$ is linearizable
on the whole disk. Moreover $|\tilde{x}_{n}|=1$. Hence the conclusion
of Lemma~\ref{lem:length_spiral} with $\lambda:=e^{i\theta}P_{\tilde{\varepsilon}}'\left(\tilde{x}_{n}\right)$
holds : \begin{eqnarray*}
\int_{-\infty}^{0}\left|\tilde{x}\left(t\right)-\tilde{x}_{n}\right|dt & \leq & \frac{B}{Re\left(\lambda\right)}\left|\tilde{x}\left(0\right)-\tilde{x}_{n}\right|\end{eqnarray*}
 with $B$ independent on $\varepsilon$, as soon as $\tilde{x}\left(0\right)\in A_{n}$.
This type of inequality is robust under analytic changes of variables.
We have $Re\left(\lambda\right)>\left|P_{\tilde{\varepsilon}}'\left(\tilde{x}_{n}\right)\right|\sin\delta=2\sin\delta|\tilde{x}_{n}|=2\sin\delta$
according to Lemma~\ref{lem:mino_deriv}. Therefore : \begin{eqnarray*}
\int_{-\infty}^{0}\left|P_{\tilde{\varepsilon}}\left(\tilde{x}\left(t\right)\right)\right|dt & \leq & 2\zeta\int_{-\infty}^{0}\left|\tilde{x}\left(t\right)-\tilde{x}_{n}\right|dt\\
 & \leq & \frac{\zeta B}{\sin\delta}\left|\tilde{x}\left(0\right)-\tilde{x}_{n}\right|\,.\end{eqnarray*}

The same argument applies for $\tilde{x}_{s}=-\tilde{x}_{n}$ so when
$\tilde{x}\left(0\right)\in\Omega_{\eta}\backslash A_{n}$ : \begin{eqnarray*}
\int_{0}^{+\infty}\left|P_{\tilde{\varepsilon}}\left(\tilde{x}\left(t\right)\right)\right|dt & \leq & \frac{\zeta B}{\sin\delta}\left|\tilde{x}\left(0\right)-\tilde{x}_{s}\right|\,.\end{eqnarray*}
 If $\tilde{x}\left(0\right)\notin\Omega_{\eta}$ then $\left|\tilde{x}\left(0\right)-\tilde{x}_{n}\right|\geq\eta>0$
and the trajectory remains outside $\Omega_{\eta}$ only for a finite
interval of time $\left[t_{-},t_{+}\right]$. The lengths of all those
trajectories are uniformly bounded by some $M=M\left(\zeta\right)>0$
for all $\tilde{\varepsilon}\in W$ with $\left|\left|\tilde{\varepsilon}\right|\right|=1$.
We can thus write\begin{eqnarray*}
\int_{t_{-}}^{t_{+}}\left|P_{\tilde{\varepsilon}}\left(\tilde{x}\left(t\right)\right)\right|dt & \leq & \frac{M}{\eta}\left|\tilde{x}\left(0\right)-\tilde{x}_{n}\right|\end{eqnarray*}
 and the same for $\tilde{x}_{s}$. When $\zeta$ is large enough
the trajectories look like the trajectories of $\dot{x}=x^{2}$ (which
are circles tangent at the origin to the real axis). Therefore $M\left(\zeta\right)\leq M'\zeta$
with $M'$ independent on large $\zeta$ and on $\tilde{\varepsilon}$.
To conclude, if we let $B_{1}=B+\frac{M\sin\delta}{\eta\zeta}$, we
obtain \begin{eqnarray*}
\int_{-\infty}^{0}\left|P_{\tilde{\varepsilon}}\left(\tilde{x}\left(t\right)\right)\right|dt & \leq & \frac{\zeta B_{1}}{\sin\delta}\left|\tilde{x}\left(0\right)-\tilde{x}_{n}\right|\end{eqnarray*}
 for all $\tilde{x}\left(0\right)\in\zeta\ww{D}$ and all $\zeta>0$,
for all $\tilde{\varepsilon}\in W$ with $\left|\left|\tilde{\varepsilon}\right|\right|=1$.

Back to the original coordinates we find : \begin{eqnarray*}
\int_{-\infty}^{0}\left|P_{\varepsilon}\left(x\left(t\right)\right)\right|dt & =\left|\left|\varepsilon\right|\right|^{2}\int_{-\infty}^{0}\left|P_{\tilde{\varepsilon}}\left(\tilde{x}\left(t\right)\right)\right|dt & \leq\frac{\left|\left|\varepsilon\right|\right|\zeta B_{1}}{\sin\delta}\left|x\left(0\right)-x_{n}\right|\end{eqnarray*}
 for all $\left|x\left(0\right)\right|<\zeta\left|\left|\varepsilon\right|\right|$.
By letting $\zeta:=r\left|\left|\varepsilon\right|\right|^{-1}$ we
obtain $C:=rB_{1}$.

\medskip{}

We deal now briefly with the general case $k>1$ in a similar way.
Inside each component $A$ of $\Omega_{\eta}$ we apply the recursion
hypothesis since there are at most $k$ singular points lying within
$A$. The argument developed just above applies again to obtain the
bound in $\zeta\ww{D}\backslash\Omega_{\eta}$, then $\zeta\ww{D}$.
We finally derive :\begin{eqnarray*}
\int_{-\infty}^{0}\left|P_{\varepsilon}\left(x\left(t\right)\right)\right|dt & \leq & \frac{\left|\left|\varepsilon\right|\right|^{k}\zeta^{k}B_{k}}{\sin\delta}\left|x\left(0\right)-x_{n}\right|\end{eqnarray*}
 for all $x\left(0\right)\in\zeta\left|\left|\varepsilon\right|\right|\ww{D}$.
The conclusion is reached by setting $\zeta:=r\left|\left|\varepsilon\right|\right|^{-1}$. 
\end{proof}

\subsubsection{Back to the proof of Theorem~\ref{thm:secto-solve}}

~

\textbf{(1)} We write \textbf{$\gamma_{\varepsilon}\left(\overline{x},\overline{y}\right)\left(t\right)=\left(x\left(t\right),y\left(t\right)\right)$}
for \textbf{$t\leq0$}. Lemma~\ref{lem:modulus_estimate} yields
\begin{eqnarray}
\left|y(t)\right| & \leq\left|\ov{y}\right|e^{\alpha t}\label{eq:ineg_y}\end{eqnarray}
 for some $\alpha>\frac{1}{2}\cos\theta$ independent on small $\varepsilon\in W$.

Assume that for $\varepsilon\in W$ \begin{equation}
G_{\varepsilon}\left(x,y\right)=P_{\varepsilon}\left(x\right)g_{0,\varepsilon}\left(x\right)+yg_{1,\varepsilon}\left(x,y\right).\label{eq_g_1}\end{equation}
 There exists constants $M_{1},\, M_{2}>0$, both independent on $\overline{x}$
and $\varepsilon$, such that \begin{eqnarray}
\left|G_{\varepsilon}\left(x\left(t\right),y\left(t\right)\right)\right| & \leq & M_{1}\left|P_{\varepsilon}\left(x\left(t\right)\right)\right|+M_{2}\bar{y}e^{\alpha t}\,\,.\label{eq:estim-integ}\end{eqnarray}
 The integral $\int_{-\infty}^{0}G\left(x\left(u\right),y\left(u\right)\right)du$
is thus absolutely convergent according to Proposition~\ref{pro:length_multi_spiral}.
Since $\tau_{\varepsilon}\left(X_{\varepsilon}\right)=1$ it follows
that \begin{eqnarray}
\int_{t}^{0}\gamma_{\varepsilon}\left(p\right)^{*}\left(G_{\varepsilon}\tau_{\varepsilon}\right) & = & e^{i\theta}\int_{t}^{0}G\left(x\left(u\right),y\left(u\right)\right)du\label{eq:integ-rep}\end{eqnarray}
 for $t\leq0$, which completes the proof.

\medskip{}

\textbf{(2)} From Proposition~\ref{pro:leaf-is-graph} the leaf passing
through $\left(x,y\right)$ is the graph of a holomorphic function
$v\mapsto l\left(v,y\right)$. Because of (\ref{eq:integ-rep}) the
function $p\mapsto F_{\varepsilon}\left(p\right)$ is holomorphic.
Besides \begin{eqnarray}
X_{\varepsilon}\ddt F_{\varepsilon}\left(p\right) & = & \left.\left(e^{-i\theta}\frac{d}{dt}e^{i\theta}\int_{-\infty}^{t}G\left(x(u),l\left(x\left(u\right),\overline{y}\right)\right)du\right)\right|_{t=0}\\
 & = & G\left(p\right)\nonumber \end{eqnarray}
 as in Lemma~\ref{lem:integ-over-tangent}.

\textbf{(3)} Let us consider a point $p:=\left(x_{n},\ov{y}\right)$
and show that\begin{eqnarray}
\lim_{\left(x,y\right)\to p}F_{\varepsilon}\left(x,y\right) & = & \int_{\left[0,\overline{y}\right]}g_{1,\varepsilon}\left(x_{n},y\right)\frac{dy}{1+a\left(\varepsilon\right)x_{n}^{k}+yR_{\varepsilon}\left(x_{n},y\right)}\,,\label{eq:solution_on_separatrix_1}\end{eqnarray}
 where $g_{1,\eps}$ is defined in \eqref{eq_g_1}. The latter integral
is holomorphic with respect to small $\ov{y}$ and can be rewritten
(using \eqref{eq:syst_diff_path_bis}) \begin{eqnarray}
F_{\varepsilon}\left(p\right) & := & e^{i\theta}\int_{-\infty}^{0}\frac{g_{1,\varepsilon}\left(x_{n},\tilde{y}(t)\right)}{1+a\left(\varepsilon\right)x_{n}^{k}+\tilde{y}\left(t\right)R_{\varepsilon}\left(x_{n},\tilde{y}\left(t\right)\right)}\tilde{y}\left(t\right)dt\end{eqnarray}
 where $\left(\tilde{x}\left(t\right),\tilde{y}\left(t\right)\right)$
is the solution for $t\leq0$ to the system\begin{equation}
\begin{array}{lll}
\dot{\tilde{x}} & = & \frac{e^{i\theta}P\left(\tilde{x}\right)}{1+a\left(\varepsilon\right)\tilde{x}^{k}+\tilde{y}R_{\varepsilon}\left(\tilde{x},\tilde{y}\right)}\\
\dot{\tilde{y}} & = & e^{i\theta}\tilde{y}\end{array}\label{system_tilde}\end{equation}
 with initial condition $\left(\overline{x},\overline{y}\right)$.

We only need to prove that $\lim_{x\to x_{n}}F_{\varepsilon}\left(x,\ov{y}\right)=F_{\varepsilon}\left(x_{n},\ov{y}\right)$
because $F_{\varepsilon}\left(x_{n},\cdot\right)$ is holomorphic
and the integral is additive with respect to the path of integration.
There exists two constants $L_{1},L_{2}>0$ independent on $\varepsilon$
such that : \begin{eqnarray*}
\left|G\left(x,y\right)-G\left(x_{n},y\right)\right| & \leq & L_{1}\left|P_{\varepsilon}\left(x\right)\right|+L_{2}\left|y\right|\,.\end{eqnarray*}
 Since as $t\mapsto-\infty$ the modulus $\left|\tilde{y}\left(t\right)\right|$
is exponentially flat of rate $\alpha$ close to $1$ while $\left|\tilde{x}\left(t\right)-x_{n}\right|$
is exponentially flat at an order tending toward $0$ as $\varepsilon$
gets smaller, we can write $\left|\tilde{y}\left(t\right)\right|\leq L_{3}\left|P_{\varepsilon}\left(\tilde{x}\left(t\right)\right)\right|$.
Hence there exists some $L$ depending only on $r,\, r',\,\rho$ such
that, for all $t\leq0$, \begin{eqnarray*}
\left|G\left(\tilde{x}\left(t\right),\tilde{y}\left(t\right)\right)-G\left(x_{n},\tilde{y}\left(t\right)\right)\right| & \leq & L\left|P_{\varepsilon}\left(\tilde{x}\left(t\right)\right)\right|\,,\end{eqnarray*}
 from which we deduce\begin{eqnarray}
\left|\int_{-\infty}^{0}G\left(\tilde{x}(t),\tilde{y}\left(t\right)\right)-G\left(x_{n},\tilde{y}(t)\right)dt\right| & \leq L & \int_{-\infty}^{0}\left|P_{\varepsilon}\left(\tilde{x}\left(t\right)\right)\right|dt\label{eq:estim_integ_sep}\\
 & \leq & \frac{CL}{\sin\delta}\left|\ov{x}-x_{n}\right|\,,\nonumber \end{eqnarray}
 by Lemma~\ref{lem:modulus_estimate} and Proposition~\ref{pro:length_multi_spiral}.
Hence $F_{\varepsilon}$ extends continuously to $\left\{ x_{n}\right\} \times r'\ww{D}$
and we have moreover proved that\begin{eqnarray}
\left|F_{\varepsilon}\left(\ov{x},\ov{y}\right)-F_{\varepsilon}\left(x_{n},\ov{y}\right)\right| & \leq & A\left|\ov{x}-x_{n}\right|\,.\end{eqnarray}

\begin{figure}
\includegraphics[width=8cm]{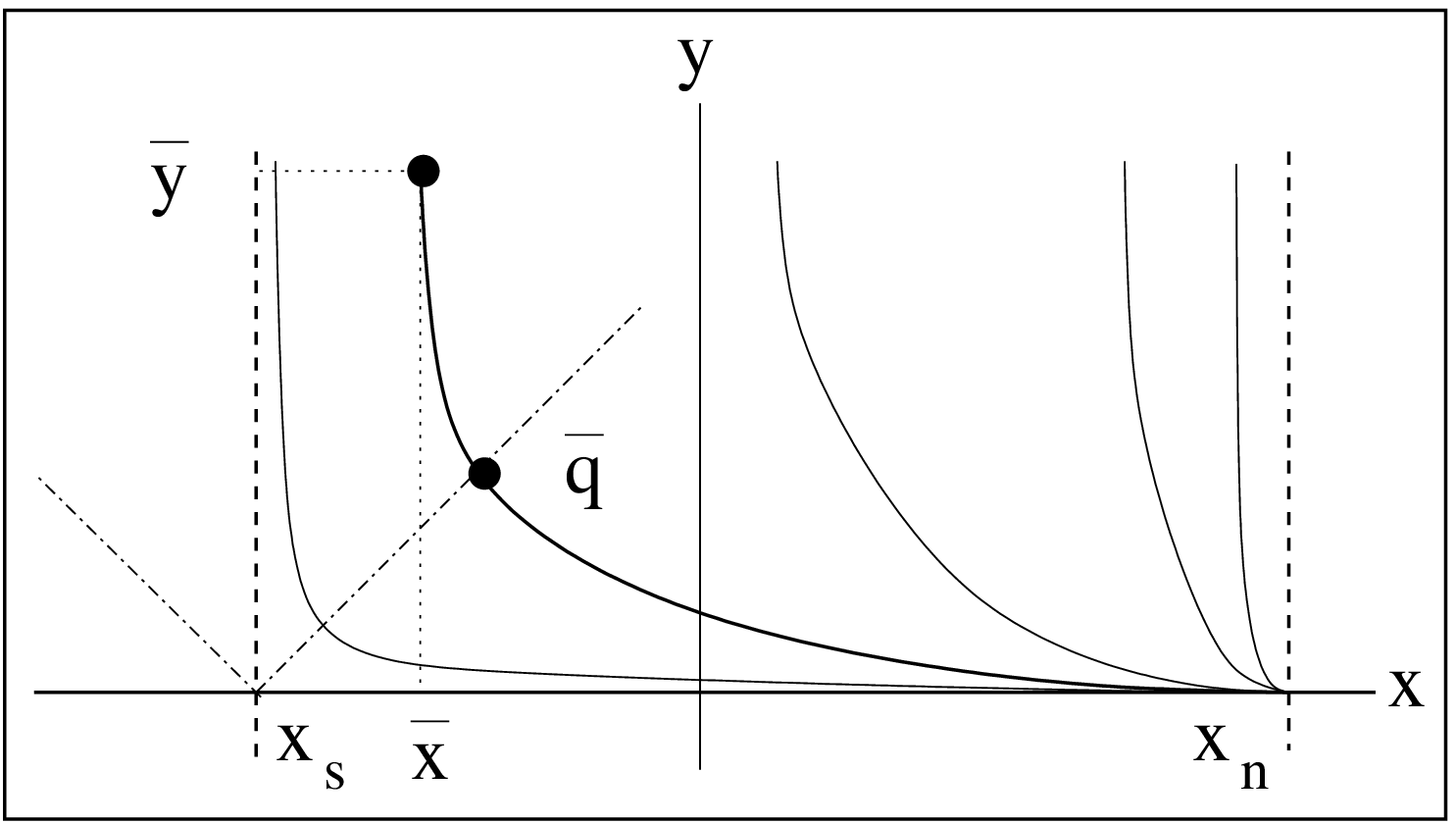}

\caption{\label{fig:integ_estim}} 
\end{figure}

We must now deal with the separatrix $\left\{ x_{s}\right\} \times r'\ww{D}$
using much the same argument. Consider an asymptotic path $\gamma$
linking the two singularities within $\sect{\varepsilon}{}\times\left\{ 0\right\} $
that is, more precisely, $\gamma\left(\infty\right)=p_{s}$ and $\gamma\left(-\infty\right)=p_{n}$.
The estimate (\ref{eq:estim-integ}) shows that the integral $\int_{\gamma}G\tau_{\varepsilon}$
converges, therefore we claim that $F_{\varepsilon}$ extends continuously
on $\left\{ x_{s}\right\} \times r'\ww{D}$ to \begin{eqnarray}
F_{\varepsilon}\left(x_{s},\ov{y}\right) & := & \int_{\gamma}G\tau_{\varepsilon}+\int_{\left[0,\ov{y}\right]}g_{1,\eps}\left(x_{s},y\right)\frac{dy}{1+a\left(\varepsilon\right)x_{s}^{k}+yR_{\varepsilon}\left(x_{s},y\right)}\,.\label{eq:solution_on_separatrix_2}\end{eqnarray}
 \medskip{}
 If $\ov{y}=0$ the result is trivial. For a given $\ov{y}\neq0$
we can choose $\ov{x}-x_{s}$ sufficiently small so that the path
$\gamma_{\varepsilon}\left(\overline{x},\overline{y}\right)$ crosses
the real analytic set $\left\{ \left|x-x_{s}\right|=\left|y\right|\right\} $
at some point $\overline{q}=\overline{q}\left(\ov{x}\right)$ (see
Figure~\ref{fig:integ_estim}). Because of (\ref{eq:ineg_y}) we
know that $\overline{q}$ is unique. Let $\overline{s}=\overline{s}\left(\overline{x}\right)\leq0$
be such that $\overline{q}=\left(\tilde{x}\left(\overline{s}\right),\tilde{y}\left(s\right)\right)$.
Obviously \begin{eqnarray}
\lim_{\ov{x}\to x_{s}}\overline{s} & = & -\infty\,,\\
\lim_{\ov{x}\to x_{s}}\tilde{x}\left(\overline{s}\right) & = & x_{s}\,.\end{eqnarray}
 Define $\overline{t}\leq0$ such that $x\left(\overline{t}\right)=\tilde{x}\left(\overline{s}\right)$
($\ov{t}$ (\emph{resp}. $\ov{s}$) is the time for $(\ov{x},\ov{y})$
to reach $\ov{q}$ in \eqref{eq:syst_diff_path_bis} (\emph{resp.}
\eqref{system_tilde})). On the one hand\begin{eqnarray}
\int_{-\infty}^{\overline{t}}\left|G\left(x\left(t\right),y\left(t\right)\right)-G\left(x\left(t\right),0\right)\right|dt & \leq & L'\left|y\left(\overline{t}\right)\right|\alpha^{-1}e^{\alpha\overline{t}}\\
 & \leq & \frac{2L'}{\cos\theta}\left|\tilde{x}\left(\overline{s}\right)-x_{s}\right|\nonumber \end{eqnarray}
 whereas on the other hand \begin{eqnarray*}
\int_{\overline{s}}^{0}\left|G\left(\tilde{x}\left(t\right),\tilde{y}\left(t\right)\right)-G\left(x_{s},\tilde{y}\left(t\right)\right)\right|dt & \leq & L\int_{0}^{-\overline{s}}\left|P_{\varepsilon}\left(\tilde{x}\left(t+\overline{s}\right)\right)\right|dt\\
\
 & \leq & L\int_{0}^{+\infty}\left|P_{\varepsilon}\left(\tilde{x}\left(t+\overline{s}\right)\right)\right|dt\\
 & \leq & \frac{CL}{\sin\delta}\left|\tilde{x}\left(\overline{s}\right)-x_{s}\right|\,.\end{eqnarray*}
 From the fact that $t\mapsto\left|y\left(t\right)\right|$ is exponentially
increasing we infer $\tilde{x}\left(\overline{s}\right)-x_{s}=O\left(\overline{x}-x_{s}\right)$,
which ends the proof.

\medskip{}

\textbf{(4)} From (\ref{eq:solution_on_separatrix_1}) and (\ref{eq:solution_on_separatrix_2})
we deduce that $y\mapsto F_{j,\varepsilon}^{\pm}\left(x_{\#},y\right)$
belongs to $\mathcal{O}_{b}\left(\left\{ x_{\#}\right\} \times r'\ww{D},W\right)$.
Because of the estimate obtained in (3) we deduce that $F_{j,\varepsilon}^{\pm}$
belongs to $\mathcal{O}_{b}\left(\ssect{j,\varepsilon}{\pm},W\right)$.

\medskip{}

\textbf{(5)} Let $h$ be an asymptotic homology between $\gamma=h\left(-\infty,\cdot\right)$
and $\gamma_{\varepsilon}\left(p\right)=h\left(+\infty,\cdot\right)$
within $\ssect{\varepsilon}{}$. Because we integrate the restriction
of $G$ to a leaf (a holomorphic surface) we only need to prove that
\begin{eqnarray*}
\lim_{t\to\pm\infty}\int_{h\left(\cdot,t\right)}G\tau_{\varepsilon} & = & 0\,.\end{eqnarray*}
 But this is obvious since the length of $s\mapsto h\left(s,t\right)$
is bounded (thanks to the uniformity of the convergence $h\left(s,\cdot\right)\to h\left(\pm\infty,\cdot\right)$)
and $G\left(h\left(s,t\right)\right)\to0$ as $t\to\pm\infty$ uniformly
in $s$.

\medskip{}

\textbf{(6)} Let $F$ be another solution. Then \begin{eqnarray}
X_{\varepsilon}\ddt\left(F-F_{\varepsilon}\right) & = & 0\,.\end{eqnarray}
 Hence there must exist a bounded, holomorphic function $\left(\varepsilon,h\right)\mapsto f\left(\varepsilon,h\right)$
such that, for all $p\in\mathcal{V}_{\varepsilon}$: \begin{eqnarray}
F\left(p\right)-F_{\varepsilon}\left(p\right) & = & f\left(\varepsilon,H_{\varepsilon}\left(p\right)\right)\end{eqnarray}
 where $H_{\varepsilon}$ is the canonical first integral (see Corollary~\ref{cor:first-integ-secto}).
According to the same corollary $H_{\varepsilon}\left(\mathcal{V}_{\varepsilon}\right)=\ww{C}$,
so that $h\mapsto f\left(\varepsilon,h\right)$ must be an entire
function, thus constant.

\section{\label{orb_classification}The sectorial normalization theorem}

As this is the principal tool to the identification of moduli of analytic
classification for generic families unfolding a saddle-node of codimension
$k$ we start by some generalities on this matter.

\subsection{The principle of a modulus of analytic classification}

The heart of the paper is to identify complete moduli of analytic
classification for generic families unfolding a saddle-node of codimension
$k$, under either

\begin{itemize}
\item orbital equivalence, 
\item or conjugacy. 
\end{itemize}
The idea is the following: we start with two germs of prepared families
$Z_{\eps}=U_{\eps}X_{\eps}$ and $Z'_{\eps}=U'_{\eps}X'_{\eps}$ and
we want to decide whether they are orbitally equivalent or conjugate.
We know that the multi-parameter is canonical and that an equivalence
or conjugacy must preserve the parameter (up to the equivalence relation
\eqref{cong1}). Applying a rotation of order $k$ to one of the families,
we can suppose that the two families have the same canonical multi-parameter
and the same polynomial $P_{\eps}(x)$. We can then limit ourselves
to discuss orbital equivalence or conjugacy preserving the multi-parameter.

For that purpose we construct orbital equivalences or conjugacies
with the model family $Z_{\eps}^{M}=Q_{\eps}X_{\eps}^{M}$ over canonical
sectors. The modulus is a measure of the obstruction to glue these
in a global equivalence or conjugacy with the model. By composing
them, this provides equivalences or conjugacies between the two families
over canonical sectors. These provide a global equivalence or conjugacy
between the two families precisely when the two families have equal
moduli.

\medskip{}

Solving the orbital equivalence problem of $Z_{\eps}$ with $Z_{\eps}^{M}$
on a canonical sector is the same as solving the conjugacy problem
between $X_{\eps}$ and $X_{\eps}^{M}$ under a change of coordinates
preserving the $x$-variable. Once the sectorial center manifold has
been straightened, this change of coordinates, $\Psi_{\eps}^{N}$,
will be taken as the flow $\Phi_{Y}^{N_{\eps}}$ of the vector field
\begin{equation}
Y:=y\pp{y}\label{field_Y}\end{equation}
 for some time $N_{\eps}(x,y)$ with $N_{\eps}$ analytic and we will
have $\left(\Phi_{Y}^{N_{\eps}}\right)^{*}\left(X_{\varepsilon}^{M}\right)=X_{\varepsilon}$.

To solve the conjugacy problem we first solve the equivalence problem
and then consider the conjugacy problem between $Q_{\varepsilon}X_{\varepsilon}$
and $Z_{\eps}^{M}$. This will be done through a map $\Psi_{\eps}^{T}=\Phi_{Q_{\varepsilon}X_{\varepsilon}}^{T_{\eps}}$,
\emph{i.e.} the flow of $Q_{\eps}X_{\eps}$ for some time $T_{\eps}\left(x,y\right)$
given by an analytic function $T_{\eps}$.

\medskip{}

Both maps $N_{\varepsilon}$ (for Normal change of coordinates) and
$T_{\varepsilon}$ (for Tangential change of coordinates) will be
obtained by solving cohomological equations. The justification of
this procedure comes from the following {}``fundamental'' lemma,
which can be proved using Lie's formal formula for the flow of a vector
field : \begin{eqnarray}
\flw{X}{t}{} & = & \sum_{n\geq0}\frac{t^{n}}{n!}X\cdot^{n}Id\end{eqnarray}

\begin{lem}
\label{lem:funda} \cite{T4} Let $X$ and $Y$ be germs at $\left(0,0\right)$
of commuting, holomorphic vector fields. Consider a function $F\in\mathcal{O}\left(\mathcal{W}\right)$,
where $\mathcal{W}$ is a domain on which $X$ and $Y$ are holomorphic.
Assume that $Y\cdot F\left(\overline{p}\right)\neq-1$ for some $\overline{p}\in\mathcal{W}$.
Then the map $\psi$ defined by \begin{eqnarray}
\psi\left(x,y\right) & := & \flw{Y}{F(x,y)}{(x,y)}\end{eqnarray}
 is a local change of coordinates near $\overline{p}$ satisfying
:\begin{eqnarray}
\psi^{*}\left(X\right) & = & X-\frac{X\cdot F}{1+Y\cdot F}Y\,.\end{eqnarray}

\end{lem}
By setting $X=Y$ we deduce a useful sufficient condition (which is
generically necessary) for two vector fields inducing the same foliation
to be conjugate.

\begin{cor}
\label{cor:bridge} Let $X$ be a germ of a vector field with a singularity
at $\left(0,0\right)$ and $U,\, V$ two non-vanishing holomorphic
germs. If there exists a germ of holomorphic function $T$ such that\begin{eqnarray}
X\cdot T & = & \frac{1}{V}-\frac{1}{U}\label{hom_T}\end{eqnarray}
 then the vectors fields $UX$ and $VX$ are (locally) analytically
conjugate, the conjugacy being given by $\Phi_{UX}^{T}$, namely $\left(\Phi_{UX}^{T}\right)^{*}\left(UX\right)=VX$.
In particular $U\left(0,0\right)=V\left(0,0\right)$. 
\end{cor}
\begin{proof}
It follows from Lemma~\ref{lem:funda} applied with $X,Y\mapsto UX$.
The hypotheses of the lemma are satisfied for $\overline{p}:=\left(0,0\right)$
(indeed $X\left(0,0\right)=0$). Moreover \eqref{hom_T} is solvable
only if $\frac{1}{V\left(0,0\right)}-\frac{1}{U\left(0,0\right)}=0$. 
\end{proof}

\subsection{Sectorial normalization}

Let $Z_{\eps}$ be given in \eqref{unfold3} and its restriction on
\begin{equation}
\ssect{j,\varepsilon}{}=\ssect{j,\varepsilon}{+}\cup\ssect{j,\varepsilon}{-}.\end{equation}
 On $\ssect{j,\varepsilon}{}$ we consider the change of coordinates
$(x,y)\mapsto\Phi_{j,\eps}(x,y):=\left(x,y-S_{j,\eps}(x)\right)$
of Corollary~\ref{cor:red_center_man} and let \begin{equation}
X_{j,\eps}:=\left(\Phi_{j,\eps}\right)_{*}X_{\eps}\,\,.\label{vect_field_j}\end{equation}
 We call $\tau_{\eps}=\frac{dx}{P_{\varepsilon}}$ the corresponding
time-form. Then $X_{j,\eps}$ has the form \begin{equation}
X_{j,\eps}=P_{\eps}(x)\frac{\partial}{\partial x}+\left(y\left(1+a(\eps)x^{k}-\tilde{R}_{j,\eps}(x,y)\right)\right)\frac{\partial}{\partial y}\,\,,\label{eq-norm}\end{equation}
 with $\tilde{R}_{j,\eps}(x,y)=O(y)$.

The following proposition, and its extension in Theorem~\ref{thm:secto-norma},
gives a geometric proof of the theorem of Hukuhara-Kimura-Matuda \cite{HKM}
($\varepsilon=0$) together with its generalization to unfoldings.
For $\varepsilon=0$ one can recover the summability property using
the theorem of Ramis-Sibuya and the fact that the first-integrals
in the intersections $\ssect{j,0}{s}$ are flat of order $\frac{1}{k}$.

\begin{prop}
\label{6.cor} Let $\Phi_{j,\varepsilon}$ be the change of coordinates
$\left(x,y\right)\mapsto\left(x,y-S_{j,\varepsilon}\left(x\right)\right)$
which straightens the sectorial separatrix (see Corollary~\ref{cor:red_center_man}). 
\begin{enumerate}
\item We consider $X_{j,\eps}=X_{\eps}^{M}-\tilde{R}_{j,\eps}Y$ over $\ssect{j,\varepsilon}{}$,
given in \eqref{eq-norm}, where $Y$ is given in \eqref{field_Y}.
Let $N_{j,\eps}^{\pm}$ be the solution of $X_{j,\eps}\cdot N_{j,\eps}^{\pm}=\tilde{R}_{j,\eps}$
on $\ssect{j,\varepsilon}{\pm}$ with $N_{j,\varepsilon}^{\pm}\left(p_{j,n}\right)=0$,
and $\Psi^{N}:=\Phi_{Y}^{N_{j,\eps}^{\pm}}\circ\Phi_{j,\varepsilon}$.
Then $\left(\Psi^{N}\right)^{*}X_{\eps}^{M}=X_{\eps}$. 
\item Let $T_{j,\eps}$ be the sectorial solution with $T_{j,\varepsilon}\left(p_{j,n}\right)=0$
of \begin{equation}
X_{\varepsilon}\cdot T_{j,\eps}=\frac{1}{U_{\eps}}-\frac{1}{Q_{\varepsilon}}\end{equation}
 on $\ssect{j+1,\varepsilon}{+}\cup\ssect{j,\varepsilon}{-}$ and
let $\Psi^{T}=\Phi_{Q_{\eps}X_{\eps}}^{T_{j,\eps}}$. Then $\left(\Psi^{T}\right)^{*}\left(Q_{\eps}X_{\varepsilon}\right)=Z_{\eps}.$ 
\item Moreover the functions $N_{j,\eps}^{\pm}$ and $T_{j,\eps}$ are bounded
on $\ssect{j,\varepsilon}{\pm}$ and have continuous extensions to
$\ssect{j,\varepsilon}{\pm}\cup\left(\left\{ x_{j,s},x_{j,n}\right\} \times r'\ww{D}\right)$,
and \begin{equation}
N_{j,\eps}^{\pm}\left(p_{j,\#}\right)=T_{j,\eps}\left(p_{j,\#}\right)=0,\end{equation}
 where $\#\in\left\{ s,n\right\} $. 
\item The functions $N_{j,\eps}^{\pm}$ and $T_{j,\eps}^{\pm}$ are uniformly
bounded when $\varepsilon$ belongs to some good sector in $\Sigma_{0}$. 
\end{enumerate}
\end{prop}
\begin{proof}
~ 
\begin{enumerate}
\item We apply Lemma~\ref{lem:funda} to \eqref{vect_field_j} with $X=X_{\eps}^{M}$
and $Y=y\frac{\partial}{\partial y}$. Then $\left(\Phi_{Y}^{N_{j,\eps}^{\pm}}\right)^{*}X_{\eps}^{M}=X_{\eps}^{M}-\tilde{R}_{j,\eps}Y$
since $\left(X_{\eps}^{M}-\tilde{R}_{j,\eps}Y\right)\cdot N_{j,\eps}^{\pm}=\tilde{R}_{j,\eps}$.
The conclusion follows using Corollary~\ref{cor:red_center_man}. 
\item This is immediate from Corollary~\ref{cor:bridge}. 
\item and (4) are consequences of Theorem~\ref{thm:secto-solve}. 
\end{enumerate}
\end{proof}
We summarize these results in the sectorial normalization Theorem~\ref{thm:secto-norma}
below. To state it we require the notion of sectorial diffeomorphism.

\subsubsection{Sectorial diffeomorphisms}

\begin{defn}
\label{def:secto-diffeo} Recall that the squid sectors $\sect{j,\varepsilon}{\pm}$
actually depend on an angle $\theta$ and a width $w_{0}>0$ (see
Definition~\ref{def:squid-sectors} and Figure~\ref{fig:z_squid}). 
\begin{enumerate}
\item Let $\varepsilon\in\Sigma_{0}\cup\left\{ 0\right\} $ and $\#\in\left\{ -,+,s,n,g\right\} $
be given. Define $\sect{j,\varepsilon}{\#}\left(\eta\right):=\sect{j,\varepsilon}{\#}\cap\eta\mathbb{D}$.
A \textbf{germ of a sectorial diffeomorphism} over $\ssect{j,\varepsilon}{\#}$
is a (class of) map(s) $\Psi_{j,\eps}^{\#}:\,\tilde{\mathcal{V}}_{j,\eps}^{\#}\to\ww{C}^{2}$
holomorphic and one-to-one on \begin{equation}
\tilde{\mathcal{V}}_{j,\eps}^{\#}:=\sect{j,\varepsilon}{\#}(r_{0})\times r_{0}'\ww{D}\label{r_0}\end{equation}
 for $r_{0},r_{0}'>0$ sufficiently small, satisfying:

\begin{enumerate}
\item $\Psi_{j,\eps}^{\#}$ extends to a homeomorphism $\Psi_{j,\eps}^{\#}$
defined on $\mathcal{W}:=\left(\sect{j,\varepsilon}{\#}(r_{0})\cup\left\{ p_{j,s},p_{j,n}\right\} \right)\times r_{0}'\ww{D}$,
fixing the singularities and such that $\Psi_{j,\eps}^{\#}\left(\left\{ p_{j,*}\right\} \times r_{0}'\ww{D}\right)\subset\left\{ p_{j,*}\right\} \times\ww{C}$
for each $*\in\{s,n\}$. 
\item The image $\Psi_{j,\eps}^{\#}\left(\tilde{\mathcal{V}}_{j,\eps}^{\#}\right)$
of the squid sector $\tilde{\mathcal{V}}_{j,\eps}^{\#}$ is squeezed
between two squid sectors: \[
\mathcal{V}_{1}\subset\Psi_{j,\varepsilon}^{\#}\left(\tilde{\mathcal{V}}_{j,\eps}^{\#}\right)\subset\mathcal{V}_{2},\]
 where \begin{equation}
\mathcal{V}_{\ell}:=\sect{j,\varepsilon}{\#}(r_{\ell})\times r'_{\ell}\ww{D},\qquad\ell=1,2,\label{bounding_sectors}\end{equation}
 have same angle $\theta$ but with maybe different widths $w_{\ell}$
and $r'_{\ell}$. 
\end{enumerate}
\item Let $W\subset\Sigma_{0}$ be some good sector. A \textbf{germ of a
family of sectorial diffeomorphisms} is a family of canonical sectors
$\ssect{j,\varepsilon}{\#}$ together with a family $\left(\Psi_{j,\varepsilon}^{\#}\right)_{\varepsilon\in W\cup\left\{ 0\right\} }$
of germs of sectorial diffeomorphisms for all values of $\eps\in W\cup\left\{ 0\right\} $
with $||\eps||\leq\rho$ for some $\rho>0$ and for which we can choose
$w_{\ell},r_{\ell},r_{\ell}'$, $\ell\in\{0,1,2\}$, of \eqref{r_0}
and \eqref{bounding_sectors} independent on $\eps$. 
\end{enumerate}
\end{defn}
This implies that sectorial diffeomorphisms respect locally the fibered
squid sectors, \emph{e}.\emph{g.} neither crush them nor blow them
away along the separatrices $\left\{ p_{j,*}\right\} \times r'\ww{D}$.
This property is necessary to ensure that we are able to construct
holomorphic conjugacies on a full neighborhood of the singularities
when the moduli of two vector fields coincide. In practice we will
consider a good covering of $\Sigma_{0}$ given by Theorem~\ref{covering}
and we will construct germs of families of sectorial diffeomorphisms
depending analytically on $\eps$ on each open set of the covering.

\subsubsection{Sectorial normalization theorem}

\begin{thm}
\label{thm:secto-norma}Let $W\subset\Sigma_{0}$ be a good sector.
There exists $r,r',\rho>0$ sufficiently small so that, for any $\eps\in W\cup\left\{ 0\right\} $
with $||\eps||\leq\rho$ and associated set of canonical sectors,
the vector field $Z_{\varepsilon}$ is conjugate to its model $Z_{\eps}^{M}$
by a sectorial diffeomorphism $\Psi_{j,\varepsilon}^{\pm}$ over $\mathcal{V}_{j,\varepsilon}^{\pm}$.
The change of coordinates splits into an orbital part, namely \begin{eqnarray}
\Psi_{j,\varepsilon}^{N,\pm}\left(x,y\right) & := & \left(x,\left(y-S_{j,\varepsilon}^{\pm}\left(x\right)\right)\exp\left(N_{j,\varepsilon}^{\pm}\left(x,y\right)\right)\right)\end{eqnarray}
 transforming $Q_{\varepsilon}X_{\varepsilon}^{N}$ into $Q_{\eps}X_{\eps}$
composed with a tangential part, namely \begin{eqnarray}
\Psi_{j,\varepsilon}^{T}\left(x,y\right) & := & \flw{Q_{\varepsilon}X_{\varepsilon}}{T_{j,\varepsilon}\left(x,y\right)}{\left(x,y\right)}.\end{eqnarray}
 transforming $Q_{\eps}X_{\eps}$ into $U_{\eps}X_{\eps}$. Both $\left(\Psi_{j,\varepsilon}^{N,\pm}\right)_{\varepsilon\in W\cup\left\{ 0\right\} }$
and $\left(\Psi_{j,\varepsilon}^{T}\right)_{\varepsilon\in W\cup\left\{ 0\right\} }$
are families of sectorial diffeomorphisms, over $\left(\ssect{j,\varepsilon}{\pm}\right)_{\varepsilon}$
and $\left(\ssect{j+1,\varepsilon}{+}\cup\ssect{j,\varepsilon}{-}\right)_{\varepsilon}$
respectively. 
\end{thm}
\begin{proof}
The holomorphy of the changes of coordinates, their dependence on
$\varepsilon$ and the continuity property for families of sectorial
diffeomorphisms follow from Theorem~\ref{thm:secto-solve} and Proposition~\ref{6.cor}.

For the sake of clarity we omit to write the upper and lower indices
$\varepsilon$, $j$ and $\pm$. We prove that $\Psi^{N}$ and $\Psi^{T}$
are sectorial diffeomorphisms. The easiest case is $\Psi^{N}$ since
it preserves the $x$-coordinate. According to Theorem~\ref{thm:secto-solve}(3)
we have \begin{eqnarray}
\left|N\left(x,y\right)-N_{\#}\left(y\right)\right| & \leq & A\left|x-x_{\#}\right|\end{eqnarray}
 with $A$ independent of $\varepsilon$ and $N_{\#}\left(y\right):=N\left(x_{\#},y\right)$
for $\#\in\left\{ s,n\right\} $. Hence, if we let $\Psi^{N}=\left(Id,\psi_{1}\right)$
then \begin{eqnarray}
\left|\psi_{1}\left(x,y\right)-\psi_{1}\left(x_{\#},y\right)\right| & \leq & \left|y-S_{j}\left(x\right)\right|\left|e^{N\left(x,y\right)}-e^{N_{\#}\left(y\right)}\right|\\
 & \leq & A'\left|y-S_{j}\left(x\right)\right|\left|x-x_{\#}\right|\,.\nonumber \end{eqnarray}
 Because $y\mapsto\psi_{1}\left(x_{\#},y\right)$ is a diffeomorphism
for small $r'>0$ independently of $\varepsilon$ small, $\psi_{1}$
is well-behaved near $\left\{ x_{n}\right\} \times r'\ww{D}$.

This allows to conclude that, for maybe smaller $r,\, r'>0$, the
image $\Psi^{N}\left(\mathcal{V}\right)$ is included, and contains,
some fibered squid sector as required. Here the width of the squid
sector does not change.

We now show that $\psi_{1}$ is one-to-one, \emph{i.e.} if $\psi_{1}\left(x,y_{1}\right)=\psi_{1}\left(x,y_{2}\right)$
then $y_{1}=y_{2}$. We have\begin{eqnarray}
\left|\psi_{1}\left(x,y_{1}\right)-\psi_{1}\left(x,y_{2}\right)\right| & = & \left|y_{1}-y_{2}\right|\left|e^{N\left(x,y_{1}\right)}+\left(y_{2}-S\left(x\right)\right)\frac{e^{N\left(x,y_{2}\right)}-e^{N\left(x,y_{1}\right)}}{y_{2}-y_{1}}\right|\\
 & \geq & K\left|y_{1}-y_{2}\right|\nonumber \end{eqnarray}
 with $K>0$, since $\left|y_{2}-S\left(x\right)\right|\leq2r'$ can
be made as small as we wish whereas $e^{N\left(x,y_{1}\right)}$ remains
far from $0$.

\medskip{}

Let us now consider $\Psi^{T}:=\left(\psi_{0},\psi_{1}\right)=\flw{\ov U_{\eps}X_{\eps}^{N}}{T_{\eps}}{}$
(with another $\psi_{1}$), and prove that it is a sectorial diffeomorphism.
Let us first deal with $\psi_{0}$. We have\begin{eqnarray}
\psi_{0}\left(x,y\right) & = & \flw{Q_{\varepsilon}P_{\varepsilon}\pp{x}}{T_{\varepsilon}\left(x,y\right)}{\left(x\right)}\\
 & = & \flw{P_{\varepsilon}\pp{x}}{\ov{T}\left(x,y,\varepsilon\right)}{\left(x\right)}\nonumber \end{eqnarray}
 where $\ov{T}$ is continuous. Because $T_{\varepsilon}\left(p_{n}\right)=0$
we can assume that $\left|\ov{T}\right|$ is bounded by some arbitrary
small $\frac{1}{2}\eta$ if $r$ is sufficiently small. Hence\begin{eqnarray}
\psi_{0}\left(\flw{P_{\varepsilon}\pp{x}}{\exp(i\theta)t}{\left(x\right)},y\right) & = & \flw{P_{\varepsilon}\pp{x}}{\exp(i\theta)t+\ov{T}\left(x,y,\varepsilon\right)}{\left(x\right)}\end{eqnarray}
 so the open set $\psi_{0}\left(V\right)$ contains a squid sector
$V\left(r_{1}\right)$ with some width $w_{0}-\eta$ and is contained
in some $V\left(r_{2}\right)$ with some width $w_{0}+\eta$. On the
other hand $\flw{Q_{\varepsilon}X_{\varepsilon}}{t}{\left(x,y\right)}=\left(f\left(x,y,t\right),g(x,y,t)\right)$
with\begin{eqnarray}
g\left(x,y,t\right) & = & y+t\left((x-x_{n})(x-x_{s})O(1)+yO(1)\right)\,.\end{eqnarray}
 Here again we can conclude that $\psi_{1}$ is well behaved near
$\left\{ x_{\#}\right\} \times r'\ww{D}$ so $\Psi^{T}\left(\mathcal{V}\right)$
is included, and contains, some fibered squid sector.

\medskip{}

It only remains to show that $\Psi^{T}$ is one-to-one. Since it is
given by the flow of $Q_{\varepsilon}X_{\varepsilon}$ it sends each
leaf of $\fol{\varepsilon}^{\pm}$ into itself. It is then sufficient
to show that its restriction to each leaf $\mathcal{L}$ is one-to-one.
Assume then that $\left(x_{j},y_{j}\right)\in\mathcal{L}$ for $j\in\left\{ 1,2\right\} $
and that $\Psi^{T}\left(x_{1},y_{1}\right)=\Psi^{T}\left(x_{2},y_{2}\right)$;
because $\mathcal{L}$ is the graph of a function $l\,:\,\Omega\to\ww{C}$
(see Proposition~\ref{pro:leaf-is-graph}) if we can show that $x_{1}=x_{2}$
then $y_{1}=y_{2}$.

There exists $K>0$ independent of $\varepsilon$ and of $\left(x_{1},x_{2}\right)$,
such that we can find a path $\gamma_{0}$ linking $x_{1}$ to $x_{2}$
within $\Omega$ with a length less than $K\left|x_{1}-x_{2}\right|$.
Let $\gamma:=l\circ\gamma_{0}$ be the lift of $\gamma_{0}$ in $\mathcal{L}$;
the application of Lemma~\ref{lem:integ-over-tangent} yields\begin{eqnarray}
\psi_{0}\left(x_{2},y_{2}\right)-\psi_{0}\left(x_{1},y_{1}\right) & = & \int_{\gamma}\left(X_{\varepsilon}\ddt\psi_{0}\right)\tau_{\varepsilon}\\
 & = & 0\,.\nonumber \end{eqnarray}
 On the one hand, according to Lemma~\ref{lem:funda}, $Z_{\varepsilon}\ddt\Psi^{T}=\left(Q_{\varepsilon}X\right)\circ\Psi^{T}$
so that, according to Lemma~\ref{lem:integ-over-tangent}, \begin{eqnarray}
0=\psi_{0}\left(x_{2},y_{2}\right)-\psi_{0}\left(x_{1},y_{1}\right) & = & \int_{\gamma}\frac{Q_{\varepsilon}\circ\psi_{0}}{U_{\varepsilon}}\frac{P_{\varepsilon}\circ\psi_{0}}{P_{\varepsilon}}dx\,.\end{eqnarray}
 On the other hand, one can find a constant $K_{1}>0$ (independent
of small $\varepsilon$) such that\begin{eqnarray}
\left|\frac{Q_{\varepsilon}\circ\psi_{0}}{U_{\varepsilon}}\frac{P_{\varepsilon}\circ\psi_{0}}{P_{\varepsilon}}\left(x,y\right)-1\right| & \leq & K_{1}\left(\left|x-x_{n}\right|+\left|y\right|\right)\end{eqnarray}
 because $Q_{\varepsilon}\left(x_{n}\right)=U_{\varepsilon}\left(p_{n}\right)$
and $T_{\varepsilon}\left(p_{n}\right)=0$. We derive\begin{eqnarray}
\left|\int_{\gamma}1dx\right| & \leq & K_{1}\int_{\gamma}\left(\left|x\right|+\left|x_{n}\right|+\left|y\right|\right)\left|dx\right|\\
\left|x_{1}-x_{2}\right| & \leq & K_{1}\left(2r+r'\right)K\left|x_{1}-x_{2}\right|\nonumber \end{eqnarray}
 which, if $r,r'>0$ are sufficiently small necessarily means $x_{1}=x_{2}$. 
\end{proof}

\subsection{Canonical first integral and spaces of leaves}

~

\begin{defn}
\label{def:secto-first-integ} We use the map $\Psi_{j,\eps}^{N,\pm}$
of Theorem~\ref{thm:secto-norma} to define the \textbf{canonical}
\textbf{sectorial first integral} $H_{j,\varepsilon}^{\pm}:=H_{j,\varepsilon}^{M}\circ\Psi_{j,\varepsilon}^{N,\pm}$
of $Z_{\varepsilon}$ over $\ssect{\varepsilon}{\pm}$ as in Section~\ref{first_int_model}.
For $\eps\in\Sigma_{0}$ it is given by: \begin{equation}
H_{j,\varepsilon}^{\pm}(x,y)=\left(y-S_{j,\eps}^{\pm}(x)\right)\exp\left(N_{j,\eps}^{\pm}\left(x,y-S_{j,\eps}^{\pm}(x)\right)\right)\prod_{j=0}^{k}(x-x_{j})^{-\frac{1}{\nu_{j}}}.\label{int_prem_syst}\end{equation}

\end{defn}
\begin{cor}
\label{cor:first-integ-secto}For each $h\in\ww{C}$ the level surface
$\left(H_{j,\varepsilon}^{\pm}\right)^{-1}\left(h\right)$ is connected
and coincides with a leaf of $\fol{j,\varepsilon}^{\pm}$. Each leaf
of $\fol{j,\varepsilon}^{\pm}$ is reached in that way. For any domain
$\mathcal{W}\subset\ssect{j,\varepsilon}{\pm}$ and any analytic function
$F\in\mathcal{O}\left(\mathcal{W}\right)$ such that $X_{\varepsilon}\ddt F=0$
(or, equivalently, $F$ is constant on each leaf of $\mathcal{F}_{j,\varepsilon}^{\pm}$)
there exists a unique holomorphic function $f\in\mathcal{O}\left(H_{j,\varepsilon}^{\pm}\left(\mathcal{W}\right)\right)$
such that $F=f\circ H_{j,\varepsilon}^{\pm}$. 
\end{cor}
\begin{proof}
The first part follows immediately from Proposition~\ref{pro:first-integ-model}
since $H_{j,\varepsilon}^{\pm}=H_{j,\varepsilon}^{M}\circ\Psi_{j,\varepsilon}^{N,\pm}$
and $\Psi_{j,\varepsilon}^{N,\pm}$ is one-to-one, sending $X_{\varepsilon}^{M}$
to $X_{\varepsilon}$. On the other hand if $X_{\varepsilon}\ddt F=0$
then Lemma~\ref{lem:integ-over-tangent} implies that $F$ is constant
on any leaf of the foliation induced by $Z_{\varepsilon}$ on $\mathcal{W}$.
Hence $F$ factors as $f\circ H_{j,\varepsilon}^{\pm}$ for some function
$f\,:\, H_{j,\varepsilon}^{\pm}\left(\mathcal{W}\right)\to\ww{C}$.
However for any $\left(\ov{x},\ov{y}\right)\in\mathcal{W}$ the restriction
of $H_{j,\varepsilon}^{\pm}$ to a small disk $\left\{ \ov{x}\right\} \times\left\{ \left|y-\ov{y}\right|<\eta\right\} $
is invertible since $X_{\varepsilon}$ is transverse to the lines
$\left\{ \ov{x}\right\} \times\ww{C}$. As a conclusion $f$ must
be holomorphic at $H_{j,\varepsilon}^{\pm}\left(\ov{x},\ov{y}\right)$
since $F$ is also holomorphic at $\left(\ov{x},\ov{y}\right)$. 
\end{proof}
The next corollary is a direct consequence of Proposition~\ref{pro:first-integ-model}
and the fact that $\Psi_{j,\varepsilon}^{\pm}$ is a sectorial diffeomorphism
:

\begin{cor}
\label{cor:space-of-leaves} The space of leaves $H_{j,\varepsilon}^{\pm}\left(\ssect{j,\varepsilon}{\pm}\right)$
is biholomorphic to $\ww{C}$. Moreover: 
\begin{enumerate}
\item The space of leaves $H_{j,\varepsilon}^{+}\left(\ssect{j,\varepsilon}{s}\right)$
of the foliation induced by $X_{\varepsilon}$ on $\ssect{j,\varepsilon}{s}$
is biholomorphic to $\ww{D}$. When $\eps$ is sufficiently small
and belongs to a good sector, the size of the conformal disk $H_{j,\varepsilon}^{+}\left(\ssect{j,\varepsilon}{s}\right)$
is bounded from below and does not vanish as $\varepsilon\to0$. 
\item The spaces of leaves over $\ssect{j,\varepsilon}{n}$ and $\ssect{j,\sigma(j),\varepsilon}{g}$
are biholomorphic to $\ww{C}$. 
\end{enumerate}
Except on $\mathcal{V}_{j,\eps}^{n}$ we can choose the conformal
coordinate on $\ww{C}$ so that $0$ corresponds to the sectorial
separatrix. 
\end{cor}

\section{\label{sec:orbital}Modulus under orbital equivalence}

Corollary~\ref{cor:eq-homog} yields a necessary and sufficient condition
for a family to be orbitally equivalent to the model family through
the existence of analytic center manifolds and global solutions to
the homological equations. We want to be more precise and to quantify
{}``how far'' we are from the existence of solutions. For this we
need a {}``canonical coordinate'' on the space of leaves over the
different sub-sectors to do the measurement. This will allow to measure
how the solutions compare in the different sectors of the intersections
$\ssect{j,\varepsilon}{+}\cap\ssect{\ell,\varepsilon}{-}$. This canonical
coordinate is provided by Corollary~\ref{cor:space-of-leaves}.

\begin{defn}
\label{leaf-coordinate} 
\begin{enumerate}
\item The space of leaves over $\ssect{j,\varepsilon}{\pm}$ is $\mathbb{C}$.
A coordinate parameterizing the leaves over $\ssect{j,\varepsilon}{\pm}$
is a first integral for the system over that domain. A first integral
vanishing on the center manifold is called a \textit{leaf-coordinate
over $\ssect{j,\varepsilon}{\pm}$}. 
\item The space of leaves over $\ssect{j,\varepsilon}{s}$ is biholomorphic
to $\ww{D}$. A first integral vanishing on the center manifold is
called a \textit{leaf-coordinate over $\ssect{j,\varepsilon}{s}$}
if it extends to a leaf-coordinate over $\ssect{j,\varepsilon}{-}$. 
\end{enumerate}
\end{defn}
\begin{lem}
~ 
\begin{enumerate}
\item Given a leaf-coordinate over $\ssect{j,0}{\pm}$ for $\eps=0$, then
for each good sector $W\subset\Sigma_{0}$, as in Definition~\ref{good_covering},
the leaf-coordinate over $\ssect{j,\varepsilon}{\pm}$ can be chosen
to depend analytically on $\eps$ and such that its limit for $\eps\to0$
is the chosen leaf-coordinate over $\ssect{j,0}{\pm}$. 
\item On $\ssect{j,\sigma(j),\varepsilon}{g}$, $\ssect{j,\varepsilon}{s}$
and $\ssect{j,\varepsilon}{\pm}$, the only changes of leaf-coordinates
are the linear maps. On $\ssect{j,\varepsilon}{n}$ they are the affine
maps. 
\end{enumerate}
\end{lem}
\begin{proof}
(1) The canonical first integral $H_{j,\varepsilon}^{\pm}$ (see Definition~\ref{def:secto-first-integ})
is one leaf-coordinate which has the required analytic dependence
in $\eps$.

(2) For fixed $\eps$ any other leaf-coordinate on $\ssect{j,\eps}{\pm}$
is the composition of $H_{j,\eps}^{\pm}$ by an analytic diffeomorphism
$\varphi_{\eps}$ as stated in Corollary~\ref{cor:first-integ-secto}.
A possible choice for the leaf-coordinate is thus $\varphi_{\eps}\circ H_{j,\varepsilon}^{\pm}$.
The only global diffeomorphisms of $\mathbb{C}$ are the affine maps.
Moreover $\ssect{j,\sigma(j),\varepsilon}{g}$ is attached to one
point of saddle type. Its center manifold is unique and corresponds
to the origin in the leaf-coordinate. The same is true if the point
is a saddle-node as we restrict to one of its saddle sectors. Both
spaces of leaves are $\mathbb{C}$ and the only global diffeomorphisms
of $\mathbb{C}$ are the affine maps. Those preserving the origin
are the linear maps. 
\end{proof}

\subsection{The first part of the orbital modulus}

\begin{thm}
On any sector $\ssect{j,\varepsilon}{n}$ the change of leaf-coordinate
from $\ssect{j,\varepsilon}{-}$ to $\ssect{j+1,\varepsilon}{+}$
is an affine map $\psi_{j,\eps}^{\infty}$. If $\varepsilon$ belongs
to some good sector $W\subset\Sigma_{0}$ then $\psi_{j,\eps}^{\infty}$
depends analytically on $\eps$ and its continuous limit for $\eps\to0$
is $\psi_{j,0}^{\infty}$. For a suitable choice of the leaf-coordinate
one can choose \begin{eqnarray*}
\left(\psi_{j,\varepsilon}^{\infty}\right)'\left(0\right) & = & e^{2i\pi a\left(\varepsilon\right)/k}\,.\end{eqnarray*}

\end{thm}
\begin{proof}
This follows simply from the fact that the change of leaf-coordinate
is a global diffeomorphism of $\mathbb{C}$. It depends analytically
on $\eps$ and has the right limit for $\eps=0$ as soon as the leaf-coordinate
does. 
\end{proof}

\subsection{The cohomological equation over $\ssect{j,\varepsilon}{+}\cup\ssect{j,\varepsilon}{-}$
and the second part of the orbital modulus}

Take $\varepsilon\in\Sigma_{0}\cup\left\{ 0\right\} $. Let \begin{equation}
\ssect{j,\varepsilon}{}:=\ssect{j,\varepsilon}{+}\cup\ssect{j,\varepsilon}{-}.\end{equation}
 On $\mathcal{V}_{j,\eps}$ we consider the change of coordinates
$(x,y)\mapsto\Phi_{j,\eps}(x,y)=\left(x,y-S_{j,\eps}(x)\right)$ and
set $X_{j,\eps}=\left(\Phi_{j,\eps}\right)_{*}X_{\eps}$. Let $\tilde{R}_{j,\eps}$
be defined on $\ssect{j,\varepsilon}{}$ as in \eqref{eq-norm}. Taking
$p\in\ssect{j,\varepsilon}{s}$ we consider \begin{equation}
L_{j,\eps}(p):=\int_{\gamma_{j,\eps}^{s}(p)}\tilde{R}_{j,\eps}\,\tau_{\eps},\label{function_L}\end{equation}
 where $\gamma_{j,\eps}^{s}(p)$ is defined in Definition~\ref{def:asy_path}.

\begin{prop}
The function $L_{j,\eps}$ in \eqref{function_L} is constant on each
leaf of $\fol{j,\varepsilon}^{\pm}$. In the leaf-coordinate it is
a holomorphic map $\phi_{j,\eps}^{0}$, vanishing at the origin. On
open sets where the leaf-coordinate is analytic in $\eps$ it depends
analytically on $\eps$. It has the limit $L_{j,0}$ ($\phi_{j,0}^{0}$
in the leaf-coordinate) when $\eps\to0$. 
\end{prop}
\begin{defn}
For a prepared family of vector fields $X_{j,\eps}$ of the form (\eqref{eq_X1})
we consider the functions $L_{j,\eps}$ of \eqref{function_L} and
the associated functions $\phi_{j,\eps}^{0}$ in the leaf-coordinate
over $\ssect{j,\varepsilon}{s}$.

For each value of $\eps\in\Sigma_{0}\cup\left\{ 0\right\} $ and each
associated set of canonical sectors we have defined a ($2k+1$)-tuple
$\mathcal{N}_{\eps}=\left(a,\psi_{0,\eps},\dots,\psi_{k-1,\eps}^{\infty},\phi_{0,\eps}^{0},\dots,\phi_{k-1,\eps}^{0}\right)$.
This ($2k+1$)-tuple depends on a choice of leaf-coordinates over
the sectors $\ssect{j,\varepsilon}{\pm}$. A different choice of leaf-coordinates
over the same canonical sectors also depending analytically on $\eps$
with a continuous limit for $\eps\to0$ yields to a different ($2k+1$)-tuple
$\overline{\mathcal{N}}_{\eps}=\left(a,\overline{\psi}_{0,\eps}^{\infty},\dots,\overline{\psi}_{k-1,\eps}^{\infty},\overline{\phi}_{0,\eps}^{0},\dots,\overline{\phi}_{k-1,\eps}^{0}\right)$.
They are related by the equivalence relation \begin{equation}
\mathcal{N}_{\eps}\sim\overline{\mathcal{N}}_{\eps}\Longleftrightarrow\left(\exists c_{\varepsilon}\in\mathcal{O}_{b}\left(W_{i}\right)\right)\left(\forall j\right)\,\,\,\,\begin{cases}
\psi_{j,\eps}^{\infty}(c_{\varepsilon}h)=c_{\varepsilon}\overline{\psi}_{j,\eps}^{\infty}(h)\\
\phi_{j,\eps}^{0}(c_{\varepsilon}h)=\overline{\phi}_{j,\eps}^{0}(h).\end{cases}\label{mathcal_N_reduit}\end{equation}
 In order to take into account that changes of coordinates and parameters
of the form \eqref{sym_prep} transform a prepared family into a prepared
family we enlarge the equivalence relation \eqref{mathcal_N_reduit}.
Let $\overline{\mathcal{N}}_{\ov{\eps}}=\left(\overline{a},\overline{\psi}_{0,\ov{\eps}}^{\infty},\dots,\overline{\psi}_{k-1,\ov{\eps}}^{\infty},\overline{\phi}_{0,\ov{\eps}}^{0},\dots,\overline{\phi}_{k-1,\ov{\eps}}^{0}\right)$
\begin{equation}
\mathcal{N}_{\eps}\sim\overline{\mathcal{N}}_{\ov{\eps}}\Longleftrightarrow\left(\exists c_{\varepsilon}\in\mathcal{O}_{b}\left(W_{i}\right)\right)\left(\exists m\in\ww{Z}/k\right)\left(\forall j,h,\varepsilon\right)\,\,\,\,\begin{cases}
\eps_{\ell}=\exp(-2\pi im(\ell-1)/k)\ov\eps_{\ell}\\
a\left(\varepsilon\right)=\overline{a}\left(\overline{\varepsilon}\right)\\
\psi_{j+m,\eps}^{\infty}(c_{\varepsilon}h)=c_{\varepsilon}\overline{\psi}_{j,\ov{\eps}}^{\infty}(h)\\
\phi_{j+m,\eps}^{0}(c_{\varepsilon}h)=\overline{\phi}_{j,\ov{\eps}}^{0}(h).\end{cases}\label{mathcal_N}\end{equation}
 Note that a $2k$-tuple $\mathcal{N}_{\eps}$ depends on a good sector
$W_{i}$ in $\eps$ space for which we can construct an adequate set
of squid sectors with fixed good angle $\theta$. In order to emphasize
this dependence we will note \[
\mathcal{N}_{\eps}^{i}:=\left(a,\psi_{0,\eps}^{\infty,i},\dots,\psi_{k-1,\eps}^{\infty,i},\phi_{0,\eps}^{0,i},\dots,\phi_{k-1,\eps}^{0,i}\right).\]
 Given a good covering $\{W_{i}\}_{1\leq i\leq d}$ of $\Sigma_{0}$
in $\eps$-space we have $d$-tuples $(\mathcal{N}_{\eps}^{i})_{1\leq i\leq d}$. 
\end{defn}
\begin{thm}
\label{orbital-modulus} Given a germ of prepared family $X_{\eps}$
of the form (\ref{eq_X0}) a good covering $\{W_{i}\}_{1\leq i\leq d}$
of $\Sigma_{0}$ in $\eps$-space, the $d$ families of equivalence
classes of $2k$-tuples \[
\mathcal{N}_{\eps}^{i}=\left\{ \left(a,\psi_{0,\eps}^{\infty,i},\dots,\psi_{k-1,\eps}^{\infty,i},\phi_{0,\eps}^{0,i},\dots,\phi_{k-1,\eps}^{0,i}\right)\right\} /\sim,\]
 is a complete modulus of analytic classification for the prepared
family $X_{\eps}$ under orbital equivalence. Moreover $\mathcal{N}_{\eps}^{i}$
can be chosen to depend analytically on $\eps\in W_{i}$ and such
that its limit for $\eps\to0$ is a given $\mathcal{N}_{0}$. 
\end{thm}
We postpone the proof of the theorem till Section~\ref{Proofs}.

\begin{thm}
A complete modulus of analytic classification under orbital equivalence
of a germ of an analytic family of vector fields unfolding a saddle-node
of codimension $k$ is given by the modulus of an associated prepared
family. 
\end{thm}
\begin{cor}
\label{cor:orb_global_factorization}Let $\left(Z_{\varepsilon}\right)_{\varepsilon}$
and $\left(\overline{Z}_{\ov{\varepsilon}}\right)_{\ov{\varepsilon}}$
be two orbitally equivalent prepared families. Then there exists an
equivalence of the form $R_{m}\circ\Psi_{\varepsilon}^{N}$ where
$R_{m}$ is the rotation of angle $\frac{2i\pi m}{k}$ and a change
of parameter $\varepsilon_{\ell}=\exp(-2\pi im(\ell-1)/k)\ov\eps_{\ell}$.
The change of coordinates $\Psi_{\varepsilon}^{N}$ preserves the
$x$-variable and is a conjugacy between $\left(R_{m}^{*}\left(X_{\varepsilon}\right)\right)_{\ov{\varepsilon}}$
and $\left(\overline{X}_{\overline{\varepsilon}}\right)_{\overline{\varepsilon}}$. 
\end{cor}
This statement will result directly from the proof of Theorem~\ref{orbital-modulus}
given further below.

\section{\label{sec:conjugacy}Modulus under conjugacy}

Two families of vector fields can only be conjugate if they are orbitally
equivalent. The modulus of an analytic family under conjugacy is constructed
by adding a time part to the modulus $(\mathcal{N}_{\eps})$ of orbital
equivalence.

\subsection{Time-part of the modulus}

We consider the family $Z_{\eps}$ given in \eqref{unfold3}. Taking
$p\in\ssect{j,\varepsilon}{s}$ we construct $\gamma_{j,\eps}^{\infty}(p)$
and we consider \begin{equation}
T_{j,\eps}(p)=\int_{\gamma_{j,\eps}^{s}(p)}\left(\frac{1}{U_{\eps}}-\frac{1}{Q_{\eps}}\right)d\tau_{j,\eps},\label{function_T}\end{equation}
 where $\gamma_{j,\eps}^{s}(p)$ is introduced in Definition~\ref{def:asy_path}.

\begin{prop}
\label{pro:period_in_leaf}The function $T_{j,\eps}(p)$ in \eqref{function_T}
depends only on the leaf. In the leaf-coordinate it is a holomorphic
map $\tilde{T}_{j,\varepsilon}$. 
\end{prop}
\begin{proof}
This follows from Theorem~\ref{thm:secto-solve} and Proposition~\ref{pro:intersect_secto-trivial}. 
\end{proof}
\begin{defn}
For the vector field $Z_{\eps}$ of the form \eqref{unfold3} and
$\eps$ in a good sector $W_{i}$ in parameter space we consider the
functions $T_{j,\eps}$ of \eqref{function_T} and the associated
functions $\tilde{T}_{j,\varepsilon}$ in the leaf-coordinate over
$\ssect{j,\varepsilon}{s}$. We build the functions $\zeta_{j,\eps}:=\tilde{T}_{j,\varepsilon}-\tilde{T}_{j,\varepsilon}\left(0\right)$
as part of the time modulus of $Z_{\varepsilon}$.

For each value of $\eps\in W_{i}\cup\left\{ 0\right\} $ and each
associated set of canonical sectors we have defined a $\left(2k+1\right)$-tuple
\[
\mathcal{T}_{\eps}^{i}=\left(C_{0,\varepsilon},\ldots,C_{k,\varepsilon},\zeta_{0,\eps}^{i},\dots,\zeta_{k-1,\eps}^{i}\right).\]
 This $\left(2k+1\right)$-tuple depends on a choice of a leaf-coordinate
over the sectors $\ssect{j,\varepsilon}{+}$. A different choice of
leaf-coordinates over the same canonical sectors yields a different
$\left(2k+1\right)$-tuple $\overline{\mathcal{T}}_{\eps}^{i}=(\overline{C}_{0,\varepsilon},\ldots,\overline{C}_{k,\varepsilon},$
$\overline{\zeta}_{0,\eps}^{i},\dots,\overline{\zeta}_{k-1,\eps}^{i})$.
If we also take into account the changes of coordinates and parameters
of the form \eqref{sym_prep} sending a prepared family to a prepared
family and we let $\overline{\mathcal{T}}_{\ov{\eps}}^{i}=(\overline{C}_{0,\ov{\eps}},\ldots,\overline{C}_{k,\ov{\eps}},$
$\overline{\zeta}_{0,\ov{\eps}}^{i},\dots,\overline{\zeta}_{k-1,\ov{\eps}}^{i})$,
we introduce the following equivalence relation\begin{eqnarray*}
\left(\mathcal{N}_{\varepsilon}^{i},\mathcal{T_{\varepsilon}}^{i}\right)\sim\left(\overline{\mathcal{N}}_{\ov{\eps}}^{i},\overline{\mathcal{T}}_{\ov{\eps}}^{i}\right) & \Longleftrightarrow\mathcal{N}_{\varepsilon}^{i}\sim\overline{\mathcal{N}}_{\ov{\eps}}^{i}\mbox{ and for the same }c_{\eps}^{i}\mbox{ and }m\,\,:\, & \begin{cases}
C_{j,\eps}e^{2i\pi mj/k} & =\overline{C}_{j,\ov{\eps}}\\
\zeta_{j+m,\eps}^{i}(c_{\eps}^{i}h) & =\overline{\zeta}_{j,\ov{\eps}}^{i}(h)\,.\end{cases}\end{eqnarray*}
 where the constants $c$ and $m$ are the same as in \eqref{mathcal_N}. 
\end{defn}
\begin{thm}
\label{conjugacy-modulus} Given a prepared family $Z_{\eps}$ of
the form (\ref{eq_Z1}) and a good covering $\{W_{i}\}_{1\leq i\leq d}$
of $\Sigma_{0}$ in $\eps$-space, the $d$ families of equivalence
classes of $\left(4k+2\right)$-tuples \begin{equation}
\left\{ \left({a,\psi}_{0,\eps}^{\infty,i},\dots,\psi_{k-1,\eps}^{\infty,i},\phi_{0,\eps}^{0,i},\dots,\phi_{k-1,\eps}^{0,i},C_{0,\varepsilon},\ldots,C_{k,\varepsilon},\zeta_{0,\eps}^{i},\dots,\zeta_{k-1,\eps}^{i}\right)_{\eps\in W_{i}}\right\} /_{\sim},\label{4_tuple}\end{equation}
 is a complete modulus of analytic classification for the family $Z_{\eps}$
under conjugacy. Moreover $\left(\mathcal{N}_{\eps}^{i},\mathcal{T}_{\varepsilon}^{i}\right)$
can be chosen to depend analytically on $\eps\in W_{i}$ and such
that the limit for $\eps\to0$ is a chosen $\left(\mathcal{N}_{0},\mathcal{T}_{0}\right)$. 
\end{thm}
We postpone a more precise statement of the theorem and the proof
till Section~\ref{Proofs}.

\begin{thm}
A complete modulus of analytic classification under conjugacy of a
germ of analytic family of vector fields unfolding a saddle-node of
codimension $k$ is given by the modulus of an associated prepared
family. 
\end{thm}
\begin{cor}
\label{cor:conj_global_factorization}Let $\left(Z_{\varepsilon}\right)_{\varepsilon}$
and $\left(\overline{Z}_{\ov{\eps}}\right)_{\ov{\eps}}$ be conjugate
prepared families. Then there exists a conjugacy of the form $R_{m}\circ\Psi_{\varepsilon}^{N}\circ\Psi_{\varepsilon}^{T}$,
where each $\left(\Psi_{\varepsilon}^{\#}\right)_{\varepsilon}$ is
an analytic family of diffeomorphisms, $R_{m}$ is the rotation of
angle $\frac{2i\pi m}{k}$ and $\eps_{\ell}=\ov\eps_{\ell}\exp(-2\pi im(\ell-1)/k)$.
The change of coordinates $\Psi_{\varepsilon}^{N}$ preserves the
$x$-variable and is an orbital equivalence between $\left(R_{m}^{*}\left(X_{\varepsilon}\right)\right)_{\ov{\eps}}$
and $(\overline{X}_{\ov{\eps}})_{\ov{\eps}}$. The change of coordinates
$\Psi_{\varepsilon}^{T}$ is given by the flow $\Phi_{Z_{\eps}}^{T_{\eps}}$
of $Z_{\varepsilon}$ for the time $T_{\eps}$, where $\left(x,y,\varepsilon\right)\mapsto T_{\varepsilon}\left(x,y\right)$
is holomorphic near $\left(0,0,0\right)$. 
\end{cor}
We postpone the proof of this result till the end of Section~\eqref{Proofs}.

\subsection{Global symmetries}

Theorem~\eqref{orbital-modulus} and Theorem~\eqref{conjugacy-modulus}
will ultimate rely on the following classification of global symmetries,
described as follows :

\begin{prop}
\label{pro:facto_sym} Let $W$ be a good sector. Any germ $\left(\chi_{\varepsilon}\right)_{\varepsilon\in W\cup\left\{ 0\right\} }$
of analytic family of symmetries of $\left(Z_{\varepsilon}\right)_{\varepsilon\in W\cup\left\{ 0\right\} }$
over $\left(\ssect{\varepsilon}{}\right)_{\varepsilon\in W\cup\left\{ 0\right\} }$,
bounded with respect to $\varepsilon\in W$, is entirely determined
by an integer $m\in\ww{Z}/k$ and two maps $\alpha\,:\,\varepsilon\mapsto\alpha\left(\varepsilon\right)$
and $\beta\,:\,\varepsilon\mapsto\beta\left(\varepsilon\right)$ holomorphic
and bounded on $W$ with continuous extension to $W\cup\left\{ 0\right\} $.
We have \begin{eqnarray}
\chi_{\varepsilon} & = & \chi_{\varepsilon}^{N}\circ\chi_{\varepsilon}^{T}\circ R_{m}\end{eqnarray}
 where $R_{m}$ is the rotation $\left(x,y\right)\mapsto\left(e^{2i\pi\frac{m}{k}}x,y\right)$,
the map $\chi_{\varepsilon}^{N}$ preserves the $x$-variable and
comes from a linear change of leaf-coordinate, and $\chi_{\varepsilon}^{T}=\flw{Z_{\varepsilon}}{\beta\left(\varepsilon\right)}{}$.
The complex number $\alpha\left(\varepsilon\right)$ represents the
linear change of leaf-coordinate induced by $\chi_{\varepsilon}^{N}$
in the sectorial spaces of leaves. For $\varepsilon$ fixed the group
of all possible $\left(m,\exp\alpha\left(\varepsilon\right)\right)$
is isomorphic to the group of changes of leaf-coordinate preserving
$\left(\mathcal{N}_{\varepsilon},\mathcal{T}_{\varepsilon}\right)$
(the symmetry group of the invariants). Moreover : 
\begin{enumerate}
\item if there exists $\varepsilon\in W\cup\left\{ 0\right\} $ such that
one of the $\psi_{j,\varepsilon}^{\infty}$ is not linear then $\alpha=0$. 
\item if there exists $\varepsilon\in W\cup\left\{ 0\right\} $ such that
for all $n\in\ww{N}_{>1}$ one of the $\phi_{j,\varepsilon}^{0}$
for some $\eps$ is not of the form $h\mapsto f\left(h^{n}\right)$
then $\alpha=0$. 
\item if all $\phi_{j,\varepsilon}^{0}$ are of the form $f_{j,\eps}(h^{n})$
for some fixed maximal $n>1$ then $\alpha=2i\pi\frac{q}{n}$ for
some fixed $q\in\ww{Z}/n$ independent on $j$. 
\item If all $\psi_{j,\varepsilon}^{\infty}$ are linear and all $\phi_{j,\varepsilon}^{0}$
vanish, then $\alpha\in\germ{\varepsilon}$. 
\end{enumerate}
The families of orbital symmetries of $\left(\mathcal{F}_{\varepsilon}\right)_{\varepsilon}$
are of the same form with $\beta$ being some germ of a holomorphic
function at $\left(0,0,0\right)$. For a fixed $\varepsilon$ the
group of all possible $\left(m,\exp\alpha\left(\varepsilon\right)\right)$
is isomorphic to the symmetry group of $\mathcal{N}_{\varepsilon}$. 
\end{prop}
\begin{proof}
We endow each sectorial space of leaves over $\ssect{j,\varepsilon}{\pm}$
with the sectorial canonical first-integral $H_{j,\varepsilon}^{\pm}$
so that $H_{j+1,\eps}^{+}=\psi_{j,\eps}^{\infty}\circ H_{j,\eps}^{-}$
and $H_{j,\eps}^{-}=H_{j,\eps}^{+}\exp\left(\phi_{j,\eps}^{0}\circ H_{j,\eps}^{+}\right)$
on $\ssect{j,\varepsilon}{n}$ and $\ssect{j,\varepsilon}{s}$ respectively.
The symmetry $\chi_{\varepsilon}\left(x,y\right)=\left(A_{\varepsilon}\left(x,y\right),B_{\varepsilon}\left(x,y\right)\right)$
induces a change of leaf-coordinate $\chi_{j,\eps}^{\pm}\,:\, h\mapsto h\exp\alpha_{j,\varepsilon}^{\pm}$,
for some $\alpha_{j,\varepsilon}^{\pm}\in\ww{C}$, and a time scaling
$\xi_{j,\varepsilon}^{\pm}=\flw{Z_{\varepsilon}}{\beta_{j,\varepsilon}^{\pm}}{}$.
Since $\chi_{\varepsilon}$ is a global map, the first observation
is that all $\beta_{j,\varepsilon}^{\pm}$ must be equal to the same
$\beta\left(\varepsilon\right)$ since $\xi_{j,\varepsilon}^{\pm}\left(x,y\right)=\left(A_{\varepsilon}\left(x,y\right),\ldots\right)$.
It is also possible to show that $\exp\alpha_{j,\varepsilon}^{\pm}=\exp\alpha\left(\varepsilon\right)$
depends only on $\varepsilon\in W$. The fact that $\alpha$ and $\beta$
depend analytically on $\varepsilon\in W$ is clear enough. As $\chi_{\varepsilon}$
is bounded on $W$ with continuous extension to $\varepsilon=0$ it
is also the case for $\alpha$ and $\beta$.

\medskip{}

For the same reason, namely because $\chi_{\varepsilon}$ is a global
object, the changes of leaf-coordinate $\chi_{j,\varepsilon}^{+}$
(\emph{resp.} $\chi_{j,\varepsilon}^{-}$) must commute with $\psi_{j,\varepsilon}^{\infty}\,:\, h\mapsto e^{2i\pi a/k}h+s_{\varepsilon}$
(\emph{resp.} $\psi_{j,\varepsilon}^{0}\,:\, h\mapsto h\exp\phi_{j,\varepsilon}^{0}\left(h\right)$).
Hence \begin{eqnarray}
\phi_{j,\varepsilon}^{0}\left(h\right) & = & \phi_{j,\varepsilon}^{0}\left(h\exp\alpha\left(\varepsilon\right)\right)\\
s_{\varepsilon} & = & s_{\varepsilon}\exp\alpha\left(\varepsilon\right)\,.\end{eqnarray}

We now discuss several cases:

\medskip{}
 \textbf{(i)} If any of the $\psi_{j,0}^{\infty}$ is nonlinear then
necessarily the same is true of $\psi_{j,\eps}^{\infty}$ for $\varepsilon\neq0$.
Thus $\exp\alpha\left(\varepsilon\right)\equiv1$ and we can choose
$\alpha\left(\varepsilon\right)=0$.

\noindent \medskip{}
 \textbf{(ii)} If all $\psi_{j,0}^{\infty}$ are linear but one $\psi_{j,\eps}^{\infty}$
is nonlinear, then it is nonlinear for all values of $\eps$ on a
dense open subset of $W$. For these values of $\eps$ we have $\alpha\left(\varepsilon\right)=0$.
By analytic continuation this is the case for all values of $\eps$
in $W$.

\noindent \medskip{}
 \textbf{(iii)} If for any $n>1$ there exists $j$ such that one
$\phi_{j,\eps}^{0}$ is not of the type $\phi_{j,\eps}^{0}(h)=f_{j}(h^{n})$
for some analytic germ of function $f_{j}$ and at least one value
of $\eps$ then $\alpha\left(\varepsilon\right)=0$.

\noindent \medskip{}
 \textbf{(iv)} If all $\phi_{\ell,\eps}^{0}$ are of the type $\phi_{\ell,\eps}^{0}(h)=f_{\ell,\eps}(h^{n})$
with $n>1$ and $n$ is maximal with this property, then the only
symmetries are of the form $\alpha\left(\varepsilon\right)=2\pi i\frac{q}{n}$.

\noindent \medskip{}
 \textbf{(v)} If all $\psi_{j,\eps}^{\infty}$ are linear and all
$\phi_{j,\eps}^{0}$ are zero then there is no constraint and $\alpha\left(\varepsilon\right)\in\mathbb{C}$.

Checking the remaining statements is straightforward. 
\end{proof}
\begin{cor}
$\left(Z_{\varepsilon}\right)_{\varepsilon}$ is orbitally equivalent
to $\left(X_{\varepsilon}^{M}\right)_{\varepsilon}$ if, and only
if, the symmetry group of $\mathcal{N}_{\varepsilon}$ is infinite
for all $\varepsilon$. 
\end{cor}

\section{Proofs of Theorems~\ref{orbital-modulus} and \ref{conjugacy-modulus}}

\label{Proofs}

\begin{defn}
\label{orbital_equivalence} Two germs of $k$-parameter analytic
families of vector fields $\left(Z_{\eps}\right)_{\varepsilon}$ (\emph{resp.}
$\left(\ov Z_{\ov{\eps}}\right)_{\overline{\varepsilon}}$) unfolding
a saddle-node of codimension $k$ at the origin for $\eps=0$ (\emph{resp.}
$\ov{\eps}=0$) are \textbf{\textit{\emph{orbitally equivalent}}}
if there exists a germ of analytic map \begin{equation}
K=(g,\Psi,\xi):\quad(\eps,x,y)\mapsto(g(\eps),\Psi(\eps,x,y),\xi(\eps,x,y))\end{equation}
 fibered over the parameter space where 
\begin{enumerate}
\item $g:\eps\mapsto\ov{\eps}=g(\eps)$ is a germ of an analytic diffeomorphism
preserving the origin; 
\item there exists a representative $\Psi_{\eps}(x,y)=\Psi(\eps,x,y)$ which
is an analytic diffeomorphism in $(\eps,x,y)$ on a small neighborhood
of the origin in $(\eps,x,y)$-space; 
\item there exists a representative $\xi_{\eps}(x,y)=\xi(\eps,x,y)$ depending
analytically on $(\eps,x,y)$ in a small neighborhood of the origin
in $(\eps,x,y)$-space with values in $\mathbb{C}_{\neq0}$; 
\item the change of coordinates $\Psi_{\eps}$ and time scaling $\xi_{\eps}$
is an equivalence between $Z_{\eps}$ and $\ov Z_{g(\eps)}$ over
a ball of small radius $r>0$ : \begin{equation}
\ov Z_{g(\eps)}\left(x,y\right)=\xi\left(\varepsilon,x,y\right)\left(\left(\Psi_{\eps}\right)_{*}Z_{\eps}\right)\left(x,y\right).\label{transform}\end{equation}

\end{enumerate}
The families are \textbf{\textit{\emph{conjugate}}} if it is possible
to choose $K=(g,\Psi,\xi)$ with $\xi\equiv1$. 
\end{defn}

\subsection{Proof of Theorem~\ref{orbital-modulus}}

\noindent If two families are orbitally equivalent under an equivalence
$\Psi$ preserving the parameter modulo a rotation of order $k$ then
they have the same modulus. Let $K$ be indeed an orbital equivalence
between two families $\left(Z_{\varepsilon}\right)_{\varepsilon}$
and $\left(\overline{Z}_{\varepsilon}\right)_{\varepsilon}$, which
we can assume to be prepared. After possibly applying a rotation of
order $k$ to one of the systems and changing the parameter accordingly
we can assume that the $x$-component of $\Psi$ is tangent to the
identity and that $g=Id$. Let $W\subset\Sigma_{0}$ be a good sector
with associated canonical sectors. Over $\ssect{j,\varepsilon}{\pm}$
the map $\Psi$ induces a change of leaf-coordinate, given by linear
and invertible maps $h\mapsto\ov{h}=c_{j,\varepsilon}^{\pm}h$. Hence
$c_{j+1,\varepsilon}^{+}\psi_{j,\varepsilon}^{\infty}\left(h\right)=\overline{\psi}_{j,\varepsilon}^{\infty}\left(c_{j,\varepsilon}^{-}h\right)$
and $\phi_{j,\varepsilon}^{0}\left(h\right)=\overline{\phi}_{j,\varepsilon}^{0}\left(c_{j,\varepsilon}^{j}h\right)$.
Because all the $(\psi_{j,\eps}^{\infty})'(\infty)$ have been normalized
and all $\phi_{j,\eps}(0)=0$ for fixed $\varepsilon$ all the $c_{j,\varepsilon}^{\pm}$
agree: $c_{j,\eps}^{\pm}=c_{\eps}$. Moreover $c_{\eps}$ depends
analytically on $\eps\in W_{i}$ with continuous limit at $\eps=0$.
Hence the moduli coincide. \medskip{}

Conversely we consider two prepared families with same modulus. We
can apply a rotation of order $k$ (and the corresponding change of
parameter) and a change of leaf-coordinate in the modulus so that
the functions $\psi_{j,\eps}^{\infty,i}$ and $\phi_{j,\eps}^{0,i}$
be exactly the same for the two families. Because a change of leaf
is given by some functions $c_{\eps}^{i}\in\mathcal{O}_{W_{i}}$ which
are bounded and bounded away from $0$ there exist common $r,r',\rho$
valid for the two families. We now look for an equivalence preserving
the parameter. The strategy is the following: we start by constructing
an equivalence between the two families on a good open covering $\{W_{i}\}$
of $\Sigma_{0}$ given in Definition~\ref{good_covering} and we
show the existence of an equivalence depending analytically on $\eps\neq0$
in each $W_{i}$ and continuously on $\eps$ in $W_{i}$ near $\eps=0$.
Using that the equivalences are bounded and the symmetries of the
system we correct to an equivalence depending analytically on $\eps$.

\medskip{}
 On each fixed open set $W_{i}$ we drop the upper index $i$. We
can suppose that the two families have the same representative of
the modulus $\mathcal{N}_{\eps}=\overline{\mathcal{N}_{\eps}}$. On
each $\ssect{j,\varepsilon}{\pm}$ we get first integrals given by
the canonical leaf-coordinates $H_{j,\eps}^{\pm}$ and $\ov H_{j,\eps}^{\pm}$.
Each first integral yields a change of coordinates on $\ssect{j,\varepsilon}{\pm}$
transforming $X_{\eps}^{M}$ into $X_{\eps}$ of the form $\Psi_{j,\eps}^{\pm}\,:\,(x,y)\mapsto(x,\tilde{y})$
where \begin{equation}
\tilde{y}=H_{j,\eps}^{\pm}\prod_{j=0}^{k}(x-x_{j})^{\frac{1}{\nu_{j}}}.\end{equation}
 Then \begin{equation}
\Psi_{j,\eps}^{\pm}\left(x,y\right)=\left(x,\left(y-S_{j,\varepsilon}\left(x\right)\right)\exp N_{j,\varepsilon}^{\pm}\left(x,y\right)\right)\label{decomposition}\end{equation}
 which is univalued on the sector. Similar changes of coordinates
$\ov\Psi_{j,\eps}^{\pm}$ exist for $\ov X_{\eps}$. We define an
equivalence between the vectors fields $X_{\eps}$ and $\ov X_{\eps}$
as \begin{equation}
\Psi_{\eps}:=\left(\ov\Psi_{j,\eps}^{\pm}\right)^{-1}\circ\Psi_{j,\eps}^{\pm}.\label{map-psi}\end{equation}
 This change of coordinates is well defined as the two vector fields
have the same modulus. On the sectors $\ssect{j,\varepsilon}{n}$
the result follows from $H_{j+1,\eps}^{+}=\psi_{j,\eps}^{\infty}\circ H_{j,\eps}^{-}$
and \eqref{decomposition}. On the sectors $\mathcal{V}_{j,\eps}^{s}$
the result follows similarly from $H_{j,\eps}^{-}=H_{j,\varepsilon}^{+}\exp\left(\phi_{j,\eps}^{0}\circ H_{j,\eps}^{+}\right)$.
On the sectors $\ssect{j,\sigma(j),\varepsilon}{g}$ the result follows
from the fact that the linear maps transforming one first integral
to the other are identical for the two families. Let us show that
they are completely determined by the $\phi_{j,\varepsilon}^{0}\left(0\right)$
and $(\psi_{j,\eps}^{\infty})'(\infty)$. Indeed, we look at the decomposition
$r\mathbb{D}=\cup_{j=0}^{k-1}(V_{j,\eps}^{+}\cup V_{j,\eps}^{-})\cup\left\{ x_{0},\dots,x_{k}\right\} $
(see for instance Figure~\ref{fig_sectors}). Making one turn in
the positive (\emph{resp.} negative) direction around each point $x_{\ell}$
of node (\emph{resp.} saddle) type yields a correspondence map $k_{\ell,\eps}$
on the leaves which is analytic when written in a leaf-coordinate.
If $x_{\ell}$ is of node type (\emph{resp.} saddle type), then this
map is a composition of some of the linear maps with some of the $\psi_{j,\eps}^{\infty}$
(\emph{resp.} $\psi_{j,\eps}^{0}$), where $\psi_{j,\eps}^{0}$ is
defined as \[
\psi_{j,\eps}^{0}(h)=h\exp(\phi_{j,\eps}^{0}(h))\]
 and is tangent to the identity. $\psi_{j,\eps}^{0}$ is a correspondence
map from $\mathcal{V}_{j,\eps}^{-}$ to $\mathcal{V}_{j,\eps}^{+}$
over $\mathcal{V}_{j,\eps}^{s}$. The multiplier of the correspondence
map $k_{\ell,\eps}$ at the fixed point $x_{\ell}$ is given by $\exp(-\frac{2\pi i}{\nu_{\ell}})$.
On the other hand it is given by the product of the multipliers of
the linear maps together with those of the maps $\left(\psi_{j,\eps}^{0}\right)'(0)$
(\emph{resp}. $(\psi_{j,\eps}^{\infty})'(\infty)$) arising in the
decomposition. This yields a system allowing to find the multipliers
of the linear maps. We need to show that this system has a unique
solution. This comes from the structure of gate sectors $V_{j,\sigma(j),\eps}^{g}$
discussed in Lemma~\ref{lemma_DS} and also studied by Oudkerk \cite{O}.
Identifying a gate sector to a segment between two singular points,
the resulting graph of the gate sectors is a tree (an explanation
follows below). Then we start solving for the multipliers of the linear
maps by the ends of the trees, where we can find one multiplier at
a time and move towards the inside of the graph, until all multipliers
are found.

Let us describe why the graph is a tree. The separating graph $\Gamma$
(the union of the separatrices from infinity) allows to divide $D'=r\ww{D}\setminus\Gamma$
into connected components, such that the intersection of the closure
of each connected component with $r\mathbb{S}^{1}$ is exactly $\partial V_{j,\eps}^{+}\cup\partial V_{\sigma(j),\eps}^{-}$,
which yielded the map $\sigma$ defined in \eqref{def_sigma}. In
each connected component of $D'$ it is possible to draw a curve joining
$\partial V_{j,\eps}^{+}$ to $\partial V_{\sigma\left(j\right),\eps}^{-}$
(Figure~\ref{fig:fig_sigma}). This curve cuts exactly one gate sector.
If we were having a cycle of gate sectors some of these curves would
cut more than one gate sector.

A second argument to show that the graph is a tree is the following.
The graph of the gate sectors has $k+1$ vertices and $k$ edges.
Moreover, from its construction, it is easy to see that it is connected.
Indeed two adjacent boundary sectors share a singular point, so their
respective attached gate sectors, each corresponding to an edge, share
a common vertex. If we have cycles in the graph, then necessarily
the number of edges should be at least as large as the number of vertices.
Hence there are no cycles.

\medskip{}
 It is of course possible to define a map $\Psi_{\eps}$ as in \eqref{map-psi}
for any value of $\eps\in\Sigma_{0}\cup\left\{ 0\right\} $. Moreover
for $\eps\in W_{i}$ it is possible to choose $\Psi_{\eps}^{i}$ depending
analytically on $\eps\neq0$ in $W_{i}$ and having the same limit
$\Psi_{0}$ for $\eps\to0$. The last step of the proof is to build
a global $\Psi_{\eps}$ depending analytically on $\eps$ on a full
neighborhood of $\eps=0$ from the $\Psi_{\eps}^{i}$ defined for
$\eps\in W_{i}$. The ideas are similar to those of the addendum of
\cite{MRR}, namely to use the symmetries of the system. Indeed on
$W_{i,i'}=W_{i}\cap W_{i'}$ \begin{equation}
\chi_{i,i',\varepsilon}:=\left(\Psi_{\eps}^{i'}\right)^{-1}\circ\Psi_{\eps}^{i}\end{equation}
 is a symmetry of $X_{\eps}$ on a full neighborhood of the origin
in $(x,y)$-space preserving the $x$-variable. These symmetries have
been described in Proposition~\ref{pro:facto_sym}. They are given
by analytic maps $\alpha_{i,i'}(\eps)$ corresponding to linear changes
of leaf-coordinate, with $\alpha_{i,i'}(0)=0$. By Proposition~\ref{pro:facto_sym}
$\alpha_{i,i'}(\eps)\equiv0$, which implies $\chi_{i,i',\eps}=id$,
except in the case where all $\psi_{j,\eps}^{\infty}$ are linear
and all $\phi_{j,\eps}^{0}\equiv0$. (If these properties are true
for some $W_{i}$ then they are true for all the other good sectors
of the covering). In the latter case the center manifold $y=S_{\eps}(x)$
is a global analytic bounded map for $\eps\in\Sigma_{0}$, hence for
all $\eps$ with $||\eps||\leq\rho$. The linear changes of leaf coordinates
are induced by linear changes $y\mapsto cy$ in the $y$-coordinate
over the model. Looking at the constructions of the functions $N_{\eps}^{i}$
and $N_{\eps}^{i'}$ (\emph{resp}. $\ov N_{\eps}^{i}$ and $\ov N_{\eps}^{i'}$)
over sectors $W_{i}$ and $W_{i}'$, it is clear that their values
coincide on the separatrices. Indeed the value of $N_{\eps}^{i}(x_{j,\#},y)$,
$\#\in\{s,n\}$, is given by the integral of $\tilde{R}_{j,\eps}$
on the segment $[0,y]$ in $\{x=x_{j,\#}\}$ (the horizontal part
of the integral from $p_{j,s}$ to $p_{j,n}$ vanishing, since inside
the center manifold). Moreover, on $x=x_{j,\#}$, $\tilde{R}_{j,\eps}(x_{j,\#},y)=R_{2,\eps}(x_{j,\#},y)$.
Hence $\alpha_{i,i'}(\varepsilon)=0$, which implies that the equivalences
$\Psi_{i,\eps}$ between $X_{\eps}$ and $\ov X_{\eps}$ defined over
the sectors $W_{i}$ glue into a global bounded equivalence $\Psi_{\eps}$
defined for $\eps\in\Sigma_{0}$.

\noindent \medskip{}
 Hence we have defined a global map $\Psi_{\eps}$ on $\Sigma_{0}$.
As it is bounded, it is possible to extend it to a full neighborhood
$W$ of the origin in $\eps$-space.\hfill{}$\Box$

\subsection{Proof of Theorem~\ref{conjugacy-modulus}}

The strategy is similar to that of Theorem~\ref{orbital-modulus}.

For the direct part, if two families are conjugate then we can bring
them to the prepared form \eqref{unfold3} with same $X_{\eps}$ and
the two families have the form $Z_{\eps}=X_{\eps}U_{\eps}$ and $\ov Z_{\eps}=X_{\eps}\ov U_{\eps}$
with same temporal normal form $Q_{\varepsilon}X_{\varepsilon}$.
Then there exists a conjugacy $\Psi_{\varepsilon}$ between these
two forms which is a symmetry of the foliation. Proposition~\ref{pro:facto_sym}
describes those maps and we obtain the existence of a holomorphic
map $T_{\varepsilon}$ such that $\flw{Z_{\varepsilon}}{T_{\varepsilon}}{}$
conjugates $Z_{\varepsilon}$ to $\overline{Z}_{\varepsilon}$. According
to Lemma~\ref{lem:funda} we have $X_{\varepsilon}\cdot T_{\varepsilon}=\frac{1}{U_{\varepsilon}}-\frac{1}{\overline{U}_{\varepsilon}}$
and according to Corollary~\ref{cor:eq-homog} we obtain, for all
$p\in\ssect{j,\varepsilon}{s}$, \begin{eqnarray}
\int_{\gamma_{j,\varepsilon}^{s}\left(p\right)}\left(\frac{1}{U_{\varepsilon}}-\frac{1}{Q_{\varepsilon}}+\frac{1}{Q_{\varepsilon}}-\frac{1}{\overline{U}_{\varepsilon}}\right)\,\tau_{\varepsilon} & = & I\left(j\right)\end{eqnarray}
 independently on $p$ (we recall that $\tau_{\varepsilon}=\frac{dx}{P_{\varepsilon}}$
is the time-form associated to $X_{\varepsilon}$). With the notations
of Proposition~\ref{pro:period_in_leaf} this implies \begin{eqnarray*}
\tilde{T}_{j,\varepsilon}\left(h\right) & = & \tilde{\overline{T}}_{j,\varepsilon}\left(h\right)+I\left(j\right)\end{eqnarray*}
 in the leaf coordinate so that $\zeta_{j,\varepsilon}-\overline{\zeta}_{j,\varepsilon}=0$.

\medskip{}

Conversely, let us suppose that two prepared families $Z_{\eps}$
and $\ov Z_{\eps}$ have the same modulus and same polynomial $P_{\eps}$.
From Theorem~\ref{orbital-modulus} we know that the two families
are orbitally equivalent and (after possibly applying a rotation of
order $k$ and the corresponding change of parameter) that the equivalence
preserves the parameter, so we can suppose that they have the form
$Z_{\eps}=U_{\eps}X_{\eps}$ and $\ov Z_{\eps}=\ov U_{\eps}X_{\eps}$.
We look for a conjugacy of the form $\Phi_{U_{\eps}X_{\eps}}^{T_{\eps}}$
where $T_{\eps}$ is holomorphic over $\mathcal{V}_{\eps}$ and satisfies
\begin{equation}
X_{\eps}\cdot T_{\eps}=\frac{1}{U_{\eps}}-\frac{1}{\ov U_{\eps}}\,\,.\end{equation}

As before we consider a good open covering $\{W_{i}\}_{1\leq i\leq d}$
of $\Sigma_{0}$ and we construct analytic functions $T_{\eps}^{i}$
depending analytically on $\eps$ over $W_{i}$ and having the same
limit when $\eps\to0$ inside $W_{i}$. For a fixed leaf-coordinate
over $\ssect{j,\varepsilon}{\pm}$ we have $\zeta_{j,\varepsilon}=\overline{\zeta}_{j,\varepsilon}$
so that, for any $p\in\ssect{j,\varepsilon}{s}$, \begin{equation}
\int_{\gamma_{j,\eps}^{s}(p)}\left(\frac{1}{U_{\eps}}-\frac{1}{\ov U_{\eps}}\right)\;\tau_{\eps}=I\left(j\right)\end{equation}
 where $I\left(j\right)$ is constant. Applying once more Corollary~\ref{cor:eq-homog},
while using the fact that the graph of gate sectors is a tree, gives
the existence of $T_{\varepsilon}^{i}$. \medskip{}

\noindent The last step is to build from the $T_{\eps}^{i}$ a global
$T_{\eps}$ depending analytically on $\eps\in\Sigma_{0}$. We proceed
as in Theorem~\ref{orbital-modulus} and correct $\flw{Z_{\varepsilon}}{T_{\varepsilon}^{i}}{}$
by composing with a symmetry $\chi_{i,\varepsilon}$ of $Z_{\eps}$
over $W_{i}$. As described in Proposition~\ref{pro:facto_sym} any
family of symmetries over $W_{i}\cap W_{i'}$ which does not exchange
leaves are given by analytic maps $\varepsilon\mapsto\beta_{i,i'}\left(\varepsilon\right)$
with $\beta_{i,i'}(0)=0$. We want to find functions $\beta_{i}$
such that $\beta_{i,i'}=\beta_{i}-\beta_{i'}$. Of course this is
the first Cousin problem, which is solvable since $\Sigma_{0}$ is
a Stein manifold but this is not sufficient as we need to show that
the $\beta_{i}$ are bounded. So we proceed as follows. On each sector
$W_{i}$ we have constructed a family of functions $T_{\eps}^{i}$
defined on $r\ww{D}\times r'\ww{D}$ and conjugating $Z_{\eps}$ and
$\ov Z_{\eps}$ for $\eps\in W_{i}$. Hence these functions differ
from a constant $\beta_{i,i'}\left(\varepsilon\right)$ for each $\varepsilon\in W_{i}\cap W_{i'}$
on $W_{i}\cap W_{i'}$. This constant is calculated for instance as
$T_{\eps}^{i}(0,0)-T_{\eps}^{i'}(0,0)$. We let $\beta_{i}(\eps):=T_{\eps}^{i}(0,0)$.
Defining \[
T_{\eps}:=T_{\eps}^{i}-\beta_{i}(\eps),\]
 yields the required map $T_{\eps}$ so that $\Phi_{Z_{\eps}}^{T_{\eps}}$
conjugates $Z_{\eps}$ with $\ov Z_{\eps}$. \hfill{}$\Box$

\subsection{Proof of Corollaries~\ref{cor:orb_global_factorization} and \ref{cor:conj_global_factorization}}

Since $\left(Z_{\varepsilon}\right)_{\varepsilon}$ and $\left(\overline{Z}_{\overline{\varepsilon}}\right)_{\overline{\varepsilon}}$
are conjugate they have the same moduli, thus are orbitally equivalent.
From the proof of Theorem~\ref{orbital-modulus} we get $\Psi_{\varepsilon}^{N}$
and $R_{m}$, while the proof of Theorem~\ref{conjugacy-modulus}
done just above provides us with $\Psi_{\varepsilon}^{T}=\Phi_{Z_{\eps}}^{T_{\eps}}$.\hfill{}$\Box$

\section{Perspectives, applications and questions}

\label{sec:Questions}

\subsection{Reading the dynamics from the modulus}

\label{subsec:reading}

The modulus allows to read the dynamics of the system. A first case
was presented in Example~\ref{ex_center_man}. This example was not
finished. Indeed we gave sufficient conditions for the stable manifold
of $(x_{2},0)$ to coincide with the weak invariant manifold of $(x_{0},0)$,
but they were not necessary. We now can give the necessary and sufficient
condition.

\begin{figure}
\includegraphics[width=8cm]{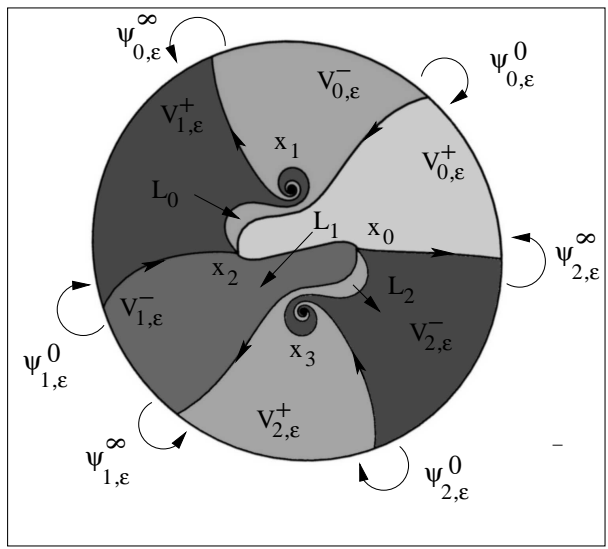}

\caption{\label{fig:dynamics}Interpretation of the orbital invariants in
terms of the global dynamics around singular points. The linear transformations
$L_{j}$ are the changes of leaf-coordinates over the gate sectors.
The direction of the arrows yield their direction. For instance $L_{1}$
is the change from the leaf-coordinate on $V_{0,\eps}^{+}$ to the
one on $V_{1,\eps}^{-}$.} 
\end{figure}

We also discuss other cases coming from Figure~\ref{fig:dynamics}
which is the Figure~\ref{fig:center} of Example~\ref{ex_center_man}.
We introduce the transition maps: $\psi_{j,\eps}^{0}:\mathcal{V}_{j,\eps}^{-}\rightarrow\mathcal{V}_{j,\eps}^{+}$
defined over $\mathcal{V}_{j,\eps}^{s}$ by \begin{equation}
\psi_{j,\eps}^{0}\left(h\right):=h\exp(\phi_{j,\eps}^{0}(h)).\end{equation}
 Note that the $\psi_{j,\eps}^{\infty}$ (\emph{resp.} $\psi_{j,\eps}^{0}$)
are transition maps when we move in the anticlockwise (\emph{resp.}
clockwise) direction. A hidden motivation for this choice is that
it is the direction for which the dynamics of the holonomy is going
forward: the iterates of a point under the holonomy map move in that
direction (\cite{R2}). We locate the area of action of each transition
map in Figure~\ref{fig:dynamics}.

\begin{example}
\label{ex:11_1} End of Example~\ref{ex_center_man}. We introduce
the Lavaurs maps $L_{i}$ (see Figure~\ref{fig:dynamics}). These
maps are the changes of leaf-coordinates over the gate sectors. 
\begin{description}
\item [{{{{\textmd{(5)}}}}}] To give a necessary and sufficient condition
for the stable manifold of $(x_{2},0)$ to coincide with the weak
invariant manifold of $(x_{0},0)$ we need to characterize the weak
invariant manifold of $(x_{0},0)$. It is the only leaf which is not
ramified at the point. Hence it is the fixed point of the first return
map of leaves when one makes a positive turn around $(x_{0},0)$.
We choose to start this return map in a sector where the stable manifold
of $(x_{2},0)$ corresponds to the zero leaf-coordinate.

So the necessary condition is given for instance by : \[
L_{1}\circ\psi_{2,\eps}^{\infty}\circ L_{2}\circ\psi_{1,\eps}^{\infty}(0)=0.\]

\item [{{{{\textmd{(6)}}}}}] We can for instance read the dynamics
of $(x_{2},0)$ for the particular values of the parameters for which
it is a saddle point (the ratio of eigenvalues is in $\mathbb{R}_{\leq0}$).
Then the return map for leaves if given by : \[
k_{2,\eps}=\psi_{0,\eps}^{0}\circ L_{0}\circ\psi_{1,\eps}^{0}\circ L_{1}.\]
 In particular $(x_{2},0)$ is orbitally linearizable if and only
if $k_{2,\eps}$ is linearizable. The parametric resurgence phenomenon
described in \cite{R2} also appears here. Indeed, let us recall that
$(\psi_{j,\eps}^{0})'(0)=1$ which comes from the fact that $\phi_{j,\eps}^{0}(0)=0$.
For instance let $\eps_{n}$ be a sequence of values of $\eps$ such
that $\lim\eps_{n}=0$, $L_{0}$ and $L_{1}$ are fixed and $k_{2,\eps}'(0)=\exp(2\pi i\frac{p}{q})$.
If \[
k_{2,0}=\psi_{0,0}^{0}\circ L_{0}\circ\psi_{1,0}^{0}\circ L_{1}\]
 is non-linearizable then so is the case for $k_{2,\eps_{n}}$ as
soon as $n$ is sufficiently large. The non-linearizability of $k_{2,0}$
can be seen from the non vanishing of a coefficient of the normal
form. In the particular case where $k_{2,\eps}'(0)=1$, this is the
case as soon as $k_{2,0}$ is nonlinear. 
\end{description}
\end{example}

\subsection{Extending the discussion beyond $\Sigma_{0}$}

We have made an extensive description of the family $Z_{\eps}$ for
the values of $\eps$ in $\Sigma_{0}$. Such a description can also
be made for the other values of $\eps$ and the paper of Douady and
Sentenac \cite{DS} already contains the necessary adjustments. Indeed
here, when $\eps\in\Sigma_{0}$ the sectors are constructed as strips
in $z$-space. When $\eps\notin\Sigma_{0}$ and some of the singular
points are saddle-nodes such a decomposition is $2k$ sectors still
exist, but some of the strips in $z$-plane are replaced by half-spaces.
Each sector is again adherent to two points, one of saddle type and
one of node type, using the remark that a saddle-node can be of saddle
type or of node type when restricted to a domain over a sector. The
center manifold theorem (Theorem~\ref{center_manifold}) is still
valid in this context and we get $k$ center manifolds on $k$ sectors
attached to sectors $\partial V_{j,\eps}^{\pm}$ of the boundary $|x|=r$.
The construction of asymptotic paths (Theorem~\ref{thm:secto-trivial})
can be performed in full generality. Similarly Theorem~\ref{thm:secto-solve},
where we solve cohomological equations on sectors, remains true for
all values of $\varepsilon$.

In our discussion we have worked with a finite open covering $\mathcal{W}=\{W_{i}\}$
of $\Sigma_{0}$. In this way we have avoided discussing the stratification
of the complement of $\Sigma_{0}$. The $W_{i}$ are constructed as
cones on open sets in the sphere $\{||\eps||=\rho\}$. A subset of
$\mathcal{W}$ is necessary to cover the neighborhood of a value $\eps_{0}\notin\Sigma_{0}$.
The spatial organization of these sectors around $\eps_{0}$ is an
interesting question for a future work. For instance, if $\eps_{0}$
is a regular point of a stratum of codimension $1$, then the intersection
of these sectors with a section transverse to the stratum gives an
open covering of the section minus $\eps_{0}$.

\subsection{The link with the holonomy of the strong separatrix}

In \cite{RC} the modulus of analytic equivalence under orbital equivalence
of a generic 1-parameter family unfolding a planar vector field with
a resonant saddle is given in terms of the modulus of analytic equivalence
of the family of holonomy maps corresponding to one separatrix. Then
a time part is added to the modulus to give a modulus of analytic
equivalence under conjugacy. The approach of \cite{RC} could obviously
have been extended to the case of the saddle-node. We have preferred
a geometric approach, based on the asymptotic homology of the leaves,
as it is the fact that the space of leaves over the canonical sectors
is $\mathbb{C}$ which yields that the maps $\psi_{j,\eps}^{\infty}$
of the modulus are affine maps.

If we consider a section $y=1$ of the strong separatrix, then it
can be proved as in \cite{R2} that all leaves over a canonical sector
$\ssect{j,\varepsilon}{\pm}$ intersect $y=1$ and that different
points of intersection belong to the same orbit of the holonomy map.
So we have a correspondence between the space of leaves over the canonical
sectors and the orbit spaces of the holonomy maps. Hence two germs
of generic families of vector fields with a saddle-node of codimension
$k$ at the origin and same formal parameter are orbitally equivalent
if and only if the families of unfoldings of the holonomies of their
strong separatrices are conjugate. The same is true for conjugacy
if we add to the holonomies the times needed to compute them by following
the flow of the vector field.

It is worth here citing the work \cite{R} of Javier Ribón which studies
the modulus of analytic classification of 1-parameter families unfolding
of resonant diffeomorphisms.

\subsection{Questions and directions for future research}

\begin{enumerate}
\item The most important question coming from our work is to identify the
modulus space, both for the problem of orbital equivalence and for
the problem of conjugacy. The dependence on $\eps$ of the components
of the moduli is a highly non trivial question. In an upcoming paper
with Reinhard Sch\"{a}fke we propose to prove that in the case $k=1$
the moduli $\phi_{0,\varepsilon}^{0,i}$, $\psi_{0,\varepsilon}^{\infty,i}$
and $\zeta_{0,\varepsilon}^{i}$ represent $\frac{1}{2}$-sums of
formal power series $\sum A_{n}\left(h\right)\varepsilon^{n}$ with
$A_{n}$ holomorphic, as was earlier suspected. For a given value
of $k$, this requires in particular to describe the relationships
between the different $\mathcal{N}_{\varepsilon}^{i}$ (\emph{resp.}
$(\mathcal{N}_{\varepsilon}^{i},\mathcal{T}_{\varepsilon}^{i})$)
on all intersections $W_{i}\cap W_{i'}$ of two good sectors in $\eps$-space. 
\item As for the time part of the modulus, the problem addressed is whether
it is possible to {}``unfold'' a result of \cite{T2} (an adaptation
of Ramis-Sibuya theorem) stating that, given $k$ functions holomorphic
on the space of leaves of the canonical sectors $\ssect{j,0}{s}$,
it is possible to find a holomorphic function $G_{0}$ such that the
obstructions to solve $Z_{0}\cdot F_{0}=G_{0}$ are precisely the
given functions. 
\end{enumerate}

\section*{Acknowledgements}

We are very grateful to Colin Christopher and Adrien Douady for lengthy
discussions on the preparation part.


\begin{thebibliography}{10}
\bibitem{CMR2}C.~Christopher, P.~Marde\v{s}i\'{c} and C.~Rousseau,
{}``Normalizability, synchronicity and relative exactness for vector
fields in $\mathbb{C}^{2}$'', \emph{Journal of Dynamical and Control
Systems}, \textbf{10} (2004), 501--525.

\bibitem{DS} A.~Douady and P.~Sentenac, {}``Champs de vecteurs
polynomiaux sur $\mathbb{C}$'', prépublication, Paris 2005.

\bibitem{G}A.~A.~Glutsyuk, {}``Confluence of singular points and
nonlinear Stokes phenomenon'', \emph{Trans. Moscow Math. Soc.}, \textbf{62}
(2001), 49--95.

\bibitem{HKM}M.~Hukuhara, T.~Kimura and T.~Matuda, {}``Équations
différentielles ordinaires du premier ordre dans le champ complexe'',
\emph{Pub. Math. Soc. Japan}, vol.7, (1961).

\bibitem{I}Y.~Ilyashenko, {}``Nonlinear Stokes phenomena'', in
\emph{Nonlinear Stokes phenomena}, Y. Ilyashenko editor, Advances
in Soviet Mathematics, vol. 14, Amer. Math. Soc., Providence, RI,
(1993), 1--55.

\bibitem{IY} Y.~Ilyashenko and S.~Yakovenko, {}``Stokes phenomena
in smooth classification problems'', in \emph{Nonlinear Stokes phenomena},
Y. Ilyashenko editor, Advances in Soviet Mathematics, vol. 14, Amer.
Math. Soc., Providence, RI, (1993), 235--287.

\bibitem{K}V.~Kostov, {}``Versal deformations of differential forms
of degree $\alpha$ on the line'', \emph{Functional Anal. Appl.}
\textbf{18} (1984), 335--337.

\bibitem{MRR}P.~Marde\v{s}i\'{c}, R.~Roussarie and C.~Rousseau,
{}``Modulus of analytic classification for unfoldings of generic
parabolic diffeomorphisms'', \emph{Moscow Mathematical Journal},
\textbf{4} (2004), 455--502.

\bibitem{MR1} J.~Martinet and J.-P.~Ramis, Problèmes de modules
pour des équations différentielles non linéaires du premier ordre,
\emph{Publ. Math., Inst. Hautes Etud. Sci.} \textbf{55} (1982), 63--164.

\bibitem{MV} Y.~I.~Meshcheryakova and S.~M.~Voronin, {}``Analytic
classification of typical degenerate elementary singular points of
the germs of holomorphic vector fields on the complex plane, \emph{Izvestiya
V.U.Z. Mat.} \textbf{46} (2002), 11-14.

\bibitem{O}R.~Oudkerk, {}``The parabolic implosion for $f_{0}(z)=z+z^{\nu+1}+O\left(z^{\nu+2}\right)$'',
thesis, University of Warwick, (1999).

\bibitem{R} J.~Ribón, {}``Modulus of analytic classification for
unfoldings of resonant diffeomorphisms'', preprint IMPA, 2006.

\bibitem{R1} C.~Rousseau, {}``Normal forms, bifurcations and finiteness
properties of vector fields'', in \emph{Normal forms, bifurcations
and finiteness properties of vector fields}, NATO Sci. Ser. II Math.
Phys. Chem., 137, Kluwer Acad. Publ., Dordrecht, 2004, 431--470.

\bibitem{R2} C.~Rousseau, {}``Modulus of orbital analytic classification
for a family unfolding a saddle-node'', \emph{Moscow Mathematical
Journal} \textbf{5} (2005). 243--267.

\bibitem{RC}C.~Rousseau and C.~Christopher, {}``Modulus of analytic
classification for the generic unfolding of a codimension one resonant
diffeomorphism or resonant saddle'', Annales de l'Institut Fourier
(2007).

\bibitem{T1}L.~Teyssier, {}``Équation homologique et cycles asymptotiques
d'une singularité n\oe ud-col'', \emph{Bulletin des Sciences Math}é\emph{matiques},
vol. 128, \textbf{3} (2004), 167--187.

\bibitem{T2}L.~Teyssier, {}``Analytical classification of singular
saddle-node vector fields'', \emph{Journal of Dynamical and Control
Systems}, vol. 10, \textbf{4} (2004), 577--605.

\bibitem[17]{T4}L.~Teyssier, {}``Équation homologique et classification
analytique des germes de champs de vecteurs holomorphes de type noeud-col'',
thesis, Université de Rennes, \#2869 (2003).

\bibitem[18]{VG}S.~M.~Voronin and A.~A.~Grintchy, {}``An analytic
classification of saddle resonant singular points of holomorphic vector
fields in the complex plane'', \emph{Journal of Dynamical and Control
Systems}, \textbf{2} (1996), 21--53.

\bibitem[19]{Y}J.-C.~Yoccoz, {}``Théorème de Siegel, nombres de
Bruno et polynômes quadratiques'', \emph{Astérisque}, \textbf{231}
(1995), 3--88.

\end{thebibliography}
\end{document}